\documentclass{amsart}
\usepackage{graphicx}
\usepackage{amssymb}
\usepackage{epstopdf}
\usepackage{nicefrac}
\usepackage{enumerate}
\usepackage{hyperref}
\usepackage{float}
\usepackage{color}
\usepackage{pst-all}
\usepackage{tabularx}

\newtheorem{thm}{THEOREM}[section]

\newtheorem{cor}[thm]{COROLLARY}

\newtheorem{defn}[thm]{DEFINITION}
\newtheorem{ex}[thm]{EXAMPLE}

\newtheorem{lemma}[thm]{LEMMA}
\newtheorem{prob}[thm]{PROBLEM}
\newtheorem{prop}[thm]{PROPOSITION}

\newtheorem{remark}[thm]{REMARK}

\newcommand{\ds}{\displaystyle}

\newcommand{\F}{{\mathcal F}} 
\newcommand{\FfM}{{\mathcal F}_{\fM}} 
\newcommand{\FP}{{\mathcal F}_{\cP}} 

 \newcommand{\bA}{{\bf A}}
 
 \newcommand{\bF}{{\bf F}}
 \newcommand{\bH}{{\bf H}}
 \newcommand{\bK}{{\bf K}}
 
 \newcommand{\cC}{{\mathcal C}}
 \newcommand{\cD}{{\mathcal D}}
\newcommand{\cG}{{\mathcal G}}
\newcommand{\cH}{{\mathcal H}}
\newcommand{\cL}{{\mathcal L}}

\newcommand{\cP}{{\mathcal P}}
\newcommand{\cQ}{{\mathcal Q}}
\newcommand{\cS}{{\mathcal S}}

\newcommand{\e}{{\varepsilon}} 
\newcommand{\G}{\Gamma}
\newcommand{\mR}{{\mathbb R}}
\newcommand{\mS}{{\mathbb S}}
\newcommand{\mT}{{\mathbb T}}

\newcommand{\ovpi}{{\overline{\pi}}}

\newcommand{\oU}{{\overline U}}

\newcommand{\mZ}{{\mathbb Z}}
\newcommand{\whG}{{\widehat{G}}}
\newcommand{\whH}{{\widehat{H}}}
\newcommand{\whE}{{\widehat{E}}}
\newcommand{\whTheta}{{\widehat{\Theta}}}
\newcommand{\whPhi}{{\widehat{\Phi}}}

\newcommand{\fM}{{\mathfrak{M}}}

\newcommand{\fG}{{\mathfrak{G}}}
\newcommand{\fT}{{\mathfrak{T}}}
\newcommand{\fX}{{\mathfrak{X}}}

\newcommand{\wtzM}{\widetilde{M_0}}
 \newcommand{\oG}{{\overline{\Phi(G)}}}

\newcommand{\vp}{{\varphi}}

\newcommand{\cU}{{\mathcal U}}
\newcommand{\cT}{{\mathcal T}}

\newcommand{\cGF}{\cG_{\F}} 

 \newcommand{\whL}{\widehat{L}}
 \newcommand{\whM}{\widehat{M}}
 \newcommand{\whfM}{\widehat{\fM}}
 \newcommand{\whcS}{\widehat{\cS}}
 \newcommand{\whF}{\widehat{\F}}
 \newcommand{\whpi}{\widehat{\pi}}
 \newcommand{\whiota}{\widehat{\iota}}

\newcommand{\whg}{{\widehat{g}}}
\newcommand{\whp}{{\widehat{p}}}
\newcommand{\whq}{{\widehat{q}}}
\newcommand{\whcP}{{\widehat{\cP}}}
\newcommand{\whPi}{{\widehat{\Pi}}}
\newcommand{\whfX}{{\widehat{\fX}}}

\newcommand{\wtSigma}{{\widetilde{\Sigma}}}
\newcommand{\wtF}{{\widetilde{\F}}}
 
 \newcommand{\psg}{{\rm pseudo}{\star}{\rm group}}

\input xy
\xyoption {all}

\begin{document}

\title{Molino theory for matchbox manifolds}

   \begin{abstract}    
   A matchbox manifold is a foliated space with totally disconnected transversals, and an equicontinuous matchbox manifold is the generalization of   Riemannian foliations for smooth manifolds in this context. In this paper, we develop the Molino theory for  all equicontinuous matchbox manifolds. Our work extends the Molino theory developed in the work of {\'A}lvarez L\'opez and  Moreira Galicia which required the hypothesis that the holonomy actions for these spaces satisfy the strong quasi-analyticity condition. The methods of this paper are based on the authors' previous works on the structure of weak solenoids, and provide many new properties of the Molino theory for the case of totally disconnected transversals, and examples to illustrate these properties. In particular, we show that the Molino space need not be uniquely well-defined, unless the global holonomy dynamical system is stable, a notion defined in this work. We show that examples in the literature for the theory of weak solenoids provide examples for which the strong quasi-analytic condition fails. 
   Of particular interest is a new class of examples of equicontinuous minimal Cantor actions by finitely generated groups, whose construction relies on a result of Lubotzky. These examples have non-trivial Molino sequences, and   other interesting properties.
   
   \end{abstract}

\author{Jessica Dyer}
\author{Steven Hurder}
 \author{Olga Lukina}
 \email{jessica.dyer@tufts.edu, hurder@uic.edu, lukina@uic.edu}
\address{JD: Department of Mathematics, Tufts University, School of Arts and Sciences, 503 Boston Avenue, Bromfield-Pearson, Medford, MA 02155}
\address{SH \& OL: Department of Mathematics, University of Illinois at Chicago, 322 SEO (m/c 249), 851 S. Morgan Street, Chicago, IL 60607-7045}

\thanks{Version date: October 12, 2016; revised December 23, 2016.}



\maketitle



\section{Introduction}\label{sec-intro}

A smooth foliation $\F$ of a connected compact manifold is a smooth decomposition of $M$ into leaves, which are connected submanifolds of $M$ with constant leaf dimension $n$ and codimension $q$, where $m = n+q$ is the dimension of $M$. This structure is defined by a finite covering of $M$  by   coordinate charts whose image is the product space $(-1,1)^n \times (-1, 1)^q \subset \mR^m$, 
such that the leaves are mapped into linear planes of dimension $n$, and the transition functions between charts preserve these planes. The space $(-1,1)^q$ is called the local transverse model for $\F$.
A smooth foliation $\F$ is said to be \emph{Riemannian}, or \emph{bundle-like},  if there exists a Riemannian metric on the normal bundle $Q \to M$ which is invariant under the transverse holonomy transport along the leaves of $\F$. This condition was introduced by Reinhart in \cite{Reinhart1959}, and is a very strong assumption to impose on a foliation.  The Molino theory for Riemannian foliations gives a complete structure theory for the geometry and dynamics of this class of foliations on compact smooth manifolds  \cite{Haefliger1989,MoerdijkMrcunbook2003,Molino1982,Molino1988}.

An \emph{$n$-dimensional foliated space} $\fM$, as   introduced by Moore and Schochet in \cite{MS2006}, is a continuum - a compact connected metrizable space - 
with  a continuous decomposition of $\fM$ into leaves, which are connected manifolds  with constant leaf dimension  $n$.
Moreover, the decomposition  has a local product structure analogous to that for smooth foliations \cite{CandelConlon2000,MS2006};
that is, every point of $\fM$ has an open neighborhood homeomorphic
to the open subset   $(-1,1)^n \subset \mR^n$ times an open subset of a Polish space  $\fX$, which is said to be the \emph{local
transverse model}. Thus, $\fM$ has a foliation denoted by $\FfM$ whose leaves are the  
maximal path-connected components, with respect to the fine topology on
$\fM$ induced by the plaques of the local product structure.  
 
An \emph{equicontinuous foliated space} is the topological analog of a   Riemannian foliation. In this case, the transverse holonomy pseudogroup associated to the foliation is assumed to act via an equicontinuous collection of local homeomorphisms on the transverse model spaces. The transverse holonomy maps are not assumed to be differentiable,    so there is no natural normal bundle associated to a foliated space, and the standard methods for showing an analog of the Molino theory do not apply. In a series of papers,  {\'A}lvarez L{\'o}pez and Candel \cite{ALC2009,ALC2010}, and {\'A}lvarez L{\'o}pez and Moreira Galicia \cite{ALM2016} formulated a \emph{topological Molino theory} for equicontinuous foliated spaces, which is a partial generalization of  the Molino theory for smooth Riemannian foliations.   They formulated the notion of  \emph{strongly quasi-analytic} ``regularity'' for a foliated space, which is a    condition on the pseudogroup associated to the foliation, as discussed in  Section~\ref{sec-analytic}.  The topological Molino theory in \cite{ALM2016} applies to foliated spaces which satisfy the strongly quasi-analytic condition.

The topological Molino theory for an equicontinuous foliated space $\fM$ with \emph{connected}   transversals   essentially reduces to the smooth theory, by the results of \cite{ALC2010,ALM2016,AB2016}. In contrast, when the transversals to $\FfM$ are \emph{totally disconnected}, and we then say that $\fM$ is a matchbox manifold, the development of a Molino theory in \cite{ALM2016} does not address several key issues, which can be seen as the result of using techniques developed for the smooth theory in the context of totally disconnected spaces. In this work, we apply a completely different approach to developing a topological Molino theory for the case of  totally disconnected transversals. The techniques we use were developed in the authors'   works \cite{Dyer2015,DHL2016a,DHL2016b}. They  are used  here  to  develop a topological Molino theory for     matchbox manifolds in full generality, and to  reveal the far greater complexity of the theory in this case. In particular, we show by our results and examples that the classification of equicontinuous matchbox manifolds via Molino   theory is far from complete.

We recall in Section~\ref{sec-foliatedspaces} the definitions of  a foliated space $\fM$, and of a \emph{matchbox manifold}, which is a foliated space   whose 
local transverse models for the foliation $\F_{\fM}$ are totally disconnected. The terminology ``matchbox manifold'' follows  the usage introduced in continua theory \cite{AO1991,AO1995,AM1988}.   A matchbox manifold with $2$-dimensional leaves is a lamination by surfaces, as defined in \cite{Ghys1999,LM1997}.  If all leaves of $\fM$ are dense, then it is called a \emph{minimal matchbox manifold}.  A compact minimal set $\fM \subset M$ for a foliation $\F$ on a   manifold $M$ yields a foliated space with foliation $\F_{\fM} = \F | \fM$. If the minimal set is exceptional, then $\fM$ is   a minimal matchbox manifold. It is an open problem  to determine which minimal matchbox manifolds are homeomorphic to exceptional minimal sets of   $C^r$-foliations of   compact smooth manifolds, for $r \geq 1$. For example,  the issues associated with this problem are discussed in  \cite{Cass1985,ClarkHurder2011,Hurder2016}.

It was shown in \cite[Theorem~4.12]{ClarkHurder2013} that an equicontinuous    matchbox manifold  $\fM$ is minimal; that is, every leaf is dense in $\fM$. This result generalized a result of  Auslander \cite{Auslander1988} for equicontinuous group actions.   
Examples of equicontinuous matchbox manifolds are given by \emph{weak solenoids}, 
which are discussed in Section~\ref{sec-solenoids}.  Briefly, a weak solenoid $\cS_\cP$ is the inverse limit of a sequence of covering maps $\cP = \{p_{\ell +1} \colon M_{\ell +1} \to M_{\ell} \mid \ell \geq 0\}$, called a \emph{presentation} for $\cS_\cP$, where $M_\ell$ is a compact connected manifold without boundary, and $p_{\ell+1}$ is a finite-to-one covering space.    
The results of  \cite{ClarkHurder2013}  reduce  the study of equicontinuous matchbox manifolds   to the study of weak solenoids:

\begin{thm}  \cite[Theorem~1.4]{ClarkHurder2013}\label{thm-equicontinuous}
An equicontinuous    matchbox manifold $\fM$ is   homeomorphic to a weak solenoid.  
\end{thm}
The idea of the proof of this result is to choose a clopen transversal $V_0 \subset \fM$, then associated to the induced holonomy action of $\FfM$ on $V_0$, 
 one defines (see Proposition~\ref{prop-AFpres}) a  chain of subgroups of finite index,
$\ds \cG =  \{G_0 \supset G_1  \supset \cdots \}$, 
where $G_0$ is the fundamental group of the first shape approximation $M_0$ to $\fM$, where $M_0$ is a compact manifold without boundary. Then $\fM$ is shown to be homeomorphic to the inverse limit of the infinite chain of coverings of $M_0$ associated to the subgroup chain $\cG$.
 
 The theory of inverse limits for covering spaces, as developed for example in \cite{FO2002,McCord1965,Rogers1970,RT1971a,RT1971b,Schori1966}, reduces many questions about the classification of weak solenoids    to   properties of the group chain   $\cG$ associated with the presentation $\cP$.  
Thus, every equicontinuous matchbox manifold $\fM$ admits a presentation which determines its homeomorphism type. 
In  Section~\ref{subsec-weaksols}, the notion of a weak solenoid $\cS_{\cP}$ with presentation $\cP$  is recalled, and the notion of a dynamical partition of the transversal space $V_0$ is introduced in Section~\ref{subsec-partitions}.
  As discussed in Section~\ref{subsec-homeos}, the homeomorphism constructed in the proof of Theorem~\ref{thm-equicontinuous} is well-defined up to return equivalence for the action of the respective holonomy pseudogroups \cite[Section~4]{CHL2013c}.  Thus, we  are interested in   invariants for group chains that are independent of the choice of the chain, up to the corresponding notion of return equivalence for group chains. This is the  approach   we use in this work to formulate and study ``Molino theory'' for weak solenoids.

 Section~\ref{sec-inverselimitmodel} introduces the group chain model for the holonomy action of weak solenoids, following the approach in \cite{Dyer2015,DHL2016a,DHL2016b}. Section~\ref{sec-homogeneous} then recalls results in the literature about homogeneous matchbox manifolds and the associated group chain models for their holonomy actions, which are fundamental for developing the notion of a ``Molino space''.   Section~\ref{sec-ellischains}   introduces the notion of the Ellis group associated to the holonomy action of a weak solenoid.    Ellis semigroups were developed in the works   \cite{Auslander1988,EllisGottschalk1960,Ellis1960,Ellis1969,Ellis2014}, 
  and also appeared in the work  \cite{ALC2010}. A key point of our approach is to use this concept as the foundation of our development of a topological Molino theory.

A key aspect of the Molino space for a foliation is that it is \emph{foliated homogeneous}. A continuum $\fM$ is said to be \emph{homogeneous} if given any pair of points $x,y \in \fM$, then there exists a homeomorphism $h \colon \fM \to \fM$ such that $h(x) = y$. 
A homeomorphism   $\varphi \colon \fM \to \fM$    preserves the path-connected components, hence  a homeomorphism of a matchbox manifold preserves the foliation $\FfM$ of $\fM$. It follows that if $\fM$ is   homogeneous, then it is also foliated homogeneous. 
Our first result is that every equicontinuous    matchbox manifold admits a foliated homogeneous ``Molino space''.
  \begin{thm}\label{thm-molino}
Let $\fM$ be an equicontinuous matchbox manifold, and let $\cP$ be a presentation of $\fM$, such that $\fM$ is homeomorphic to a solenoid $\cS_\cP$. Then there exists a homogeneous matchbox manifold $\whfM$ with foliation $\whF$, called a ``Molino space'' of $\fM$, and a  compact totally disconnected  group $\cD$ (the discriminant group for $\cP$ as defined in Section~\ref{subsec-disc}) such that there exists a  fibration
 \begin{equation}\label{eq-molinoseqM}
\cD  \longrightarrow \whfM \stackrel{\whq}{\longrightarrow} \fM \ , 
\end{equation}
where   the restriction of $\whq$ to each leaf in $\whfM$ is a covering map of some leaf in $\fM$.   We say that \eqref{eq-molinoseqM} is a \emph{Molino sequence} for $\fM$.
 \end{thm}
 
 The construction of the spaces in \eqref{eq-molinoseqM} is given in Section~\ref{sec-molino}.   The homeomorphism type of the fibration \eqref{eq-molinoseqM} depends on the choice of a homeomorphism of $\fM$ with a  weak solenoid $\cS_{\cP}$, and this in turn depends on the choice of the presentation $\cP$ associated to $\fM$ and a section $V_0 \subset \fM$, as discussed   in Section~\ref{subsec-homeos}.   Examples show that the topological isomorphism type of $\cD$ may depend on the choice of the section $V_0$, and the sequence \eqref{eq-molinoseqM} need not be an invariant of the homeomorphism type of $\fM$. This motivates the introduction of the following definition. 
  
 \begin{defn}\label{defn-stablesolenoid}
A matchbox manifold $\fM$ is said to be \emph{stable} if the topological type of the sequence \eqref{eq-molinoseqM} is well-defined  by choosing a sufficiently small transversal $V_0$ to the foliation $\FfM$ of $\fM$.  A matchbox manifold $\fM$ is said to be \emph{wild}, if it is not stable.
  \end{defn}

In Section~\ref{subsec-stableMM} we discuss the relation between the above definition, and the notion  of a stable group chain as given in Definition~\ref{def-stableGC}.   Our next   result concerns the existence of stable matchbox manifolds.

\begin{prop}\label{prop-stablesolenoids}
Let $\fM$ be an equicontinuous matchbox manifold, and suppose $\fM$ admits a transverse section $V_0$ with presentation $\cP$, such that the group $\cD$ in the Molino sequence \eqref{eq-molinoseqM} is finite. Then $\fM$ is stable.
\end{prop}

Proposition~\ref{prop-stablesolenoids} is proved in Section~\ref{sec-molino}. Theorem~\ref{thm-stableCantor} shows that every separable Cantor group $\cD$ can be realized as the discriminant of a stable weak solenoid, but we do not know of a general criteria for when a weak solenoid whose discriminant is a Cantor group must be stable.

The Molino space $\whfM$ is always a homogeneous matchbox manifold. By the results in \cite{DHL2016a},  $\fM$ is homogeneous if and only if for some section $V_0$, the fibration \eqref{eq-molinoseqM} has trivial fibre $\cD$. Each leaf of a homogeneous foliated space   has trivial germinal holonomy, and thus the properties of holonomy for a matchbox manifold $\fM$ are closely related to its non-homogeneity. 
 Section~\ref{sec-holonomy} considers the germinal holonomy groups associated to the global holonomy action for a matchbox manifold. 
 
 Of special importance is the notion of  \emph{locally trivial germinal holonomy}, introduced by Sacksteder and Schwartz \cite{SackstederSchwartz1965}, and used in the work by Inaba  \cite{Inaba1977,Inaba1983} in his  study of Reeb stability for non-compact leaves in smooth foliations.  A leaf $L_x$ in a matchbox manifold $\fM$, which intersects a transversal section $V_0$ at a point $x$, has locally trivial germinal holonomy, if there is an open neighborhood $U \subset V_0$ of $x$, such that the holonomy pseudogroup acts trivially on $U$. A leaf with locally trivial germinal holonomy has trivial germinal holonomy, but the converse need not be true.   In particular,   we prove the following result in Section~\ref{sec-holonomy}. We say that a leaf $L_x$ has   finite $\pi_1$-type if it's fundamental group is finitely generated. A matchbox manifold $\fM$ has finite $\pi_1$-type if all leaves in the foliation $\FfM$ have finite $\pi_1$-type.
 
 \begin{lemma}\label{lemma-germinalholonomy}
 Let $\fM$ be an equicontinuous matchbox manifold with finite $\pi_1$-type. Let $L_x$ be a leaf with trivial germinal holonomy. Then $L_x$ has locally trivial germinal holonomy.
 \end{lemma}
 
 The statement of Lemma~\ref{lemma-germinalholonomy} is implicit in the authors' work \cite{DHL2016b}. 
 The notion of locally trivial germinal holonomy, and  the germinal holonomy properties of equicontinuous matchbox manifolds, turns out to be important in the study of topological Molino theory. Since a weak solenoid is a foliated space, by a fundamental result of Epstein, Millet and Tischler \cite{EMT1977} it contains leaves with trivial germinal holonomy.  A \emph{Schori solenoid} is an example of a weak solenoid, and was first constructed in \cite{Schori1966}. Each leaf in the foliation  of a Schori solenoid   is a surface of infinite genus. 
  
 \begin{prop}\label{Schori-nonloctrivial}
 The Schori solenoid  contains leaves which have trivial germinal holonomy, but do not have locally trivial germinal holonomy.
 \end{prop}
 
 Proposition~\ref{Schori-nonloctrivial} is proved in Section~\ref{sec-analytic}. Proposition~\ref{Schori-nonloctrivial} shows that the condition of finite generation of the fundamental group is essential for the conclusion of Lemma~\ref{lemma-germinalholonomy}. Another result, proved  in Section~\ref{sec-holonomy}, relates the existence of leaves with non-trivial holonomy with non-triviality of the fibre $\cD$ in the Molino sequence \eqref{eq-molinoseqM}.
 
 \begin{thm}\label{thm-nontrivkernel}
  Let $\fM$ be an equicontinuous matchbox manifold. If $\fM$ has a leaf with non-trivial holonomy, then the Molino sequence \eqref{eq-molinoseqM} is non-trivial for any choice of section $V_0 \subset \fM$.
 \end{thm}
 
 The example in Fokkink and Oversteegen \cite{FO2002} and new examples in Section~\ref{sec-construction} 
 show that non-trivial holonomy is not a necessary condition for \eqref{eq-molinoseqM} to be  non-trivial, as one can construct non-homogeneous equicontinuous matchbox manifolds with simply connected leaves.

{\'A}lvarez L\'opez and Moreira Galicia \cite{ALM2016}  investigated Molino theory in the case when the closure of the pseudogroup of an equicontinuous foliated space (in the compact-open topology) satisfies the condition of strong quasi-analyticity (SQA). Geometrically, this means that the pseudogroup action is \emph{locally determined}, that is, if a holonomy map acts trivially on an open subset of it's domain, then it is trivial everywhere on it's domain. A natural problem is to determine which classes of equicontinuous matchbox manifolds are SQA.
This question is studied in Section~\ref{sec-analytic}. 

Note that for equicontinuous actions on Cantor sets the compact-open topology, the uniform topology and the topology of pointwise convergence coincide.
The following result is proved in Section~\ref{sec-analytic}. The set $V_n$ in the statement below  is a   partition set of $V_0 \subset \fT$ as defined in   Proposition~\ref{prop-AFpres}.

\begin{thm}\label{thm-sqaactions}
Let $\fM$ be an equicontinuous matchbox manifold which has finite $\pi_1$-type. Then there exists a transverse section $V_0$, such that the action of the holonomy pseudogroup on this section is SQA. In addition, if $V_0$ can be chosen so that the fibre $\cD$ in the Molino sequence \eqref{eq-molinoseqM} is finite, then there exists a section $V_n \subset V_0$ such that the \emph{closure} of the pseudogroup action on $V_n$ is SQA as well.
\end{thm}
 
On the other hand, there are equicontinuous matchbox manifolds which do not satisfy SQA condition.

\begin{thm}\label{thm-schorinotsqa}
 For every transverse section $V_0$ in the Schori solenoid, the holonomy pseudogroup associated to the section is not SQA.
\end{thm}

Theorem~\ref{thm-molino} proves that the Molino space exists for any matchbox manifold $\fM$, including those who do not admit section with SQA holonomy pseudogroup. Thus, for equicontinuous matchbox manifolds, our results are more general than in \cite{ALM2016}.

Analyzing the results of Lemma~\ref{lemma-germinalholonomy} and Theorem~\ref{thm-sqaactions}, one concludes that the condition of finite $\pi_1$-type, imposed on a matchbox manifold $\fM$, and the condition of finiteness of the fibre $\cD$ in the Molino sequence \eqref{eq-molinoseqM}, are quite strong and force the holonomy pseudogroup to possess various nice properties, such as locally trivial germinal holonomy and the SQA condition. 

It is   natural to ask  how diverse is the class of examples with finite fibre $\cD$ in the Molino sequence?
The authors' work \cite{DHL2016a} constructed new examples of equicontinuous matchbox manifolds with finite fibre $\cD$, which are weakly normal, that is, restricting to a smaller transverse section one can arrange that the Molino sequence \eqref{eq-molinoseqM} has  a trivial fibre. One of these examples is also described in Example~\ref{ex-product} in this paper. Rogers and Tollefson \cite{RT1972} constructed an example of a weak solenoid which turns out to be stable and have finite fibre $\cD$, where the non-triviality of $\cD$ is due to the presence of a leaf with non-trivial holonomy. This example illustrates Proposition~\ref{prop-stablesolenoids} and Theorem~\ref{thm-nontrivkernel}.

 The concluding  Section~\ref{sec-construction} gives the construction of a variety of new classes of examples which illustrate the concepts and results of this work. We first give in Section~\ref{subsec-lenstra} a reformulation of the constructions of the discriminant groups in Section~\ref{sec-ellischains}, in terms of closed subgroups of inverse limit groups, and is analogous  with  a construction attributed to Lenstra in \cite{FO2002}.  This alternate formulation is of strong interest in itself, as it gives a deeper understanding  of the Molino spaces introduced in this work.  This construction can be applied to the examples constructed by   Lubotzky in \cite{Lubotzky1993} showing  the existence of various products of  torsion groups in the profinite completion of torsion-free groups, as recalled in  Section~\ref{subsec-lubotzky}.  We then give three applications of   these results,  which are included in Section~\ref{subsec-stableexamples}. The first construction is based on the conclusions of Theorem~\ref{thm-lubotzky1}.

\begin{thm}\label{thm-stablefinitefibre}
Fix an integer $n \geq 3$. Then there exists a  finite index, torsion-free subgroup $G \subset {\bf SL}_n(\mZ)$ of the $n \times n$ integer matrices, such that given any finite group $F$ of cardinality $|F|$ which satisfies $4(|F| + 2) \leq n$,   there exists an  irregular  group chain $\cG_F$ in $G$ with the properties:
\begin{enumerate}
\item The discriminant group of $\cG_F$ is isomorphic to $F$;
\item The group chain $\cG_F$ is stable, with constant discriminant group isomorphic to $F$;
\item The kernel $K(\cG_F^{\whg})$ of each conjugate $\cG_F^{\whg}$ of this group chain is trivial. 
\end{enumerate}
\end{thm}
The terminology used in Theorem~\ref{thm-stablefinitefibre} will be explained in later sections, where we will show that given such a group chain, one can construct matchbox manifolds with the following properties:

 \begin{cor}\label{cor-stablefinitefibre}
Let $F$ be a finite group. Then  there exists a non-homogeneous  matchbox manifold $\fM$ such that every leaf of $\FfM$ has trivial germinal holonomy, and 
  for any sufficiently small   transverse section in $\fM$, its Molino sequence   is non-trivial with fiber group $\cD \cong   F$.
 \end{cor}
 Note that it follows by Theorem~\ref{thm-sqaactions} that   for the examples constructed in the proof of Corollary~\ref{cor-stablefinitefibre}, there is a 
 section $V  \subset \fM$ such that the closure of the pseudogroup action on $V$ satisfies the SQA condition of {\'A}lvarez L\'opez and Moreira Galicia \cite{ALM2016}.

   The next two constructions are based on the conclusions of Theorem~\ref{thm-lubotzky2} of Lubotzky. Again, the terminology used in the statements will be explained in later sections.

\begin{thm}\label{thm-stableCantorfibre}
There exists a  finite index, torsion-free finitely-generated group $G$  such that given any   separable profinite group $K$,   there exists an irregular   group chain $\cG_K$ in $G$ such that:
\begin{enumerate}
\item The discriminant group of $\cG_K$ is isomorphic to $K$;
\item The group chain $\cG_F$ is stable, with constant discriminant group isomorphic to $K$. 
\end{enumerate}
\end{thm}
 \begin{cor}\label{cor-stableCantorfibre}
Let $K$ be a Cantor group. Then  there exists a non-homogeneous  matchbox manifold $\fM$ such that,  
  for any sufficiently small   transverse section in $\fM$, its Molino sequence   is non-trivial with fiber group $\cD \cong   K$.  
 \end{cor}

 Finally, Theorem~\ref{thm-nvr} gives the first   examples of equicontinuous matchbox manifolds which are not virtually regular. The \emph{virtually regular} condition was introduced in the work  \cite{DHL2016b}, and is defined in Definition~\ref{def-vitreg}. As the terminology suggests, this notion is related to the homogeneity  properties of finite-to-one coverings of a matchbox manifold $\fM$. 
 
The concluding   Section~\ref{subsec-open} lists some open problems.

   \section{Equicontinuous Cantor foliated spaces}\label{sec-foliatedspaces}

  In this section, we recall   background concepts about foliated spaces,  and introduce the group chains associated to  their  equicontinuous Cantor holonomy actions.  
  
\subsection{Equicontinuous Cantor foliated spaces}
Recall that an \emph{$n$-dimensional matchbox manifold} $\fM$ is a compact connected metrizable topological space such that every point $x \in \mathfrak{M}$ has an open neighborhood $U \subset \fM$ such that there is a homeomorphism
\begin{equation}\label{eq-coordinatechart}
\varphi_x \colon  \oU_x \to [-1,1]^n \times \fT_x \ ,
\end{equation}
where $\fT_x$ is a totally disconnected space. The homeomorphism $\varphi_x$ is called a \emph{local foliation chart}, and the space $\fT_x$ is called a \emph{local transverse model}.  As usual in foliation theory, one can choose a finite atlas $\cU = \{(\varphi_i,U_i)\}_{1 \leq i \leq \nu}$ of local charts, such that the intersections of the path-connected components in $U_x \cap U_y$ are connected and simply connected, and the images $\cT_i =\varphi_i^{-1}( \{0\} \times \fT_i)$ are disjoint. The leaves of the foliation $\FfM$ of $\fM$ are defined to be the path connected components of $\fM$, which are then a union of the path connected components (the plaques) in the open sets $U_i$.

A matchbox manifold is \emph{(topologically) minimal} if   each leaf $L \subset \fM$ is dense in $\fM$.

We   require the matchbox manifold $\fM$ to be \emph{smooth}; that is, the transition maps
   $$\varphi_j \circ \varphi_i^{-1}  \colon  \varphi_i^{-1}(\overline{U}_i \cap \overline{U}_j) \to \varphi_j(\overline{U}_i \cap \overline{U}_j)$$
are $C^{\infty}$-maps in the first coordinate $x \in [-1,1]^n$, and the restrictions to plaques depend continuously  on $y \in   \fT_i$, in the $C^\infty$-topology on leaves,  for $1 \leq i,j \leq \nu$.   

 Let $pr_2 \colon [-1,1]^n \times \fT_i \to  \fT_i$ be the projection onto the second factor, then    $\pi_i = pr_2 \circ \varphi_i \colon  \oU_i \to \fT_i$ for $1 \leq i \leq \nu$, are the  local defining maps for the foliation $\FfM$.
 Denote by $\fT_{i,j} = \pi_i(U_i \cap U_j)$ for    $1 \leq i , j \leq \nu$. 
Since the path-connected components of the charts are either disjoint, or have a connected intersection, there is a well-defined change-of-coordinates homeomorphism
  \begin{align}\label{eq-transitionmaps} 
  h_{i,j} = \pi_j \circ \pi_i^{-1} \colon  \fT_{i,j} \to \fT_{j,i} 
  \end{align}
 with domain $\fT_{i,j}$ and range $\fT_{j,i}$. Let $\cGF^1 = \{(h_{i,j},\fT_{i,j}) \mid 1\leq i,j \leq \nu\}$.
   Denote $\fT = \fT_1 \cup \cdots \cup \fT_\nu$. Then the collection of maps $\cGF^1$ generates the \emph{holonomy pseudogroup} $\cGF$ 
   acting on the transverse space  $\fT$. 
   The construction and properties of $\cG_{\F}$ is described in full detail in \cite[Section~3]{ClarkHurder2013}.
   
   For the study of the dynamical properties of $\FfM$, it is useful to introduce also the collection of maps $\cGF^* \subset \cGF$, defined as follows.
   Let $\cG_0 \subset \cGF$ denote the collection  consisting of all possible compositions of homeomorphisms in $\cGF^1$. Then $\cGF^*$   consists of all possible restrictions of homeomorphisms in $\cG_0$ to open subsets of their domains.  
  The collection of maps $\cGF^*$    is closed under the operations of compositions, taking inverses, and restrictions to open sets, and  is called a \emph{$\psg$} in the works \cite{ALM2016,Matsumoto2010}, while   $\cGF^*$ is called a \emph{localization} of $\cG_0$ in the work \cite{ALM2016}.

\begin{remark}\label{rmk-pseudostar}
 {\rm 
 The standard definition of a pseudogroup \cite{CandelConlon2000} requires the pseudogroup to be closed under the operations of composition, taking inverses,  restriction to open subsets, and of combination of maps. A combination of two local homeomorphisms $h_1$ and $h_2$, with possibly disjoint domains $D(h_1)$ and $D(h_2)$ and with disjoint ranges, is a homeomorphism $h$ defined on $D(h_1) \cup D(h_2)$ where $h|D(h_1) = h_1$ and $h|D(h_2) = h_2$. However, allowing such arbitrary glueings of maps is unnatural. For example, a composition $h_{j,k} \circ h_{i,j}$ can be associated with the existence of a leafwise path $\gamma_x \colon [0,1] \to L_x \in \fM$ with $\gamma_x(0) \in U_i$ and $\gamma_y(1) \in U_k$, where $L_x$ is a leaf such that $\pi_i(x) \in D(h_{j,k} \circ h_{i,j})$. If $\pi_i(y) \in D(h_{j,k} \circ h_{i,j})$, then the path $\gamma_x$ can be lifted to a nearby leaf $L_y$ to a `parallel' path $\gamma_y$ with $\gamma_y(0) \in U_i$ and $\gamma_y(1) \in U_k$. Thus a holonomy transformation $h_{j,k} \circ h_{i,j}$ has a geometric meaning as the transverse transport in leaves along a leafwise path. 
 Therefore, in the definition of   $\cG_0$ and $\cGF^*$ (and of a $\psg$ in \cite{Matsumoto2010}), one does not allow combinations of local homeomorphisms, unless such homeomorphisms can be obtained by restrictions to open subsets of maximal domains of elements in $\cG_0$. 
  }
\end{remark}

Let $d_\fM$ be a metric on $\fM$, and denote by $d_{\cT_i}$ the restriction of $d_\fM$ to the embedded image $\cT_i$ of the transversal $\fT_i$, $1 \leq i \leq \nu$. For each $1 \leq i \leq \nu$, consider the pullback $d_{\fT_i}$ of $d_{\cT_i}$ along the embedding. Then define a metric $d_\fT$ on $\fT$ by the following formula:
  \begin{align*} d_\fT(x,y) = \left\{ \begin{array}{ll} d_{\fT_i}(x,y), & \textrm{if }x,y \in \fT_i \textrm{ for some }i, \\ \infty, & \textrm{otherwise}. \end{array} \right. \end{align*}
For a homeomorphism $\gamma \in \cGF^*$, denote by $D(\gamma)$ and $R(\gamma)$ the domain and the range of $\gamma$ respectively.

\begin{defn}\label{defn-equicontinuity}
The action of the $\psg$ $\cGF^*$ on the transversal $\fT$ is \emph{equicontinuous} if for all $\e>0$ there exists $\delta>0$ such that for all $\gamma \in \cGF^*$, if $x,x' \in D(\gamma)$ and $d_\fT(x,x') < \delta$, then $d_\fT(\gamma(x),\gamma(x')) < \e$.
\end{defn}

The following notion is used in the statement of various results in this work.
\begin{defn}\label{def-topfinitetype}
A path-connected topological space $X$ is said to have finite $\pi_1$-type if the fundamental group $\pi_1(X , x)$ is a finitely-generated group, for the choice of some basepoint $x \in X$. A matchbox manifold $\fM$ is said to have finite $\pi_1$-type if each leaf $L \subset \fM$ is a space of finite $\pi_1$-type.
\end{defn}
  
 \subsection{Suspensions}\label{subsec-suspensions}
 
 There is a well-known construction which yields a foliated space  from a group action, called the \emph{suspension construction}, as discussed in   \cite[Chapter 3]{CandelConlon2000} for example. We state this construction in the restricted context which we use in this work.
  
Let $X$ be a Cantor space, and $H$ a finitely-generated group, and assume there is given  an action $\vp \colon H \to Homeo(X)$. Suppose that $H$ admits a generating set $\{g_1 , \ldots , g_k\}$, then there is a   homomorphism $\alpha_k \colon \mZ * \cdots * \mZ \twoheadrightarrow H$ of the free group on $k$ generators onto $H$, given by mapping generators to generators. Of course, the map $\alpha_k$ will have non-trivial kernel, unless $H$ happens to be a free group. Next, let $\Sigma_k$ be a compact surface without boundary of genus $k$. Then for a choice of basepoint $x_0 \in \Sigma_k$ set   $G = \pi_1(\Sigma_k, x_0)$. Then there is a   homomorphism $\beta_k  \colon G \to \mZ * \cdots * \mZ$ onto the free group of $k$ generators.   
Denote the composition of these   maps by $\Phi = \vp \circ \alpha_k \circ \beta_k$ to obtain the   homomorphism 
$\ds \Phi \colon G = \pi_1(\Sigma_k, x_0) \to  \mZ * \cdots * \mZ \to H \to  Homeo(X)$.  

Now, let $\wtSigma_k$ denote the universal covering space of $\Sigma_k$, equipped with the right action of $G$ by covering transformations. Form the product space $\wtSigma_k \times X$ which has a foliation $\wtF$ whose leaves are   the slices $\wtSigma_k \times \{x\}$ for each $x \in X$. Define a left action of $G$ on $\wtSigma_k \times X$, which for $g \in G$ is given  by $(y,x) \cdot g = (y \cdot g, \Phi(g^{-1})(x))$. For each $g$,  this action preserves the foliation $\wtF$, so we obtain a foliation $\FfM$ on the quotient space $\fM = (\wtSigma_k \times X)/G$.  Note that all leaves of $\FfM$ are surfaces, which are in general non-compact.

Note that $\fM$ is a foliated Cantor bundle over $\Sigma_k$, and the holonomy of this bundle $\pi \colon \fM \to \Sigma_k$ acting on the fiber $V_0 = \pi^{-1}(x_0)$ is canonically identified with the action $\Phi \colon G \to Homeo(X)$. Consequently,  if the action $\Phi$ is minimal in the sense of topological dynamics \cite{Auslander1988}, then the foliation $\FfM$ is minimal. If the action $\Phi$ is equicontinuous in the sense of topological dynamics \cite{Auslander1988}, then $\FfM$ is an equicontinuous foliation in the sense of Definition~\ref{defn-equicontinuity}.
 
 There is a variation of the above construction, where we assume that $G$ is a   \emph{finitely-presented} group, and there is given a homomorphism   $\Phi \colon G \to Homeo(X)$. In this case, it is a well-known folklore result (for example, see \cite{Massey1991}) that there   exists a closed connected $4$-manifold $B$ such that for a choice of basepoint $b_0 \in B$, then $\pi_1(B, b_0)$ is homeomorphic to $G$.  Then the suspension construction can be applied to the homomorphism $\Phi \colon \pi_1(B, b_0) \to Homeo(X)$, where we replace $\Sigma_k$ above with $B$,  and the space $\wtSigma_k$ with the universal covering $\widetilde{B}$ of $B$. The resulting foliated space $\fM$ will have holonomy given by the map $\Phi$. 
 
 In summary, the suspension construction translates results about equicontinuous minimal Cantor actions to  results about equicontinuous matchbox manifolds.

\section{Weak solenoids}\label{sec-solenoids}

In this section, we first recall the construction procedure for \emph{(weak) solenoids}, and describe some of their properties.
In Section~\ref{subsec-partitions}, we discuss the construction from \cite{ClarkHurder2013} which associates a group chain to an equicontinuous  matchbox manifold, which leads to a more precise statement of Theorem~\ref{thm-equicontinuous}. Then in Section~\ref{subsec-homeos},   we make some observations about the conclusion of Theorem~\ref{thm-equicontinuous}  which are important when considering the definition of the Molino space for matchbox manifolds.

 \subsection{Weak solenoids}\label{subsec-weaksols}
  
Let $n \geq 1$, then  for each $\ell \geq 0$, let $M_{\ell}$ be a   compact connected simplicial complex of dimension $n$. A  \emph{presentation}   is a collection $\cP = \{ p_{\ell+1} \colon M_{\ell+1} \to M_{\ell} \mid \ell \geq 0\}$, where each map $p_{\ell +1}$  is a proper surjective map of   simplicial complexes with discrete fibers, which     is called a \emph{bonding} map. For   $\ell \geq 0$ and $x \in M_{\ell}$,  the preimage $\{p_{\ell +1}^{-1}(x) \} \subset M_{\ell +1}$  is compact and discrete, so the cardinality $\#\{p_{\ell +1}^{-1}(x) \}<\infty$. For a presentation  $\cP$  defined in this generality, the cardinality of the fibers  of the maps $p_{\ell +1}$  need not be constant in either  $\ell$ or $x$.  
 
    Associated to a presentation $\cP$ is an inverse limit space,  
\begin{eqnarray} 
\cS_{\cP} & \equiv &  \lim_{\longleftarrow} ~ \{ p_{\ell +1} \colon M_{\ell +1} \to M_{\ell}\} \label{eq-presentationinvlim} \\
& = &  \{(x_0, x_1, \ldots ) \in \cS_{\cP}  \mid p_{\ell +1 }(x_{\ell + 1}) =  x_{\ell} ~ {\rm for ~ all} ~ \ell \geq 0 ~\} ~ \subset \prod_{\ell \geq 0} ~ M_{\ell} ~ . \nonumber
\end{eqnarray}
The set $\cS_{\cP}$ is given  the relative  topology, induced from the product (Tychonoff) topology, so that $\cS_{\cP}$ is itself compact and connected. 

\begin{defn}\label{defn-solenoid}
The inverse limit space $\cS_{\cP}$ in \eqref{eq-presentationinvlim} is called a \emph{(weak) solenoid}, if for  each  $\ell \geq 0$ the space $M_{\ell}$ is a compact connected manifold without boundary, and   $p_{\ell +1}$ is a proper covering map of degree $m_{\ell+1} > 1$. 
\end{defn}

Weak solenoids are a generalization of $1$-dimensional (Vietoris) solenoids, described in Example~\ref{ex-Vietoris} below. Weak solenoids were originally considered  in the papers by McCord \cite{McCord1965}, Rogers and Tollefson \cite{RT1971a,RT1972} and Schori  \cite{Schori1966}, and later by Fokkink and Oversteegen \cite{FO2002}.

\begin{ex}\label{ex-Vietoris}
{\rm 
Let  $M_{\ell} = \mS^1$ for each $\ell \geq 0$, and let the map $p_{\ell +1}$ be a proper covering map of degree $m_{\ell+1} > 1$ for $\ell \geq 0$. Then $\cS_{\cP}$ is an example of a  {classic $1$-dimensional solenoid},  discovered independently by    van~Dantzig \cite{vanDantzig1930} and   Vietoris   \cite{Vietoris1927}. If $m_{\ell+1} = 2$ for $\ell \geq 0$, then $\cS_{\cP}$ is called the \emph{dyadic} solenoid.
}
\end{ex}

Let $\cS_\cP$ be a weak solenoid as in Definition~\ref{defn-solenoid}. For each $\ell \geq 1$, the composition
\begin{equation}\label{eq-coverings}
q_{\ell} = p_{1} \circ \cdots \circ p_{\ell -1} \circ p_{\ell} \colon M_{\ell} \to M_0 ~ 
\end{equation}
is a finite-to-one covering map of the base manifold $M_0$. For each $\ell \geq 0$, projection onto the $\ell$-th factor in the product $\ds \prod_{\ell \geq 0} ~ M_{\ell}$ in \eqref{eq-presentationinvlim} yields a   fibration map denoted by $\Pi_{\ell} \colon \cS_{\cP}  \to M_{\ell}$. For $\ell =0$ this yields the fibration $\Pi_0 \colon \cS_{\cP} \to M_0$, and for $\ell \geq 1$ we have  
 \begin{align}\label{eqn-bundle}
 \Pi_0 =   q_{\ell} \circ \Pi_{\ell}  \colon \cS_{\cP} \to M_0 \ .
 \end{align}
A choice of a basepoint $x_0 \in M_0$ fixes a fiber  $\fX_0 = \Pi_0^{-1}(x_0)$, which is a Cantor set by the assumption  that the fibers of each map $p_{\ell+1}$ have cardinality at least $2$. McCord showed in \cite{McCord1965} that \eqref{eqn-bundle} is a fibre bundle over $M_0$ with a Cantor set fibre, and   the solenoid $\cS_\cP$ has a local product structure as in \eqref{eq-coordinatechart}.
The path-connected components of $\cS_\cP$ thus define a foliation denoted by $\FP$.  We then have:

\begin{prop}\label{prop-solenoidsMM}
Let   $\cS_{\cP}$ be a weak solenoid, whose    base space $M_0$   is a compact manifold of dimension $n \geq 1$. Then   $\cS_{\cP}$ is  a minimal matchbox manifold of dimension $n$ with foliation   $\FP$. 
\end{prop}

Denote by   $G_0 =  \pi_1(M_{0}, x_{0})$ the fundamental group of $M_0$ with basepoint $x_0$, and choose a point $ x \in \fX_0$ in the fibre over $x_0$. This defines basepoints  $x_{\ell} = \Pi_{\ell}(x) \in M_{\ell}$ for $\ell \geq 1$. 

Let $y \in \fX_0$ be another point, and set $y_{\ell} = \Pi_{\ell}(y) \in M_{\ell}$, and note that $y_0 = x_0$ by construction. We will interchangeably write $y = (y_{\ell})$ to denote a point in $\fX_0$ or $\cS_\cP$. Let $L_y$ denote the leaf of $\FP$ containing $y$. 
Then the restriction $\Pi_0|_{L_y} \colon  L_y \to M_0$ of the bundle projection to each path-connected component $L_y$ is a covering map. 
For $g = [\gamma_0] \in G_0$, let $\gamma_{\ell} \colon [0,1] \to M_{\ell}$ be a lift of $\gamma_0$ with the starting point $\gamma_{\ell}(0) = y_{\ell}$. Define a homeomorphism
   $ h_g \colon  \fX_0 \to \fX_0$ by $h_g(y_{\ell}) = (\gamma_{\ell}(1))$. Thus there is a representation
  \begin{align}\label{eqn-globalholonomy} 
  \Phi_0 \colon G_0 \to Homeo(\fX_0) \colon \gamma \to h_g \ , 
  \end{align}
called the \emph{global holonomy map} of the solenoid $\cS_\cP$.

 \subsection{Dynamical partitions}\label{subsec-partitions} 
It was shown in \cite[Theorem~4.12]{ClarkHurder2013} that an equicontinuous    matchbox manifold  $\fM$ is minimal, that is, every leaf is dense in $\fM$. This result  generalizes to pseudogroups a corresponding result of   Joe~Auslander  for equicontinuous group actions in \cite{Auslander1988}.  
  It follows that for any clopen subset $V_0 \subset \fT$, the restricted $\psg$ $\cG_{V_0}^* = \cGF^*|V_0$ is \emph{return equivalent} to the $\psg$ $\cGF^*$ on $\fT$, where return equivalence is defined and studied in   \cite[Section~4]{CHL2013c}. Thus, for the study of the dynamical properties of $\FfM$ one can restrict to the study of $\cG_{V_0}^*$.  The following result   is based on the constructions in \cite{ClarkHurder2013}.

 \begin{prop}\label{prop-AFpres}\cite{ClarkHurder2013}
Let $\fM$ be a matchbox manifold with totally disconnected transversal $\fT$ and equicontinuous holonomy $\psg$ $\cGF^*$ on $\fT$,   let $x \in \fT$ be a point, and let $W \subset \fT$ be a clopen (closed and open) neighborhood of $x$. Then there exists a clopen subset $x \in V_0 \subset W$ and a descending chain of clopen sets $V_0 \supset V_1\supset \cdots $ of $\fT$ with $  \{x\} = \bigcap_{\ell} V_{\ell}$, such that:
\begin{enumerate}
\item The restriction $\cGF^*|V_0$ is generated by a group $G_0$ of transformations of $V_0$.
\item \label{property-1} For each $\ell \geq 1$ the collection $\cQ_{\ell} =\{g \cdot V_{\ell}\}_{g \in G_0}$ is a finite partition of $V_0$ into clopen sets. 
\item \label{property-2} We have ${ diam}(g \cdot V_{\ell} ) < 2^{-\ell}$ for all $g \in G_0$ and all $\ell \geq 0$.
\item The collection of elements which fix $V_{\ell}$, that is,
 $$G_{\ell}^x = \{g \in G_0 \mid g \cdot V_{\ell}= V_{\ell}\},$$
 is a subgroup of finite index in $G_0$. More precisely, $|G_0:G_{\ell}^x| = { card}(\cQ_{\ell})$.
\end{enumerate}
\end{prop}

  There are many  choices involved in the construction of the partitions   $\cQ_{\ell}$ and consequently the stabilizer groups $G_{\ell}^x$:
  \begin{enumerate}
\item The choice of a transverse section $V_0 \subset \fT$, which results in the choice of the group $G_0$.
\item The choice of a basepoint $x \in V_0$.
\item Given $V_0$, $x$ and $G_0$, there is freedom to choose clopen sets $V_1 \supset V_2 \supset \ldots$, which results in the choice of the sequence of groups 
$G_0 = G_0^x \supset G^x_1 \supset G^x_2 \supset \ldots$. 
\end{enumerate}

Thus, the algebraic and geometric data  encoded by these choices must be considered up to suitable notions of equivalence, which will be introduced   in Section~\ref{subsec-chainequiv}.

\subsection{Homeomorphisms}\label{subsec-homeos}

Let $\fM$ be a matchbox manifold with totally disconnected transversal $\fT$ and equicontinuous holonomy $\psg$ $\cGF^*$ acting on $\fT$,   let $x \in \fT$ be a point, and let   $\{V_{\ell +1} \subset V_{\ell} \mid \ell \geq 0\}$ be a descending chain of clopen subsets of $\fT$ with $x \in V_{\ell}$ for all $\ell \geq 0$, as introduced in 
 Proposition~\ref{prop-AFpres}, where $G_0$ is a group of transformations of $V_0$, and  $G_{\ell}$ denotes the stabilizer subgroup of $G_0$  of the set $V_{\ell}$.
 
The basic idea of the proof of  Theorem~\ref{thm-equicontinuous}, is that if we choose the section $V_0 \subset \fM$ sufficiently small and appropriately chosen, then there is a compact manifold $M_0$ and a fibration $\Pi_0' \colon \fM \to M_0$ for which the inverse image $(\Pi_0')^{-1}(x_0) = V_0$ where $x_0 = \Pi_0'(x)$. Moreover, the restriction of the map $\Pi_0'$ to the leaves of $\FfM$ are coverings of $M_0$. The definition of the map $\Pi_0'$ requires the highly technical results of \cite{CHL2013a} to define a transverse Cantor foliation $\cH_0$ to $\FfM$, so that the quotient space $M_0 = \fM/\cH_0$ is a compact manifold, and then $\Pi_0'$ is the projection along the leaves of the transverse foliation $\cH_0$, or better said the equivalence classes defined by the leaves of $\cH_0$. Then $V_0$ is the $\cH_0$-equivalence class of the point $x \in V_0 \subset \fM$.

  Let $V_{\ell} \subset V_0$ be the clopen set in Proposition~\ref{prop-AFpres} and     $\ds G_{\ell}^x = \{g \in G_0 \mid g \cdot V_{\ell}= V_{\ell}\}$   the isotropy subgroup of $V_{\ell}$.  Then there is a Cantor subfoliation $\cH_{\ell}$ of $\cH_0$ such that $V_{\ell}$ is the $\cH_{\ell}$-equivalence class of $x$.  
  Moreover,  there is a quotient map $\Pi_{\ell}' \colon \fM \to \fM/\cH_{\ell} \equiv M_{\ell}$ where is $M_{\ell}$ is identified with the covering of $M_0$ associated to the subgroup $G_{\ell} \subset G_0 = \pi_1(M_0, x_0)$. Note that the fiber   $(\Pi_{\ell}')^{-1}(x_0) = V_{\ell}$ and the monodromy action of $G_0$ on $V_0$ partitions $V_0$ into the translates of $V_{\ell}$.
  There is then a quotient covering map $q_{\ell} \colon M_{\ell} \to M_0$, and as in \eqref{eqn-bundle}, we have
  \begin{align}\label{eqn-bundle2}
 \Pi_0' =   q_{\ell} \circ \Pi_{\ell}'  \colon \fM \to M_0 \ .
 \end{align}
For each $\ell \geq 0$ let $p_{\ell +1} \colon M_{\ell +1} \to M_{\ell}$ be the quotient map defined by expanding the equivalence classes of $\fM$ defined by $\cH_{\ell +1}$ to the equivalence classes   by $\cH_{\ell}$. Then we obtain a collection of covering maps $\cP = \{ p_{\ell+1} \colon M_{\ell+1} \to M_{\ell} \mid \ell \geq 0\}$ which defines a weak solenoid $\cS_{\cP}$. As the diameters of the clopen partition sets $V_{\ell}$ tends to $0$ as $\ell$ increases, it is then standard that the collection of maps 
$\ds \{\Pi_{\ell}' \colon \fM \to M_{\ell} \mid \ell \geq 0 \}$ induces a foliated homeomorphism $\Pi_0^* \colon \fM \to \cS_{\cP}$.

In later sections, we will also consider the presentations $\cP_n$ obtained  by truncating the initial  $n$ terms   in   the   presentation $\cP$. That is, for $n \geq 0$ we have 
\begin{equation}\label{eq-Pn}
\cP_n = \{p_{\ell +1}' \colon M_{\ell +1}' \to M_{\ell}' \mid \ell \geq 0\} ~ , ~ {\rm where} ~ M_{\ell}' = M_{\ell + n} ~ {\rm and} ~  p_{\ell +1}' = p_{\ell + n+1} \ .
\end{equation}
It is a basic property of inverse limit spaces \cite{McCord1965,Rogers1970} that for $n \geq 1$ and $m \geq  0$, there is a homeomorphism  $\sigma_n \colon \cS_{\cP_{m+n}} \cong \cS_{\cP_n}$, 
 where the homeomorphism is given by the ``shift in coordinates'' map $\sigma_n$ in the inverse sequences defining these spaces. Also, by the same reasoning as above, there is a foliated homeomorphism $\Pi_{n}^* \colon \fM \to \cS_{\cP_n}$ and we have a commutative diagram of fibrations: 
    \begin{align}\label{diag-homeos} 
   \xymatrixcolsep{4pc}
   \xymatrix{
   \fM  \ar[d]_{\Pi_{n+m}^*} \ar[r]^{=} &  \fM  \ar[d]^{\Pi_m^*}  \\
   \cS_{\cP_{n+m}}   \ar[r]^{\sigma_n} & \cS_{\cP_m} 
   }  
   \end{align}
Note that if the presentation $\cP$ is constructed using the holonomy of $\FfM$ acting on the transversal $V_0 \subset \fM$, then for   $n > 0$ and $m \geq 0$, the map $\sigma_n \colon \cS_{\cP_{n+m}} \to \cS_{\cP_{m}}$ satisfies $\sigma_n(V_{m+n}) \subset V_n$. That is, the induced map on $\fM$ sends the   transversal $(\Pi_{m+n}^*)^{-1}(V_{m+n}) \subset \fM$ into the 
  transversal $(\Pi_{n}^*)^{-1}(V_{n}) \subset \fM$. On the other hand,  given a homeomorphism $h \colon \fM \to \fM$ there is no reason it should map the transversal $V_0$ into itself. In particular, the induced map 
\begin{equation}\label{eq-fiberpreserving}
(\Pi_n^*) \circ  h \circ  (\Pi_{m+n}^*)^{-1}  \colon \cS_{\cP_{m+n}} \to \cS_{\cP_n} 
\end{equation}
  on weak solenoids need not be fiber preserving. On the other hand, as discussed in \cite{FO2002}, there is always a map $h' \colon \fM \to \fM$ which is homotopic to $h$ such that the induced map as in \eqref{eq-fiberpreserving} maps a clopen subset of $V_{m+n}$ into a clopen subset of $V_n$. Thus, by allowing sufficiently large values of $n$ and $m$ and choice of basepoints in the range and domain, we can always ensure that a given homeomorphism of $\fM$ induces a fiber-preserving map between the weak solenoids $\cS_{\cP_{m+n}}$ and  $\cS_{\cP_n}$.

  \section{Group chain models}\label{sec-inverselimitmodel}

Let $\cS_{\cP}$ be a weak solenoid defined by a presentation $\cP$, with basepoint $x \in  \fX_0 \equiv \Pi_0^{-1}(x_0) \subset \cS_{\cP}$. 
For  $G_0 = \pi_1(M_0, x_0)$, let $\Phi_0 \colon G_0 \to Homeo(\fX_0)$ be the holonomy action in \eqref{eqn-globalholonomy}.

The following ``combinatorial  model''    for the action \eqref{eqn-globalholonomy} allows for a deeper analysis of the relation between the action $\Phi_0$ and the algebraic structure of $G_0$. For each $\ell \geq 1$, recall that  
\begin{equation}\label{eq-imahes}
G^x_{\ell} = {\rm image}\left\{  (q_{\ell} )_{\#} \colon  \pi_1(M_{\ell}, x_{\ell})\longrightarrow G_{0}\right\}  
\end{equation}
  denotes  the image of the induced map $(q_{\ell} )_{\#} $ on fundamental groups. In this way, associated to the presentation $\cP$ and basepoint $x \in \fX_0$,  we obtain a descending chain of subgroups of finite index
  \begin{equation}\label{eq-descendingchain}
\cG^x \colon  G_{0} \supset G^x_{1} \supset G^x_{2} \supset \cdots \supset G^x_{\ell} \supset \cdots   \ .
\end{equation}
Each quotient  $X_{\ell}^x = G_{0}/G_{\ell}^x$ is a finite set equipped with a left $G_0$-action, and there are surjections $X_{\ell +1}^x \to X_{\ell}^x$ which commute with the action of $G_0$.  The inverse limit,  
\begin{equation}\label{eq-Galoisfiber}
X_{\infty}^x = \lim_{\longleftarrow} ~ \{ p_{\ell +1} \colon X_{\ell +1}^x \to X_{\ell}^x\} = \{(eG_0, g_1G^x_1, \ldots) \mid g_{\ell} G^x_{\ell} = g_{\ell+1}G^x_{\ell}\} ~ \subset \prod_{\ell \geq 0} ~ X_{\ell}^x  
\end{equation}
is then a totally disconnected compact perfect set, so is   a Cantor set. The   fundamental group $G_0$ acts on the left on  $X_{\infty}^x$ via the coordinate-wise multiplication on the product in \eqref{eq-Galoisfiber}. We denote this   Cantor action by $(X_{\infty}^x , G_0 , \Phi_x)$.

\begin{lemma}\label{lem-fibremodel}
There is a homeomorphism $\tau_x \colon  \fX_0 \to X^x_\infty$ equivariant with respect to the action \eqref{eqn-globalholonomy} of $G_0$ on $\fX_0$ and $\Phi_x$ on   $X_{\infty}$; that is, $\tau_x \circ h_g(y) = \Phi_x(g) \circ \tau_x(y)$ for all $y \in \fX_0$.
\end{lemma}

In particular, this allows us to conclude that the action $\Phi_0$ of $G_0$ on the fibre of the solenoid $\cS_{\cP}$ is minimal. Indeed, the left action of $G_0$ on each quotient space $X_{\ell}^x$ is transitive, so the orbits are dense in the product topology on $X_{\infty}^x$.

\begin{remark}\label{rmk-basepoints}
{\rm 
The group chain \eqref{eq-Galoisfiber} and the homeomorphism in Lemma~\ref{lem-fibremodel} depend on the choice of a point $x \in \fX_0$. For a different basepoint $y \in \fX_0$ in the fibre over $x_0$, let   $\tau_x(y) = (g_iG_{\ell}^x) \in X_\infty^x$, then the group chain $\cG^y$ associated to $y$ is given by a   chain of conjugate subgroups in $G_0$, where 
 $G_{\ell}^y = g_{\ell} G_{\ell}^x g_{\ell}^{-1}$ for $\ell \geq 0$.
 The group chains $\cG^y$ and $\cG^x$  are said to be \emph{conjugate chains}.
  The composition $\tau_y \circ \tau_x^{-1} \colon X_{\infty}^x \to X_{\infty}^y$ gives a topological conjugacy between the minimal Cantor actions  $(X_{\infty}^x , G_0 , \Phi_x)$ and $(X_{\infty}^y , G_0 , \Phi_y)$.
The map $\tau_x \colon  \fX_0 \to X_{\infty}^x$ can be viewed as ``coordinates'' on the inverse limit space $\fX_0$, and the composition  $\tau_y \circ \tau_x^{-1}$ as a ``change of coordinates''. Properties of the minimal Cantor action 
$(X_{\infty}^x , G_0 , \Phi_x)$ which are independent of the choice of these coordinates are thus properties of the topological type of $\cS_{\cP}$. 
}
\end{remark}

  \subsection{Equivalence of group chains}\label{subsec-chainequiv}
Fokkink and Oversteegen \cite{FO2002} and the authors \cite{DHL2016a} studied equivalences of group chains associated to a given   equicontinuous minimal Cantor system $(V_0,G_0,\Phi)$. We now briefly recall the key results. 

Denote by $\fG$   the collection of all possible subgroup chains in $G_0$.
Then there are two equivalence relations on $\fG$. The first was introduced by Rogers and Tollefson in \cite{RT1971b}.
\begin{defn} \cite{RT1971b}\label{defn-greq}
In a finitely generated group $G_0$, two group chains $\{G_{\ell}\}_{\ell \geq 0}$ and $\{H_{\ell}\}_{\ell \geq 0}$ with $G_0=H_0$ are \emph{equivalent}, if and only if, there is a group chain $\{K_{\ell}\}_{\ell \geq 0}$ and infinite subsequences $\{G_{\ell_k}\}_{k \geq 0}$ and $\{H_{j_k}\}_{k \geq 0}$ such that $K_{2k} = G_{\ell_k}$ and $K_{2k+1} = H_{j_k}$ for $k \geq 0$.
\end{defn}

The next definition was introduced by Fokkink and Oversteegen in \cite{FO2002}.

\begin{defn} \cite{FO2002}\label{conj-equiv}
Two group chains $\{G_{\ell}\}_{\ell \geq 0}$ and $\{H_{\ell}\}_{\ell \geq 0}$ in $\fG$ are \emph{conjugate equivalent} if and only if there exists a sequence $(g_{\ell}) \subset G_0$ for which  the compatibility condition   $g_{\ell}G_{\ell} = g_{\ell +1} G_{\ell}$ for all $\ell\geq 0$ is satisfied, and 
such that the group chains $\{g_{\ell} G_{\ell} g_{\ell}^{-1}\}_{\ell \geq 0}$ and $\{H_{\ell}\}_{\ell \geq 0}$ are equivalent. 
\end{defn}

The dynamical meaning of the equivalences in Definitions~\ref{defn-greq} and \ref{conj-equiv} is given by the following theorem, which follows from results  in \cite{FO2002}; see also \cite{DHL2016a}.  

\begin{thm}\label{equiv-rel-11}
Let $\{G_{\ell}\}_{\ell \geq 0}$ and $\{H_{\ell}\}_{\ell \geq 0}$ be group chains in $G_0$, with $H_0 = G_0$, and let 
\begin{eqnarray*}
G_\infty & = &  \lim_{\longleftarrow}\{G_0/G_{\ell+1} \to G_0/G_{\ell}\}  \ , \\
H_\infty & = &  \lim_{\longleftarrow}\{G_0/H_{\ell+1} \to G_0/H_{\ell}\}  \ .
\end{eqnarray*}
Then:
\begin{enumerate}
\item \label{er-item1} The group chains $\{G_{\ell}\}_{\ell \geq 0}$ and $\{H_{\ell}\}_{\ell \geq 0}$ are \emph{equivalent} if and only if there exists a homeomorphism $\tau \colon  G_\infty \to H_\infty$ equivariant with respect to the $G_0$-actions on $G_\infty$ and $H_\infty$, and such that $\phi(eG_{\ell}) = (eH_{\ell})$.
\item The group chains $\{G_{\ell}\}_{\ell \geq 0}$ and $\{H_{\ell}\}_{\ell \geq 0}$ are \emph{conjugate equivalent} if and only if there exists a homeomorphism $\tau \colon  G_\infty \to H_\infty$ equivariant with respect to the $G_0$-actions on $G_\infty$ and $H_\infty$.
\end{enumerate}
\end{thm}

That is, an equivalence of two group chains corresponds to the existence of a \emph{basepoint-preserving} equivariant homeomorphism between their inverse limit systems, while a conjugate equivalence of two group chains corresponds to the existence of a equivariant conjugacy between their inverse limit systems, which need not preserve the basepoint. 

Let  $\fG(\Phi_0)$ denote the class of group chains in $G_0$ which are \emph{conjugate equivalent} to the group chain $\{G^x_{\ell}\}_{\ell \geq 0}$ with basepoint $x$.
  The following result gives a geometric interpretation    the conjugate equivalence class $\fG(\Phi_0)$ of a group chain $\{G_{\ell}^x\}_{\ell \geq 0}$.
  
\begin{prop}\label{prop-conjequiv}
For an equicontinuous minimal Cantor action $(V_0,G_0,\Phi_0)$, let $\{G^x_{\ell}\}_{\ell \geq 0}$ be a group chain with partitions $\{\cQ_{\ell}\}_{\ell \geq 0}$ and basepoint $x$, as in Proposition~\ref{prop-AFpres}.
 Then a group chain $\{H_{\ell}\}_{\ell \geq 0}$  is in $\fG(\Phi_0)$ if and only if there exists a collection of $G_0$-invariant partitions $\cS_{\ell} = \{g \cdot U_{\ell}\}_{g \in G_0}$ of $V_0$, where $U_{\ell} \subset V_0$ is a clopen set, and ${\bigcap_{\ell} U_{\ell} = \{y\} \subset V_0}$, such that $H_{\ell} = H_{\ell}^y$ is the isotropy group  at $U_{\ell}$ of the action of $G_0$ on the partition $\cS_{\ell}$, for all $\ell \geq 0$.
\end{prop}

  \subsection{Kernels of group chains}\label{subsec-kernel}
The following notion is  important for the study of group chains. 
\begin{defn}\label{def-kernel}
The \emph{kernel} of a group chain  $\cG = \{G_{\ell}\}_{\ell \geq 0}$ is the subgroup of $G_0$ given by
\begin{equation}
K(\cG) = \bigcap_{\ell \geq 0} ~ G_{\ell} \ . 
\end{equation}
\end{defn}

The following property is immediate from the definitions.
\begin{lemma}\label{lem-equivkernels}
Suppose that the group chains $\cG = \{G_{\ell}\}_{\ell \geq 0}$ and $\cH = \{H_{\ell}\}_{\ell \geq 0}$ with $G_0=H_0$ are  equivalent, then $K(\cG) = K(\cH) \subset G_0$.
\end{lemma}
If the chains $\cG$ and $\cH$ are only conjugate equivalent, then the kernels need not be equal. 
 
An infinite  group $G_0$ which admits a group chain $\cC = \{C_{\ell}\}_{\ell \geq 0}$ where each $C_{\ell}$ is a \emph{normal} subgroup of $G_0$, and such that $\bigcap C_\ell = \{e\}$, where $e$ denotes the identity element in $G_0$, is said to be \emph{residually finite}. It is an elementary fact that given any group chain $\cG = \{G_{\ell}\}_{\ell \geq 0}$ in $G_0$, there is an associated core group chain $\cG$ for which $C_{\ell} \subset G_{\ell}$ with $C_{\ell}$ normal in $G_0$, for all $\ell > 0$, as will be discussed in Section~\ref{subsec-ellischains} below. Thus, if the group chain $\cG^x = \{G^x_{\ell}\}_{\ell \geq 0}$ introduced above has $K(\cG^x)$ the trivial group, then $G_0$ must be a residually finite group.  On the other hand, there are many classes of groups which are not residually finite, and thus any group chain for these groups must have non-trivial kernels. For example, many of the types of Baumslag-Solitar groups are not residually finite \cite{Levitt2015a,Levitt2015b,Meskin1972}, so every equicontinuous minimal Cantor system defined by an action of one of these groups will have non-trivial kernels.

The kernel $K(\cG^x)$ has an interpretation in terms of the topology of the leaves of the foliation $\FP$ of a weak solenoid. 
Let $(V_0,G_0,\Phi_0)$ be the holonomy action  for a weak solenoid $\cS_{\cP}$ with   presentation $\cP$ and basepoint $x \in V_0$, and   let $\cG^x = \{G_i^x\}_{i \geq 0}$ be the group chain at $x$.
 Recall that    the restriction of the bundle projection  $\Pi_0|_{L_x} \colon  L_x \to M_0$  to  the leaf $L_x$  containing $x$ is a covering map. 
Let  $\wtzM$ be the universal cover of $M_0$. Then by   standard arguments of covering space theory (see also McCord \cite{McCord1965}) there is a homeomorphism
\begin{equation}\label{eq-covering}
 \wtzM/ K(\cG^x) \to L_x \  .
\end{equation}
Now let $y \in  \fX_0$ be another point. Then by Remark~\ref{rmk-basepoints},  the group chain   associated to $y$ is given by $\cG^y = \{g_i G_{i}^x g_i^{-1}\}_{i \geq 0}$  where $\tau_x(y) = (g_i G_i^x)$.
If $y$ is in the orbit of $x$ under the $G_0$-action, then we can take $g_i = g$ for some $g \in G_0$, and thus  $K(\cG^y) = g K(\cG^x) g^{-1}$; that is, the kernels of $\cG^x$ and $\cG^y$ are conjugate, which corresponds to the fact that the fundamental group of the leaf $L_x$ at differing basepoints are conjugate.   
If $y$ is not in the orbit of $x$, then the relationship between $K(\cG^x)$ and $K(\cG^y)$ depends on the dynamical properties of the solenoid. 

In particular, in Section~\ref{sec-holonomy} we relate the algebraic properties of the kernels $K(\cG^y)$ with the germinal holonomy groups of the foliation $\FP$. Recall from Section~\ref{sec-intro} that a manifold $L$ has $\pi_1$-finite type if it's fundamental group is finitely generated. A matchbox manifold $\fM$ has finite $\pi_1$-type if all leaves in $\FfM$ have finite $\pi_1$-type. The following statement is immediate from the above discussion.

 \begin{lemma}\label{lemma-finitetype}
An equicontinuous matchbox manifold $\fM$ has finite $\pi_1$-type if and only if, for the associated group chain $\cG^x = \{G^x_{\ell}\}_{\ell \geq 0}$,  for all $\cG^y  \in \fG(\Phi)$, the kernel $K(\cG^y)$ is a finitely generated subgroup of $G_0$. 
 \end{lemma}

We next give two examples to illustrate the above concepts.
    
\begin{ex}\label{ex-constantkernel}
{\rm
Let $\cS_\cP$ be a Vietoris solenoid, as in Example~\ref{ex-Vietoris}, where $m_{\ell} >1$ is the degree of $p_{\ell}$. Choose $x \in \cS_\cP$ so that $\Pi_{\ell}(x) = 0$ for $\ell \geq 0$. Then $G_0 = \mZ$, and the group $G_{\ell}^x = \tilde{m}_{\ell} \mZ$, where $\tilde{m}_{\ell} = m_1m_2 \cdots m_{\ell}$ is the product of the degrees of the coverings. Then the kernel $K(\cG^x) = \{0\}$, and the path-connected component $L_x$ is homeomorphic to the real line. Let $y \in \fX_0$ be any other point in the fibre. Since $\mZ$ is abelian, any subgroup conjugate to $G_{\ell}^x = \tilde{m}_{\ell} \mZ$ is equal to it. It follows that $K(\cG^y) = \{0\}$, and $L_y$ is homeomorphic to the real line for any $y \in \fX_0$.

More generally, suppose $\cS_\cP $ is an $n$-dimensional solenoid, and $G_{\ell}^x$ is a normal subgroup of $G_0$ for all $\ell \geq 1$. Then for any $y \in \fX_0$ we have $\cG^y = \cG^x$, and so $K(\cG^y) = K(\cG^x)$. It follows that all leaves in $\cS_\cP$ are homeomorphic. The Vietoris solenoid $\cS_\cP$ is of finite $\pi_1$-type.
}
\end{ex}

\begin{ex}\label{ex-RT}
{\rm
This example is due to Rogers and Tollefson \cite{RT1972}. Consider a map of the plane, given by a translation by $\frac{1}{2}$ in the first component, and by reflection in the second component, i.e.
  $$r \times i \colon  \mR^2 \to \mR^2 ~ {\rm where}  ~  (x,y) \mapsto (x+\frac{1}{2},-y).$$
This map commutes with translations by the elements in the integer lattice $\mZ^2 \subset \mR^2$, and so induces the map $r \times i \colon  \mT^2 = \mR^2/\mZ^2 \to \mT^2$ of the torus. This map is an involution, and the quotient space $K = \mT^2/(x,y) \sim r\times i(x,y)$ is homeomorphic to the Klein bottle.

Consider the double covering map $L \colon  \mT^2 \to \mT^2$ given by  $L(x,y) = (x,2y)$. The inverse limit ${\ds \mT_\infty = \lim_{\longleftarrow} \{L \colon  \mT^2 \to \mT^2 \} }$ is a solenoid with $2$-dimensional leaves. Let $x_0 = (0,0) \in M_0 = \mT^2$. The fundamental group $G_0 = \mZ^2$ is abelian, so for any $x,y \in \fX_0$ the kernels $K(\cG^x) = K(\cG^y) \cong \mZ$, and every leaf is homeomorphic to an open two-ended cylinder.

The involution $r \times i$ is compatible with the covering maps $L$, and so it induces an involution $(r \times i)_\infty \colon \mT_\infty \to \mT_\infty$, which is seen to have a single fixed point $(0,0, \ldots) \in \mT_\infty$, and permute other path-connected components. Let $p \colon  K \to K$ be the double covering of  the Klein bottle by itself, given by $p(x,y) = (x,2y)$,  and consider the inverse limit space ${\ds K_\infty = \lim_{\longleftarrow} \{p \colon  K \to K \}}$. Note that taking the quotient by the involution $r \times i$ is compatible with the covering maps $L$ and $p$; that is, $p \circ (r \times i ) =  L$, and so induces the map $i_\infty \colon  \mT_\infty \to K_\infty$ of the inverse limit spaces. Under this map, the path-connected component of the fixed point $(0,0, \ldots)$ is identified so as to become a non-orientable one-ended cylinder. The image of any other path-connected component is an orientable $2$-ended cylinder.

Let $x =(x_{\ell})\in K_\infty$ for  $x_{\ell}  \in K$. Then $G_0 = \pi_1(K,x_0) = \langle a,b \mid bab^{-1} = a^{-1} \rangle$. Fokkink and Oversteegen \cite{FO2002} computed the kernel $K(\cG^x) = \langle b \rangle $ of the group chain $\cG^x$. They also computed kernels for group chains at any other basepoint $y \in \fX_0$ and found that either $K(\cG^y)$ is conjugate to $\langle b \rangle$, or $K(\cG^y)$ is equal to $\langle b^2 \rangle$. This example has finite $\pi_1$-type.
}
\end{ex}

  \section{Homogeneous solenoids and actions} \label{sec-homogeneous}

  In this section, we review the results from various works about the criteria for homogeneity of matchbox manifolds. These data will be of use later, when we give the proof Theorem~\ref{thm-molino}.

A continuum $\fM$ is said to be \emph{homogeneous} if given any pair of points $x,y \in \fM$, then there exists a homeomorphism $h \colon \fM \to \fM$ such that $h(x) = y$. 
A homeomorphism   $\varphi \colon \fM \to \fM$    preserves the path-connected components, hence   preserves the foliation $\FfM$ of $\fM$. It follows that if $\fM$ is   homogeneous, then it is also foliated homogeneous.  

By \cite[Theorem~5.2]{ClarkHurder2013} a homogeneous matchbox manifold $\fM$ is equicontinuous. Hence by Theorem~\ref{thm-equicontinuous} above, which is proved in \cite[Theorem~1.4]{ClarkHurder2013}, the foliated space $\fM$ is homeomorphic to a weak solenoid $\cS_{\cP}$. We   restrict our attention to equicontinuous foliated spaces, so consider the problem of giving conditions for when a weak solenoid $\cS_{\cP}$ is homogeneous, which is thus equivalent to asking for criteria when an  equicontinuous matchbox manifold is homogeneous. This is one of the original  motivating problems in the study of solenoids,   to obtain necessary and sufficient conditions that the solenoid $\cS_{\cP}$ is homogeneous \cite{FO2002,Rogers1970,RT1971a,Schori1966}. In this section, we recall the relevant results of these previous works,  and of the authors in \cite{Dyer2015,DHL2016a,DHL2016b}.

\subsection{Regular actions}

An \emph{automorphism} of $(V_0, G_0, \Phi_0)$ is a homeomorphism $h \colon V_0 \to V_0$ which commutes with the   $G_0$-action on $V_0$. Denote by ${Aut}(V_0, G_0, \Phi_0)$ the group of automorphisms of the action $(V_0, G_0, \Phi_0)$.  Note that ${Aut}(V_0, G_0, \Phi_0)$ is a topological group for  the compact-open topology on maps, and  is a closed subgroup of ${Homeo}(V_0)$.

\begin{defn}\label{def-regular}
The equicontinuous minimal Cantor action $(V_0,G_0,\Phi_0)$   is:
\begin{enumerate}
\item \emph{regular}  if the action of ${Aut}(V_0, G_0, \Phi_0)$ on $V_0$ has a single orbit;
\item  \emph{weakly normal} if the action of ${Aut}(V_0, G_0, \Phi_0)$ decomposes $V_0$ into a finite collection of  orbits;
\item \emph{irregular}  if the action of ${Aut}(V_0, G_0, \Phi_0)$ decomposes $V_0$ into an infinite collection of  orbits.
\end{enumerate}
\end{defn}

The terminology in Definition~\ref{def-regular} is chosen so that to be consistent with the terminology in \cite{DHL2016a,FO2002}.

Recall that    $\fG$  denotes the collection of all possible subgroup chains in $G_0$, and let 
  $\fG(\Phi_0) \subset \fG$ denote the collection of all group chains in $\fG$ which are conjugate equivalent to a given  group chain $\cG^x = \{G_{\ell}^x\}_{\ell \geq 0}$.  Theorem~\ref{equiv-rel-11} states that a group chain $\{G_{\ell}^x\}_{\ell \geq 0}$ is equivalent to the group chain $\{H_{\ell}^y\}_{\ell \geq 0}$ if and only if there exists a conjugacy $h \colon  V_0 \to V_0$ of the $G_0$-action on $V_0$, such that $h(x) = y$. Such an $h$ is   an automorphism of $(V_0,G_0,\Phi_0)$, which gives the following   result.  

\begin{thm}\label{thm-rephr}
Let $(V_0, G_0, \Phi_0)$ be an equicontinuous minimal Cantor action, and   $\{G^x_{\ell}\}_{\ell \geq 0} \in \fG$ be  a group chain associated to the action.  Then $(V_0, G_0, \Phi_0)$ is:
\begin{enumerate}
\item \emph{regular} if all group chains in $\fG(\Phi_0)$ are equivalent;
\item  \emph{weakly normal} if $\fG(\Phi_0)$ contains a finite number of   classes of equivalent group chains;
\item \emph{irregular} if $\fG(\Phi_0)$ contains an infinite number of   classes of equivalent group chains.
\end{enumerate}
\end{thm}
 
McCord   in \cite{McCord1965} studied the case  when the chain $\{G_{\ell}^x\}_{\ell \geq 0}$ consists of   normal subgroups of $G_0$. In this case, every quotient $X_{\ell}^x = G_0/G_{\ell}^x$ is a finite group, and the inverse limit $X_\infty^x$, defined by \eqref{eq-Galoisfiber}, is then a profinite group. The group $X_\infty^x$ is identified with $V_0$ as a topological space, and it acts transitively on $V_0$ on the right. The right action of $X^x_{\infty}$ commutes with the left action of $G_0$ on $X^x_{\infty}$, and thus  $X^x_{\infty} \subset {Aut}(V_0, G_0, \Phi_0)$, and  so  the automorphism group acts transitively on  $H_{\infty}$.
 McCord used this observation  in \cite{McCord1965} to show that the group ${Homeo}(\cS_{\cP})$ acts transitively on $\cS_{\cP}$, proving the following theorem.

  \begin{thm}\cite{McCord1965}\label{thm-normal}
  Let $\cS_\cP$ be a solenoid with a group chain $\{G_{\ell}^x\}_{\ell \geq 0}$, such that $G_{\ell}^x$ is a normal subgroup of $G_0$ for all $\ell \geq 0$. Then $\cS_\cP$ is homogeneous. 
 \end{thm} 
 For example, if $G_0$ is abelian, then every group chain $\{G_{\ell}^x\}_{\ell \geq 0}$ consists of normal subgroups, and the solenoid $\cS_\cP$ is homogeneous. 
 
\subsection{Weakly normal actions} 
We next consider the problem of giving necessary and sufficient conditions for when a solenoid $\cS_{\cP}$ is homogeneous.

The converse to Theorem~\ref{thm-normal} is not true.  Indeed,  
Rogers and Tollefson in \cite{RT1971b} gave an example of a weak solenoid for which the presentation yields a chain of subgroups which are not normal in $G_0$, yet the inverse limit is    a profinite group, and so the solenoid is homogeneous.  This example was the motivation for the work of  
Fokkink and Oversteegen in \cite{FO2002}, where they    gave a necessary and sufficient condition on the chain $\{G_{\ell}^x\}_{\ell \geq 0}$ for the weak solenoid to be homogeneous.
In particular, they proved the following result.
Let  $N_{G_0}(G_{\ell})$ denote the normalizer of the subgroup $G_{\ell}$ in $G_0$; that is, $N_{G_0}(G_{\ell}) = \{g \in G_0 \mid g \, G_{\ell} \, g^{-1} = G_{\ell}\}$. 

\begin{thm} \cite{FO2002} \label{thm-criteria}
Let $(V_0,G_0,\Phi_0)$ be an equicontinuous minimal Cantor   action, $x \in V_0$ be a point, and let $\{G_{\ell}^x\}_{\ell \geq 0}$ be an associated group chain with conjugate equivalence class $\fG(\Phi_0)$. Then
\begin{enumerate}
\item   $(V_0, G_0, \Phi_0)$ is regular if and only if there exists a group chain $\{N_{\ell}\}_{\ell\geq 0} \in \fG(\Phi_0)$ such that $N_{\ell}$ is a normal subgroup of $G_0$ for each $\ell \geq 0$. 
\item   $(V_0, G_0, \Phi_0)$ is weakly normal if and only if there exists $\{{G^x_{\ell}}'\}_{i\geq 0} \in \fG(\Phi_0)$ and an $n>0$ such that ${G_{\ell}^x}' \subset G_n^x \subseteq N_{G_0}({G_{\ell}^x}')$ for all $\ell \geq n$. 
\end{enumerate}
\end{thm}

In Theorem~\ref{thm-criteria}, the set $\fG(\Phi_0)$ contains group chains which are conjugate equivalent to the given chain $\{G_{\ell}^x\}_{\ell \geq 0}$. The condition that the group chain $\{N_{\ell}\}_{\ell\geq 0} $ consists of normal subgroups implies that every chain in $\fG(\Phi_0)$ is equivalent to $\{N_{\ell}\}_{\ell\geq 0} $, and so $\{G_{\ell}^x\}_{\ell \geq 0}$ is equivalent to $\{N_{\ell}\}_{\ell\geq 0} $. In statement $2)$, the condition  ${G_{\ell}^x}' \subset G_n^x \subseteq N_{G_0}({G_{\ell}^x}')$ implies that the group chain $\{{G_{\ell}^x}'\}_{\ell \geq 0}$ is equivalent to $\{G_{\ell}^x\}_{\ell \geq 0}$. Indeed, suppose that ${G_{\ell}^x}' \subset {G_m^x}' \subseteq N_{G_0}({G_{\ell}^x}')$ for some $m$. Then for $n \leq m$ and $\ell \leq n$ we have ${G_{\ell}^x}' \subset {G_n^x}' \subseteq N_{G_0}({G_{\ell}^x}')$. If $\{{G^x_{\ell}}'\}_{i\geq 0} $ is equivalent to $\{{G^x_{\ell}}\}_{i\geq 0} $, then for some $n \leq m$ we have ${G^x_{n}}' \subset {G^x_{n}} \subset {G^x_{m}}'$, which yields the statement.

 Recall that Proposition~\ref{prop-AFpres} introduced the descending chain of clopen sets $\{V_{\ell +1} \subset V_{\ell} \mid \ell \geq 0\}$ of  $V_0$ such that $V_{\ell}$ is stabilized by the action of $G_{\ell}$.  Thus, the  weak normality condition in Theorem~\ref{thm-criteria} implies that if we restrict the $G_0$ action to the clopen set   $V_n \subset V_0$, then the restricted action   $(V_n,G_n,\Phi_n)$ with associated group chain $\cG_n^x = \{ G_{\ell}^x\}_{\ell \geq n}$ is regular.  
  In the case where the group chain $\{G_{\ell} \}_{\ell \geq 0}$ is associated to a weak solenoid $\cS_{\cP}$,  restricting to the action $(V_n,G_n,\Phi_n)$ amounts to discarding the initial   manifolds $\{M_0,\ldots,M_{n-1}\}$ in the presentation $\cP$,  to obtain the presentation $\cP_n$ defined in \eqref{eq-Pn}. Then as discussed in  Section~\ref{subsec-homeos},   there is a homeomorphism  $\cS_{\cP_n} \cong \cS_{\cP}$,   where the homeomorphism is given by the ``shift'' map $\sigma_n$. Thus, $\cS_{\cP}$ is homogeneous if and only if $\cS_{\cP_n}$ is homogeneous, and so by Theorem~\ref{thm-normal} a weak solenoid whose associated group chain is  weakly normal  is homogeneous. We thus obtain the following result of  Fokkink and Oversteegen \cite{FO2002} giving  a criterion for when a weak solenoid is homogeneous.
 
 \begin{prop}\label{prop-fokking}\cite{FO2002} 
 Let $\cS_{\cP}$ be a weak solenoid, defined by a {presentation} $\cP$ with associated group chain $\{G_{\ell}^x\}_{\ell\geq 0}$. Then  $\cS_{\cP}$ is homogeneous if and only if  $\{G_{\ell}^x\}_{\ell\geq 0}$ is weakly normal. 
 \end{prop} 
   
   We also have the following property of presentations of homogeneous solenoids.
 \begin{prop}\label{prop-conjclasses}\cite{FO2002} 
 Let $\cS_{\cP}$ be a weak solenoid, defined by a {presentation} $\cP$ with associated group chain $\cG^x = \{G_{\ell}^x\}_{\ell\geq 0}$. If   $\cS_{\cP}$ is homogeneous, then the   kernel $K(\cG^x) \subset G_0$ has a finite number of conjugacy classes in $G_0$.
 \end{prop} 
   \proof
 Suppose that  $\cS_{\cP}$ is homogeneous.  Then by Theorem~\ref{thm-criteria}, there exists ${\cG^x}' = \{{G^x_{\ell}}'\}_{\ell\geq 0} \in \fG(\Phi_0)$ and an $n>0$ 
 such that ${G_{\ell}^x}' \subset G_n^x \subseteq N_{G_0}({G_{\ell}^x}')$ for all $\ell \geq n$.  Then ${G_{n}^x}' \subseteq N_{G_0}({G_{\ell}^x}')$ for all $\ell \geq n$, which implies that ${G_{n}^x}' \subset N_{G_0}(K(\cG_x')) $. Indeed, the chain $\{{G_{\ell}^x}'\}_{\ell \geq n} $ contains subgroups normal in ${G_{n}^x}'$, and it's intersection is then again normal in ${G_{n}^x}'$. Then for any $h \in {G_{n}^x}'$ we have
\begin{equation}
h \cdot K({\cG^x}') \cdot h^{-1}  =   
K({\cG^x}') \ ,
\end{equation}
and $K({\cG^x}') $ has only a finite number of conjugacy classes, at most $[G_0: {G_n^x}']$. Since $\cG^x$ is equivalent to ${\cG^x}'$, we have that ${G_0^x} = {G_0^x}' \supset {G_{1}^x} \supset {G_{1}^x}' \supset {G_{2}^x} \supset {G_{2}^x}' \supset \cdots$, and so $K({\cG^x})  = K({\cG^x}') $, which yields the statement.
   \endproof

\section{Ellis group of   equicontinuous   minimal  systems}\label{sec-ellischains}

The  \emph{Ellis (enveloping) semigroup} associated to a continuous group action $\Phi \colon G \times X \to X$   was introduced in the papers \cite{EllisGottschalk1960,Ellis1960}, and is treated in the books   \cite{Auslander1988,Ellis1969,Ellis2014}.   The construction of $\whE(X,G,\Phi)$ is abstract, and it can be a   difficult problem to calculate this group exactly. A  key problem is to understand the relation between  the algebraic properties of $\whE(X,G,\Phi)$ and dynamics of the action.
In this section, we briefly recall some basic properties   of $\whE(X,G,\Phi)$, then consider the  results for the special case of equicontinuous   minimal  systems.

 \subsection{Ellis (enveloping) group} \label{subsec-ellis}

 Let $X$ be a compact Hausdorff  topological space, and $G$ be a finitely generated group. Consider the space $X^X = {\it Maps}(X,X)$ with topology of \emph{pointwise convergence on maps}. With this topology, $X^X$ is a compact Hausdorff space. Each $g \in G$ defines an element $\whg \in Homeo(X) \subset X^X = {\it Maps}(X,X)$. Denote by $\whG$ the set of all such elements. Ellis \cite{Ellis1960} showed that the closure $\overline{\whG} \subset X^X$ has the structure of a right topological semigroup. Moreover, if the action $(X,G,\Phi)$ is equicontinuous, then the semigroup $\overline{\whG}$ is a group naturally identified with the closure  $\overline{\Phi(G)}$  of $\Phi(G)  \subset {\it Homeo}(X)$ in the \emph{uniform topology on maps}. Each element of $\overline{\Phi(G)}$ is the limit of a sequence of points in $\whG$, and we use the notation $(g_i)$ to denote a sequence $\{g_i \mid i \geq 1\} \subset G$ such that the sequence $\{\whg_i = \Phi(g_i) \mid i \geq 1\} \subset {\it Homeo}(X)$ converges in the uniform topology. 

Assume the action of $G$ on $X$ is minimal, that is,  the orbit ${\Phi(G)}(x)$ is dense in $X$ for any $x \in X$.  It then follows that the orbit of the Ellis group $\overline{\Phi(G)}(x) = X$ for any $x \in X$. That is, the group $\overline{\Phi(G)}$ acts transitively on $X$. Then for the isotropy group of the action at $x$,  
\begin{align}\label{iso-defn2}
\overline{\Phi(G)}_x = \{ (g_i) \in \overline{\Phi(G)} \mid (g_i) \cdot x = x\},
\end{align}
we have the natural identification $X \cong \overline{\Phi(G)}/\overline{\Phi(G)}_x$ of left $G$-spaces.

Given an  equicontinuous   minimal Cantor system $(X,G,\Phi)$, the  Ellis group     $\overline{\Phi(G)}$  depends only on the image  $\Phi(G) \subset Homeo(X)$.
On the other hand, the isotropy group $\overline{\Phi(G)}_x$ may depend on the point $x \in X$. 
Since the action of $\overline{\Phi(G)}$ is transitive on $X$, given any $y \in X$,   there exists $(g_i) \in \oG$ such that $(g_i) \cdot x  =y$. It follows that 
  \begin{align}\label{eq-profiniteconj}
  \oG_y = (g_i) \cdot \oG_x \cdot (g_i)^{-1}  \ .
  \end{align}
Thus,  the \emph{cardinality} of the isotropy group $\overline{\Phi(G)}_x$ is independent of the point $x \in X$, and so  the Ellis group $\overline{\Phi(G)}$ and the cardinality of $\overline{\Phi(G)}_x$ are invariants of $(X,G,\Phi)$.

\subsection{Ellis   group for group chains}\label{subsec-ellischains}
We consider the Ellis group for an equicontinuous minimal Cantor action $(V_0,G_0,\Phi)$, in terms of an associated group chain   $\cG^x = \{G_{\ell}^x\}_{\ell \geq 0}$ for $x \in V_0$. For each subgroup $G^x_{\ell}$ consider the maximal normal subgroup of $G_{\ell}^x$ which is given by
\begin{equation}\label{eq-core}
C_{\ell} \equiv   {\rm core}_{G_0}   \, G_{\ell}^x \equiv  \bigcap_{g \in {G_0}} gG_{\ell}^xg^{-1}  ~ \subseteq ~ G_{\ell}^x \ .
\end{equation}
The group $C_{\ell}$ is called the \emph{core} of $G_{\ell}$ in $G_0$. Since $C_{\ell}$ is normal in $G_0$, the quotient $G_0/C_{\ell}$ is a finite group, and  the collection  $\cC = \{C_{\ell} \}_{\ell \geq 0}$ forms a descending chain of normal subgroups of $G_0$. The inclusions of coset spaces define bonding maps $\delta^{\ell+1}_{\ell} $ for the inverse sequence of quotients $G_0/C_{\ell}$, and the inverse limit space 
\begin{eqnarray}
  C_{\infty}  & = &  \{(eG_0, g_1 C_1, \ldots) \mid g_{\ell} C_{\ell} = g_{\ell+1} C_{\ell}\} ~ \subset \prod_{\ell \geq 0} ~ G_0/C_{\ell}  \label{cinfty-coords} \\
  & \cong &  \lim_{\longleftarrow} \, \left\{\delta^{\ell+1}_{\ell} \colon   G_0/C_{\ell+1} \to G_0/C_{\ell}   \right\}  \label{cinfty-define}  
\end{eqnarray} 
is a profinite group. Let  $\whiota \colon G_0 \to C_{\infty}$ be the homomorphism defined by $\whiota(g) = (gC_{\ell})$ for $g \in G_0$. 
 Then the induced left action of $G_0$ on $C_{\infty}$ yields a minimal Cantor system, denoted by   $(C_{\infty} , G_0, \whPhi_0)$.

Also, introduce the descending chain of clopen neighborhoods of the identity $(eC_{\ell}) \in C_{\infty}$, which for $n \geq 0$ defines a neighborhood system for $C_{\infty}$:
\begin{eqnarray}  
C_{n,\infty}     & = & \{ (g_{\ell} C_{\ell}) \in C_{\infty} \mid   g_n \in C_n   \} \ ,   \label{Cinfty-definen} \\
& \cong &  \lim_{\longleftarrow} \, \left\{\delta^{\ell+1}_{\ell} \colon   C_n/C_{\ell +1}  \to C_n/C_{\ell}  \mid \ell \geq n \right\}   \ .   \label{Cinfty-definenn}
 \end{eqnarray}

\subsection{The discriminant} \label{subsec-disc}

Observe that for each $\ell \geq  0$, the    quotient group $D_{\ell}^x = G_{\ell}^x/C_{\ell} \subset G_0/C_{\ell}$. It follows that  the inverse limit space 
\begin{equation}\label{eq-discriminantdef}
\cD_x = \lim_{\longleftarrow}\, \left \{\delta^{\ell+1}_{\ell} \colon  D_{\ell+1}^x \to D_{\ell}^x \right\}
\end{equation}
is a closed subgroup of $C_\infty$. The group $\cD_x$ is called the \emph{discriminant group} of the action $(V_0,G_0,\Phi_0)$.

The   relationship between $C_\infty$ and the Ellis group of  $(V_0,G_0,\Phi_0)$ is given by the following result.

\begin{thm}[Theorem~4.4, \cite{DHL2016a}]\label{thm-quotientspace}
Let    $(V_0,G_0,\Phi_0)$ be an   equicontinuous minimal Cantor action, let $x \in V_0$,  and let $\cG^x \equiv \{G_{\ell}^x\}_{i \geq 0}$ be the associated group chain at $x$.
Then there is a natural isomorphism of topological groups $\whTheta \colon  \overline{\Phi(G_0)} \cong C_{\infty}$ such that the restriction $\whTheta \colon \overline{\Phi(G_0)}_x \cong \cD_x$.
\end{thm}

Moreover, the discriminant subgroup is simple by the next result.

\begin{prop}[Proposition~5.3, \cite{DHL2016a}]\label{prop-discrisnotnormal}
Let $(V_0,G_0,\Phi_0)$ be an equicontinuous   minimal Cantor system,   $x \in V_0$ a basepoint, and $\overline{\Phi_0(G_0)}_x$ the isotropy group of $x$. Then  
\begin{equation}\label{eq-rationalcore}
 {\rm core}_{G_0}  \overline{\Phi_0(G_0)}_x  =  \bigcap_{k \in G_0} k \ \overline{\Phi_0(G_0)}_x \ k^{-1}  
\end{equation}
 is the trivial group. Thus, the maximal normal subgroup of $\overline{\Phi_0(G_0)}_x$ in $\overline{\Phi_0(G_0)}$ is also   trivial.
\end{prop}

    We next consider the  homogeneity properties of a solenoid $\cS_{\cP}$   in terms of   $\cD_x$ (see  \cite{DHL2016a}.) 
  It follows from Proposition~\ref{prop-discrisnotnormal} that if $\cD_x$ is non-trivial, then it is not normal in $C_\infty$, and therefore the quotient $X_\infty^x = C_\infty/\cD_x$ is not a group. We thus conclude:

\begin{prop}\cite{DHL2016a}\label{prop-regularaction}
The action $(V_0,G_0,\Phi_0)$ is regular if and only if   $\cD_x$ is trivial.
\end{prop}

 Note that Proposition~\ref{prop-regularaction} does not take into account the possibility that the action of a subgroup $G_{\ell}^x$ on a smaller section $V_{\ell}$ is regular. 
 The general formulation is then as follows.
 
 \begin{cor}\label{cor-homogeneoussolenoid}
 An equicontinuous matchbox manifold $\fM$ is homogeneous if and only if it admits a transverse section $V_0$ and a presentation $\cP$ with associated group chain $\{G_{\ell}^x\}_{\ell \geq 0}$, such that the discriminant group $\cD_x$ is trivial.
 \end{cor}

\section{Molino theory for weak solenoids}\label{sec-molino}

  In this section, we obtain a ``Molino theory'' for   weak solenoids, and hence for all equicontinuous matchbox manifolds, including those for which the hypotheses of the work \cite{ALM2016} are not satisfied. There are often subtle, and not so subtle,  differences between the theory for matchbox manifolds  and  for smooth Riemannian foliations, as will be discussed further in the following sections.

 \subsection{Molino overview}
 
  Molino theory for Riemannian foliations gives a   structure theory for the geometry and dynamics of this class of foliations on compact smooth manifolds.  The S\'eminaire Bourbaki article \cite{Haefliger1989} by Haefliger gives a concise overview of the theory and its applications, and Molino's book \cite{Molino1988}  and its multiple appendices give a more detailed treatment of this theory and its applications. The book \cite{MoerdijkMrcunbook2003} is also an excellent reference to read about the essentials of Molino theory. We give a very brief summary below of some key properties of   the ``Molino space'' $\whM$ associated to a smooth Riemannian foliation $\F$ of a compact connected manifold $M$.
  
 Given a Riemannian foliation $\F$ of a compact connected manifold $M$, 
the associated \emph{Molino space} $\whM$ is a compact connected manifold with a Riemannian foliation $\whF$ whose leaves have the same dimension as those of $\F$. In the case where $\F$ is a minimal foliation, in the sense that   each leaf of $\F$ is dense in $M$, then we can assume that the foliation $\whF$ is also minimal.

Associated to a minimal Riemannian foliation $\F$  is   the \emph{structural Lie algebra} $\mathfrak{h}$, given by the algebra of    holonomy invariant normal vector fields to $\F$, and which  is well-defined up to  isomorphism.

There is a  fibration  $\whpi \colon  \whM \to M$ equipped with a fiber-preserving right action of a connected Lie group $H$  whose Lie algebra is $\mathfrak{h}$, and  for which the foliation $\whF$ is invariant under the     action of $H$. Moreover, for each leaf $\whL \subset \whM$, there is a leaf $L \subset M$ such that the restriction $\whpi \colon \whL \to L$ is the holonomy covering of $L$.    We say that  
\begin{equation}\label{eq-molinofibration}
H \longrightarrow \whM \stackrel{\whpi}{\longrightarrow} M
\end{equation}
 is a \emph{Molino sequence} for $M$, and $H$ is the structural Lie group for $\whF$.

  A key property of the Molino space $\whM$ of $\F$ is that it is \emph{Transversally Parallelizable}, or \emph{TP}. This condition states that there are non-vanishing vector fields $\{\vec{v}_1, \ldots , \vec{v}_q\}$ on $M$ which span the normal bundle to $\F$ at each $x \in M$, and the vector fields are locally projectable. As a consequence, given any pair of points $x,y \in \whM$ there exists a diffeomorphism $h \colon \whM \to \whM$ which maps leaves of $\whF$ to leaves of $\whF$, and satisfies $h(x) = y$. A foliation  $\whF$ satisfying this condition is said to be \emph{foliated homogeneous}.

\subsection{Molino sequences for  weak solenoids} 
For a  matchbox manifold,   the \emph{TP} condition cannot be defined, as the transversal space to the foliation is totally disconnected.  Thus, we need an alternative approach to defining the Molino fibration   \eqref{eq-molinofibration} in the case where the transversal space to the foliation is a Cantor set. 
The basic observation is that   the foliated homogeneous condition for $\whM$ admits a natural generalization to all foliated spaces, as  discussed for  weak solenoids   in Section~\ref{sec-homogeneous}. 
For weak solenoids,  we will see below that the structural Lie group $H$ is replaced by   the discriminant subgroup $\cD_x \subset C_{\infty}$ of Section~\ref{subsec-disc}, and the foliated homogeneous condition is a consequence of the Ellis group construction.   We now restate and prove Theorem~\ref{thm-molino}.

 \begin{thm}\label{thm-molino-restate}
Let $\fM$ be an equicontinuous matchbox manifold, and let $\cP$ be a presentation of $\fM$, such that $\fM$ is homeomorphic to a solenoid $\cS_\cP$. Then there exists a homogeneous matchbox manifold $\whfM$ with foliation $\whF$, called a ``Molino space'' of $\fM$, a  compact totally disconnected  group $\cD$, and a  fibration
 \begin{equation}\label{eq-molinoseqM2}
\cD  \longrightarrow \whfM \stackrel{\whq}{\longrightarrow} \fM \ , 
\end{equation}
where   the restriction of $\whq$ to each leaf in $\whfM$ is a covering map of some leaf in $\fM$.   We say that \eqref{eq-molinoseqM2} is a \emph{Molino sequence} for $\fM$.
 \end{thm}

 \proof Let  $V_0 \subset \fM$ be a transverse section  to the foliation $\FfM$ of $\fM$, as given in Proposition~\ref{prop-AFpres}, and let $x \in V_0$ be a choice of a basepoint. 
 Let  $G_0$ be the restricted holonomy group acting on $V_0$.
 Let $\cP = \{ p_{\ell+1} \colon M_{\ell+1} \to M_{\ell} \mid \ell \geq 0\}$ be a presentation at $x$ such that there is a homeomorphism $\fM \cong \cS_\cP$, and for $x \in V_0$ let $\cG^x = \{G_{\ell}^x\}_{\ell \geq 0}$ be the associated group chain in $G_0 = \pi_1(M_0, x_0)$. Let $\Pi_0 \colon \cS_\cP \to M_0$ and set $\fX_0 = \Pi_0^{-1}(x_0)$. Let    $\tau \colon V_0 \to \fX_0$ with $\tau_x(x) = (e G_{\ell}^x)$ be the homeomorphism defined in Lemma~\ref{lem-fibremodel}.
 
Recall that the covering map $q_{\ell}   \colon M_{\ell} \to M_0$ defined in \eqref{eq-coverings} is associated to the subgroup  $G_{\ell}^x \subset G_0 = \pi_1(M_0 , x_0)$. 
Recall that the     core subgroup  $C_{\ell} \subset G_{\ell}^x$ is the maximal normal subgroup of $G_0$ contained in $G_{\ell}^x$, and   has finite index in   $G_{\ell}^x$.
For each $\ell > 0$,  let $\whq_{\ell} \colon \whM_{\ell} \to M_0$ be the proper covering space associated to the normal subgroup $C_{\ell}$. Each inclusion $C_{\ell +1} \subset C_{\ell}$ induces a normal covering map  $\whp_{\ell+1} \colon \whM_{\ell+1} \to \whM_{\ell}$, and so yields a presentation
 $\ds \whcP = \{\whp_{\ell+1} \colon \whM_{\ell+1} \to \whM_{\ell} \mid \ell \geq 0\}$. 
 
\begin{defn}\label{def-molino}
The \emph{Molino space} associated to a weak solenoid $\cS_{\cP}$ defined    by a presentation $\cP$ is the inverse limit space associated to the presentation $\whcP$, 
\begin{equation}\label{eq-molino}
\whcS_{\cP}   \equiv \lim_{\longleftarrow} ~ \{ \whp_{\ell +1} \colon \whM_{\ell +1} \to \whM_{\ell}\} \ .
\end{equation}
Let $\whPi_0 \colon \whcS_{\cP} \to M_0$ be the projection map, with fiber $\whfX_0 = \whPi_0^{-1}(x_0)$. 
\end{defn}

 We state some of the basic properties of the space $\whcS_{\cP}$. The proofs of the following statements are omitted, as they  follow  by    arguments  analogous to the corresponding statements for  $\cS_{\cP}$.
 \begin{prop}\label{prop-molino}
Let  $\cS_{\cP}$ be a weak solenoid defined    by a presentation $\cP$, and let $\whcS_{\cP}$ be the solenoid defined by \eqref{eq-molino}. Then we have:
\begin{enumerate}
\item  There is a natural isomorphism   $\whfX_0  \cong C_{\infty}$ where $C_{\infty}$ is the profinite group defined by \eqref{cinfty-define}   ;
\item  There is a natural map of fibrations $\whq \colon \whcS_{\cP} \to \cS_{\cP}$, whose fiber over $x \in \fX_0$ is the discriminant group $\cD_x$   ; 
\item  The global holonomy of the fibration $\whPi_0 \colon \whcS_{\cP} \to M_0$ is naturally conjugate 
as $G_0$-actions with the minimal Cantor system $(C_{\infty} , G_0, \whPhi_0)$.
\end{enumerate}
 \end{prop}
 
 \begin{defn} \label{def-molinoseq}
 The \emph{Molino sequence} for the  weak solenoid $\cS_{\cP}$   is the principal fibration
 \begin{equation}\label{eq-molinoseq}
\cD_x \longrightarrow \whcS_{\cP} \stackrel{\whq}{\longrightarrow} \cS_{\cP} \ .
\end{equation}
 \end{defn}
 
 Proposition~\ref{prop-molino}(3) implies that the foliation $\whF_{\cP}$ on $\whcS_{\cP}$ is minimal, and the restrictions of $\whq$ to the leaves of $\whF_{\cP}$ are covering maps  by construction, as there is a covering map $\whM_{\ell} \to M_{\ell}$ for each $\ell \geq 1$ which induces $\whq$. Finally, the space $\whcS_{\cP}$ is homogeneous by  Proposition~\ref{prop-fokking}, as it is defined using the normal group chain $\{C_{\ell}\}_{\ell \geq 0}$.  
 
Set $\whfM = \whcS_{\cP}$ and $\cD = \cD_x$, then we have established the claims of Theorem~\ref{thm-molino-restate}.   \endproof
 
The construction of the sequence in \eqref{eq-molinoseq} may   depend on the various choices made, and this is a fundamental aspect of the ``Molino Theory'' for weak solenoids. 
We   consider    in the next Section~\ref{subsec-stability} the dependence of the discriminant group on the partition sets $V_n \subset V_0$. Then in Section~\ref{sec-holonomy},  we consider the dependence of the sequence   \eqref{eq-molinoseq} on the choice of the basepoint $x \in V_0$  and   the role of the holonomy of the leaf $L_x$ in the properties of $\cD_x$.

\subsection{Stability of the Molino sequence}\label{subsec-stability}
We next consider    the stability of the discriminant  group  for an equicontinuous  Cantor minimal system $(V_0,G_0,\Phi_0)$    when one restricts to a   section $V_n \subset V_0$.

  We start with an example that   highlights the importance  of the ``asymptotic algebraic structure'' of the group chain $\cG^x$  
   for the definition of the Molino space. 
Consider   a weak solenoid  $\cS_{\cP}$  with associated group chain $\cG^x = \{G_{\ell}^x \}_{\ell \geq 0}$ defined by the holonomy action $(V_0, G_0, \Phi_0)$ for a clopen subset $V_0 \subset \fX_0$, and suppose that $\cG^x$ is not regular. Then by Proposition~\ref{prop-regularaction}, the discriminant group $\cD_x$ is non-trivial, and thus the sequence \eqref{eq-molinoseq} has non-trivial fiber.  Now suppose that, in addition, the group chain $\cG^x$ is weakly normal. Then by Theorem~\ref{thm-criteria}, there exists some $n > 0$ such that the restricted action $(V_n, G_n, \Phi_n)$ is regular, hence the discriminant group $\cD_x^n$ for  the truncated chain $\cG_n^x = \{G_{\ell}^x\}_{\ell \geq n}$ associated to the restricted action is trivial. For the truncated presentation   $\cP_n$ defined by \eqref{eq-Pn},   we have $\whcS_{\cP_n} = \cS_{\cP_n}$ as $\cD_x^n$ is the trivial group, and $\cS_{\cP_n} \cong \cS_{\cP}$ as remarked in Section~\ref{subsec-homeos},  hence we can consider $\whcS_{\cP_n}$ as a Molino space for $\cS_{\cP}$ as well.  That is, for this choice of $V_n$ as a section, the Molino sequence \eqref{eq-molinoseq} has trivial fiber.

We next develop a comparison, for $n \geq 0$,  of the  discriminant groups $\cD_x^n$ for the    group chain $\cG_n^x $ associated to the  truncated presentation $\cP_n$ defined by \eqref{eq-Pn}.  We work with the group chain model $(X_{\infty}^x , G_0 , \Phi_x)$ of Lemma~\ref{lem-fibremodel} for 
the holonomy action $\Phi_0 \colon G_0 \to Homeo(\fX_0)$.
By   definition   \eqref{eq-discriminantdef} of the discriminant group, it suffices to consider this invariant in sufficiently small clopen neighborhoods of the identity in the core group associated with the group chains.   For $n \geq 0$, we have the   clopen neighborhoods of $\{e\} \in X_{\infty}$: 
  \begin{eqnarray}  
U_{n}     & = & \{ (g_{\ell} G_{\ell}) \in X_{\infty} \mid   g_n \in G_n^x   \} \subset X_{\infty}    \label{Ginfty-definen} \\
 & \cong &  \lim_{\longleftarrow} \, \left\{\delta^{\ell+1}_{\ell} \colon   G_n^x/G_{\ell +1}^n \to G_n^x/G_{\ell}  \mid \ell \geq n \right\}   \ .   \label{Ginfty-definenn}
    \end{eqnarray}
Note that $U_n$ is just the inverse limit group defined by the truncated group chain $\cG_n^x$. Next, we introduce the core   groups of $\cG_n^x$ for arbitrarily small neighborhoods of $\{e\} \in U_n$. For $\ell \geq n \geq 0$, set
\begin{equation}\label{eq-Mn}
E_{n,\ell} \equiv   {\rm core}_{G_n^x}   \, G_{\ell}^x \equiv  \bigcap_{g \in {G_n^x}} gG_{\ell}^x g^{-1}    \ .
\end{equation}
Note that $E_{0,\ell}= C_{\ell}$, and that for all $m \geq n \geq 0$ and   $\ell > m$, we have $\ds E_{n,\ell}  \subset E_{m,\ell} \subset G_{\ell}^x$.
 
 For $k \geq n \geq 0$,  define the clopen neighborhood $V_{n,k}$ of $\{e\}$ for the core group of $\cG_n^x$ by 
   \begin{eqnarray}  
V_{n,k}     & = & \{ (g_{\ell} E_{n,\ell})   \mid  \ell \geq k   \ , \  g_{k} \in G_{k}^x \ , \  g_{\ell +1} E_{n,\ell} = g_{\ell} E_{n,\ell}  \}     \label{Dinfty-definen}   \\
 & \cong &  \lim_{\longleftarrow} \, \left\{\delta^{\ell+1}_{\ell} \colon   G_{k}^x/E_{n,\ell} \to G_k^x/E_{n,\ell + 1}  \mid \ell \geq k  \right\}   \ .  
   \end{eqnarray}
Then $V_{n,n}$ is the core limit group, or the Ellis group,   for the  truncated group chain $\cG_n^x$,   and $\{e\} \in V_{n,k} \subset V_{n,n}$ for all $k \geq n$.
 Note also that $V_{0,0} = C_{\infty}$ is the Ellis group for $\cG^x$.

For each $\ell \geq k \geq m \geq n$, the inclusions $E_{n,\ell} \subset E_{m,\ell}$ induces   group surjections  
\begin{align}\label{eq-quotiensmn} 
  G_{k}^x/E_{n,\ell} \stackrel{\phi_{k,n,m}^{\ell}}{\longrightarrow} G_{k}^x/E_{m,\ell} \ ,
\end{align}
  so we obtain   surjective homomorphisms of profinite groups 
   $\ds \phi_{n,m} \colon   V_{n,k}  \to V_{m,k}$ for each $m > n \geq 0$.
In particular, for $k=m$, this states that the clopen neighborhood $V_{n,m}$ of $\{e\}$ in the limit core   group for $\cG_n^x$ maps onto the limit core group $V_{m,m}$ of $\cG_m^x$.
     
We consider next the  discriminant groups associated to the group chains  $\cG_n^x$ for $n \geq 0$, $\cD_x^n    \subset V_{n,n}$,  
\begin{eqnarray} 
\cD_x^n & = &    \lim_{\longleftarrow}\, \left \{\delta^{\ell+1}_{\ell} \colon  G_{\ell+1}^x/E_{n,\ell+1}  \to G_{\ell}^x/E_{n,\ell} \mid \ell \geq n\right\} \label{eq-discquotients1} \\
  & \cong &    \lim_{\longleftarrow}\, \left \{\delta^{\ell+1}_{\ell} \colon  G_{\ell+1}^x/E_{n,\ell+1}  \to G_{\ell}^x/E_{n,\ell} \mid \ell \geq m \right\}  \ , \ {\rm for} ~ m \geq n \ . \label{eq-discquotients2}
\end{eqnarray}
 It follows from  \eqref{eq-quotiensmn}  and \eqref{eq-discquotients2}  that for $m > n$,  there are    surjective homomorphisms:
\begin{equation}\label{eq-discmapsnm2}
   \cD_{x}  ~ \stackrel{~ \psi_{0,n} ~ }{\longrightarrow} ~  \cD_{x}^n ~ \stackrel{~ \psi_{n,m} ~}{\longrightarrow} ~  \cD_{x}^m \ .
\end{equation}

\begin{defn}\label{def-stableGC}
A group chain $\cG^x = \{G_{\ell}^x\}_{\ell \geq 0}$ is said to be \emph{stable} if there exists $n_0 \geq 0$ such that the maps $\ds \psi_{n,m} \colon \cD_{x}^n  \to    \cD_{x}^m$ defined in \eqref{eq-discmapsnm2} are isomorphisms for all $m \geq n \geq n_0$. Otherwise, the group chain is said to be \emph{wild}.
\end{defn}

  Theorem~\ref{thm-criteria} implies that if the group chain $\{G_{\ell}^x\}_{\ell \geq 0}$ is weakly normal, then it is stable, 
 as there exists some $n_0 \geq 0$ such that  $\cD_{x}^n$ is the trivial group for all $n \geq n_0$. This discussion and Lemma~\ref{lem-finitestable} yield Proposition~\ref{prop-stablesolenoids} of the Introduction.
 
 \begin{lemma}\label{lem-finitestable}
 If the discriminant group $\cD_x$ for   $\cG^x = \{G_{\ell}^x\}_{\ell \geq 0}$ is finite, then $\cG^x$ is stable.
 \end{lemma}
 \proof
 The map $\psi_{0,n} \colon \cD_x \to \cD_x^n$ is surjective for all $n \geq 0$, so the assumption that the cardinality   $\#\cD_x < \infty$ implies that the cardinality $\# \cD_x^n$ of the group $\cD_x^n$ is decreasing with $n$, and thus there exists $n_0 \geq 0$ such that the cardinality of its image must stabilize for   $n \geq n_0$. 
 Then for $n \geq n_0$,  the homomorphism $\psi_{n,m} \colon \cD_x^{n_0} \to \cD_x^{n}$ is an isomorphism.
 \endproof

\subsection{Stable matchbox manifolds}\label{subsec-stableMM}
 We next consider the relationship between the notion of stable  for a matchbox manifold as given in Definition~\ref{defn-stablesolenoid}, and stable for a group chain as given in Definition~\ref{def-stableGC}.

  Let $\fM$ be an equicontinuous matchbox manifold, 
   let  $V_0$ be a transverse section in $\fM$ as given in Proposition~\ref{prop-AFpres}, and let $x \in V_0$ be a choice of a basepoint. 
Let $V_{\ell}$ be defined as in Proposition~\ref{prop-AFpres}, so that $x \in V_{\ell}$ for all $\ell \geq 0$.  
 Let  $G_0$ be the group of transformations of $V_0$ which induces the restricted holonomy group acting on $V_0$, and let $G_{\ell}^x \subset G_0$ be the stabilizer group of the set $V_{\ell}$. 
  Let $\cG^x = \{G_{\ell}^x\}_{\ell \geq 0}$ be the associated group chain in $G_0 = \pi_1(M_0, x_0)$, 
  let $\cP_n$ be the   presentation \eqref{eq-Pn} associated to the truncated  group chain $\cG^x_n = \{G_{\ell}^x\}_{\ell \geq n}$,  and let $\cS_{\cP_n}$ be the inverse limit solenoid.  For each $n \geq 0$, let $\whcS_{\cP_n}$ be the homogeneous solenoid associated to the normal group chain $\{E_{n, \ell}\}_{\ell \geq n}$ defined by \eqref{eq-Mn}.

Assume that the group chain $\cG^x$ is stable in the sense of Definition~\ref{def-stableGC}. That is, there exists an  index $n_0$, such that for any $m > n \geq n_0$ restricting to the smaller sections $V_m \subset V_n \subset V_0$ with induced presentations $\cP_m$ and $\cP_n$, then the induced map $\ds \psi_{n,m} \colon \cD_x^n \to \cD_x^m$ in \eqref{eq-discmapsnm2}   is a topological isomorphism. Then we have a commutative diagram of fibrations:
   \begin{align}\label{eq-diagramkernel} 
    \xymatrixcolsep{4pc}
   \xymatrix{
   \cD_x^n  \ar[d]   \ar[r]^{\psi_{n,m}} & \cD_x^m \ar[d]  \\
   \whcS_{\cP_n} \ar[d] \ar[r]^{\widehat{\sigma}_{m-n}}  & \whcS_{\cP_m} \ar[d]   \\
  \cS_{\cP_n}   \ar[r]^{\sigma_{m-n}} & \cS_{\cP_m}    
   }  
   \end{align}
By the discussion in Section~\ref{subsec-homeos}, the shift map $\sigma_{m-n}$ is a homeomorphism, and by assumption, the   map   $\psi_{n,m} \colon \cD_n \cong \cD_m$   is a topological isomorphism. Thus  the map  
$\ds \widehat{\sigma}_{m-n} \colon    \whcS_{\cP_n} \to  \whcS_{\cP_m}$ is a homeomorphism. Hence, the Molino sequences for the presentations $\cP_n$ and $\cP_m$ yield isomorphic topological fibrations.  
Conversely, if the topological type of the Molino sequence 
\begin{equation}\label{eq-molinosequence}
\cD_x^n  \longrightarrow   \whcS_{\cP_n} \longrightarrow \cS_{\cP_n}
\end{equation}
is well-defined up to homeomorphism of fibrations,  for given $V_0$ and $n \geq 0$ sufficiently large, then   there exists $n_0 \geq 0$ such that $m > n \geq n_0$ implies that  $\ds \cD_x^n  \stackrel{\psi_{n,m}}{\longrightarrow}   \cD_x^m$ is a topological isomorphism. Thus, the map of fibers $\ds \psi_{n,m} \colon \cD_x^n \to \cD_x^m$ is a topological isomorphism, and hence $\cG^x$ is  stable.

  The following statement summarizes these conclusions.  
 \begin{thm}\label{thm-sections}
 Let $\fM$ be an equicontinuous matchbox manifold, 
   let  $V_0$ be a transverse section in $\fM$ as given in Proposition~\ref{prop-AFpres}, and let $x \in V_0$ be a choice of a basepoint. 
Let $V_{\ell}$ be defined as in Proposition~\ref{prop-AFpres}, so that $x \in V_{\ell}$ for all $\ell \geq 0$.  
 Let  $G_0$ be the restricted holonomy group acting on $V_0$, and let $G_{\ell}^x \subset G_0$ be the stabilizer group of the set $V_{\ell}$. 
  Let $\cG^x = \{G_{\ell}^x\}_{\ell \geq 0}$ be the associated group chain in $G_0 = \pi_1(M_0, x_0)$, 
  let $\cP_n$ be the   presentation \eqref{eq-Pn} associated to the truncated  group chain $\cG^x_n = \{G_{\ell}^x\}_{\ell \geq n}$,  and let $\cS_{\cP_n}$ be the inverse limit solenoid.  For each $n \geq 0$, let $\whcS_{\cP_n}$ be the homogeneous solenoid associated to the normal group chain $\{E_{\ell}^n\}_{\ell \geq n}$ defined by \eqref{eq-Mn}.

\begin{enumerate}
\item If $\cG^x $ is stable, then there exists $n_0 \geq 0$ such that for all $n \geq n_0$ the fibration \eqref{eq-molinosequence}
is a Molino sequence for $\fM \cong \cS_{\cP_n}$, and the fiber group   $\cD_{x}^n$ is well-defined up to topological isomorphism.\\
\item If   $\cG^x $ is wild, then the topological isomorphism type of the fiber  in the sequence \eqref{eq-molinosequence} does not stabilize as $n$ tends to infinity.
\end{enumerate}
\end{thm}
 
 Theorem~\ref{thm-sections} implies that the Molino sequence of a matchbox manifold $\fM$ need not be well-defined, though if the associated group chain   $\cG^x$ is stable, then $\fM$ does have  a well-defined Molino sequence.

\section{Germinal holonomy in solenoids}\label{sec-holonomy}

In this section, we investigate the relationship between the germinal holonomy groups of leaves in a solenoid, the kernels of the associated group chains, and the discriminant group of the action. 

Let $\fM$ be an equicontinuous matchbox manifold with transverse section $V_0$, let $x \in V_0$ be a point, and let $\cP = \{f^{\ell+1}_i \colon  M_{\ell+1} \to M_{\ell}\}$ be a presentation with associated group chain $\cG^x=\{G_{\ell}^x\}_{\ell \geq 0}$ in $G_0 = \pi_1(M_0,x_0)$. Then by Theorem~\ref{thm-equicontinuous}, there is a foliated homeomorphism $\fM \cong \cS_{\cP}$.  

Let ${\ds C_\infty = \lim_{\longleftarrow}\{G_0/C_{\ell+1} \to G_0/C_{\ell}\} }$, where $C_{\ell}$ is the maximal normal subgroup of $G_{\ell}^x$, $\ell \geq 0$, and let $\cD_x$ be the discriminant group at $x$. Denote by $L_x \subset \cS_{\cP}$ the leaf of $\FP$  through $x$.   
  Recall that   the kernel of $\cG^x$ is the subgroup $K(\cG^x) \subset G_0$ as defined in Definition~\ref{def-kernel}, and is the isotropy subgroup of the action $(V_0,G_0,\Phi_0)$ at $x$.

 \subsection{Locally trivial germinal holonomy}
 The following properties of a pseudogroup action are basic for understanding their dynamical properties.

\begin{defn}\label{def-holonomy} 
  Given  $g_1 , g_2 \in K(\cG^x)$, we say $g_1$ and $g_2$ have the same  \emph{germinal  holonomy} at $x$ if there exists an open set $U_x \subset V_0$ with $x \in U_x$, such that the restrictions $\Phi_0(g_1)|U_x$ and $\Phi_0(g_2)|U_x$  agree  on $U_x$.  In particular, we say that $g \in K(\cG^x)$  has \emph{trivial  germinal holonomy} at $x$ if there exists an open set $U_x \subset V_0$ with $x \in U_x$, such that the restriction $\Phi_0(g)|U_x$ is the trivial map.
\end{defn}

By straightforward checking of definitions, one can see that the notion `germinal holonomy at $x$' defines an equivalence relation on the image of the isotropy subgroup $K(\cG^x)$ under the global holonomy map $\Phi_0 \colon  G_0 \to {\it Homeo}(V_0)$. Denote by ${\rm Germ}(\Phi_0 , x)$ the quotient of $\Phi_0(K(\cG^x))$ by this equivalence relation. Thus the composition of $\Phi_0 \colon K(\cG^x) \to Homeo(V_0)$   with the quotient map  gives us a surjective map $K(\cG^x) \to {\rm Germ}(\Phi_0 , x)$. A standard argument shows that if ${\rm Germ}(\Phi_0 , x)$ is trivial, and $y$ is in the same $G_0$-orbit of $x$, then ${\rm Germ}(\Phi_0 , y)$ is trivial. This leads to the following definition.

\begin{defn}\label{def-noholonomy}
We say that a leaf $L_x$ is \emph{without holonomy}, or that $L_x$ has \emph{trivial holonomy},  if ${\rm Germ}(\Phi_0 , x)$ is trivial.   We say that ${\rm Germ}(\Phi_0 , x)$ is \emph{locally trivial}, if there exists an open set $U_x \subset V_0$ with $x \in U_x$ such that for \emph{every} $g \in K(\cG_x)$ the restriction $\Phi_0(g)|U_x$ is the trivial map.  
\end{defn}

  The distinction between the holonomy group ${\rm Germ}(\Phi_0 , x)$ being trivial, and it being locally trivial, may seem technical, but this distinction is related to fundamental dynamical properties of the foliation $\FP$ of $\cS_{\cP}$. For example, it is a key concept in  the  generalizations of the Reeb stability theorem from compact leaves to the non-compact case for codimension-one foliations, as discussed in the works of  Sacksteder and Schwartz \cite{SackstederSchwartz1965} and Inaba \cite{Inaba1977,Inaba1983}. The nomenclature ``locally trivial'' was introduced by Inaba    \cite{Inaba1977,Inaba1983}. As we see below, this distinction is also important for the study of the dynamics of weak solenoids. First, we make an elementary observation, which implies  Lemma~\ref{lemma-germinalholonomy} of the Introduction.

\begin{lemma}\label{lem-loctrivial}
Suppose that $K(\cG^x)$ is a finitely-generated. If ${\rm Germ}(\Phi_0 , x)$ is trivial, then ${\rm Germ}(\Phi_0 , x)$ is locally trivial.
\end{lemma}
\proof
Let $\{g_1, \ldots , g_k\} \subset K(\cG^x)$ be a set of generators. Then ${\rm Germ}(\Phi_0 , x)$ trivial implies that for each $1 \leq i \leq k$ there exists an open $U_i \subset V_0$ with $x \in U_i$ such that the restriction of $\Phi_0(g_i)|U_i$ is the trivial map. Then let $U_x = U_1 \cap \cdots \cap U_k$ which is an open neighborhood of $x$, and the restriction $\Phi_0(g)|U_x$ is then trivial for all $g \in K(\cG^x)$.
\endproof

 We also recall a basic result,   which  is a version of the fundamental result of Epstein, Millet and Tischler \cite{EMT1977} in the language of group actions on Cantor sets. 

\begin{thm} \label{thm-emt}\cite{EMT1977}
Let $(V_0 , G_0, \Phi_0)$ be a given action, and suppose that $V_0$ is a Baire space. 
Then the  union of all $x \in V_0$ such that  ${\rm Germ}(\Phi_0 , x)$ is the trivial group forms a  $G_{\delta}$ subset of $V_0$. In particular, there exists at least one  $x \in V_0$ such that ${\rm Germ}(\Phi_0 , x)$ is the trivial group.
\end{thm}
The following is an immediate consequence of this result and Definition~\ref{def-regular}.

\begin{cor}\label{cor-emt}
Let $(V_0 , G_0, \Phi_0)$ be a \emph{regular}  equicontinuous minimal Cantor system,   then ${\rm Germ}(\Phi_0 , x)$ is the trivial group for all $x \in V_0$.
Consequently, if $\fM$ is a homogeneous matchbox manifold, then all leaves of $\FfM$ are without germinal holonomy. 
\end{cor}

 \subsection{Algebraic conditions}

Next, we  explore   the relation between the structure of a group chain $\cG^x$ and the germinal holonomy group at $x$.
First,   note that  for a given section $V_0$ and the holonomy action $(V_0,G_0,\Phi_0)$, the assumption that the germinal holonomy group ${\rm Germ}(\Phi_0 , x)$ is trivial  need not imply that $K(\cG^x)$ is trivial, or even that it is a normal subgroup of $G_0$, as the following example shows.

\begin{ex}\label{ex-product}
{\rm
Let $\G$ be a finitely presented group, and $\{\G_{\ell}\}_{\ell \geq 0}$ be a chain of normal subgroups in $\G$ with kernel  
   $\ds \G_x = \bigcap_{\ell} \G_{\ell}$. 
 Let $H$ be a finite simple group, and let $K \subset H$ be a non-trivial subgroup. Since $H$ is simple, $K$ is not normal in $H$.

Let $G_0 = H \times \G$, and $G_{\ell} = K \times \G_{\ell}$, $\ell \geq 0$. Note that $G_{\ell}$  is a normal subgroup of $G_1 = K \times \G_1$ for all $\ell \geq 1$, but $G_{\ell}$ is not normal in $G_0$. Thus, the group chain $\{G_{\ell}\}_{\ell \geq 0}$ is weakly normal. Let  $M_0$ be a compact connected manifold without boundary, such that $\pi_1(M_0,x_0) = G_0$, where $x_0 \in M_0$ is some basepoint. Then the group chain $\cG^x = \{G_{\ell}\}_{\ell \geq 0}$ yields a presentation $\cP = \{f^{\ell+1}_{\ell} \colon  M_{\ell+1} \to M_{\ell}\}$, and the corresponding solenoid $\cS_\cP$ is homogeneous by Proposition~\ref{prop-fokking}.

By Theorem~\ref{thm-emt} $\cS_\cP$ has a leaf $L_y$ without holonomy. By Remark~\ref{rmk-basepoints}, a group chain with basepoint $y$ is given by $\cG^y = \{g_iG_ig_i^{-1}\}_{i \geq 0}$, where $g_i=(c_i,\gamma_i)$. Since the projection $G_0/G_{\ell+1} \to G_0/G_{\ell}$ restricts to the identity map on the factor $H/K$, for all $\ell \geq 0$, one can write $g_i = (c, \gamma_i)$ for some $c \in H$. Since each $\G_{\ell}$ is a normal subgroup, we have that  $g_i G_i g_i^{-1} =  cKc^{-1} \times \G_i$.   Thus,  $K(\cG^y) = cKc^{-1} \times \G_x$ is not a normal subgroup of $G_0$, since $H$ is simple.
}
\end{ex}

Next, we consider the holonomy action of the elements in $K(\cG^x)$ on $V_0$ in more detail, using the inverse limit model 
${\ds \tau_x \colon V_0 \cong X_\infty^x = \{G_0/G_{\ell+1}^x \to G_0/G_{\ell}^x\} }$. For each $n \geq 0$, denote by 
\begin{equation}\label{eq-cylinderset}
U(x,n) = \{(g_{\ell} G_{\ell}^x) \in X_\infty^x \mid g_{\ell} = e \ {\rm if } \  \ell \leq n ~ ; ~ g_{\ell} G^x_{\ell} = g_{\ell +1} G^x_{\ell} ~ {\rm for ~ all} ~  \ell \geq n \}
\end{equation}
which is a ``cylinder neighborhood'' of $(e G_{\ell}^x) \in V_0$. 
Note that $\tau_x(V_{n}) = U(x,n) $ for $n \geq 0$, where $V_{n}$ is a generating set in the partition introduced in Proposition~\ref{prop-AFpres}.

Since $K(\cG^x)$ is a subgroup of $G_0$, for each $n \geq 1$ one can consider its left action on the cosets in $G_0/G_{n}^x$. Such an action fixes the coset  $eG_{n}^x$,  thus  the action of $g \in K(\cG_{n}^x)$ fixes the neighborhood of the identity as a set, $\Phi_0(g) \colon U(x,n)  \to U(x,n)$ for $g \in G_{n}^x$, and permutes the points in $U(x,n)$. 

Now observe that the action of $g$ has trivial germinal holonomy at $x$ if for some $n_g > 0$,  $g$ acts trivially on the clopen neighborhood $U(x,n_g)$ of $x$; that is, $\Phi_0(g)|U(x,n_g)$ is the trivial map. The following algebraic characterization of elements without holonomy was obtained in \cite[Lemma 5.3]{DHL2016b}.

\begin{lemma}\cite{DHL2016b}\label{lem-triv-holonomy}
The action of $g \in K(\cG^x)$ has \emph{trivial germinal holonomy} at $x$ if and only if there exists some index $i_g \geq 0$ 
such that multiplication by $g$ satisfies $g \cdot h K(\cG^x) = h K(\cG^x)$ for all $h \in G_{i_g}$.  
That is, $h^{-1} g h \in K(\cG^x)$ for all $h \in G_{i_g}$.
\end{lemma}





In the case where the kernel $K(\cG^x)$ is finitely generated, we have the following consequence of Lemma~\ref{lem-triv-holonomy}, whose proof can be compared with that of Lemma~\ref{lem-loctrivial}.

\begin{prop}\label{prop-kernel-inclusion}
Let $\cG^x=\{G_{\ell}^x\}_{\ell \geq 0}$ be a group chain, and suppose the kernel $K(\cG^x)$ is finitely generated. Suppose that ${\rm Germ}(\Phi_0 , x)$ is the trivial group, then 
  there is an index $\ell_x \geq 0$ such that $K(\cG^x)$   is a normal subgroup of  $G_{\ell_x}^x$.
\end{prop}
\proof
Let $\{g_1, \ldots , g_k\} \subset K(\cG^x)$  be a set of generators. Then for each $1 \leq \ell \leq k$, there exists $i_{\ell} \geq 0$ 
such that   $h^{-1} g h \in K(\cG^x)$ for all $h \in G_{i_{\ell}}$. Let $\ell_x = \max\{ i_1 , \ldots , i_k\}$, then this implies that 
$h^{-1} g h \in K(\cG^x)$ for all $g \in K(\cG^x)$ and $h \in G_{\ell_x}$; that is, 
$K(\cG^x)$ is a normal subgroup of  $G_{\ell_x}^x$.
\endproof

\begin{remark}\label{rmk-Schorinonbhd}
{\rm
The condition that the kernel $K(\cG^x)$ of the group chain $\cG^x$ is finitely generated is essential. Example~\ref{ex-schori}  gives   a group chain    whose kernel at $x$ is    infinitely generated, and the germinal holonomy group ${\rm Germ}(\Phi_0 , x)$ is not {locally trivial}. 
}
\end{remark}

   Proposition~\ref{prop-kernel-inclusion} implies the following result, which is an algebraic  analog of Reeb stability.

\begin{prop}\label{prop-equalkernels}
Let   $(V_0,G_0,\Phi_0)$ be a minimal equicontinuous Cantor group action. Let $x, y \in V_0$ be such that both germinal holonomy groups 
${\rm Germ}(\Phi_0 , x)$ and ${\rm Germ}(\Phi_0 , y)$ are \emph{locally trivial}. Then for   
 associated group chains $\cG^x$ and $\cG^y$, the    kernels $K(\cG^x)$ and  $K(\cG^y)$ are conjugate   subgroups of $G_0$.
\end{prop}

\proof 
  Let $\cG^x$ and $\cG^y$ be   group chains   at $x$ and $y$, respectively, for the action  $(V_0,G_0,\Phi_0)$.
  
  Let   $\tau_x \colon  \fX_0 \to X^x_\infty$  and $\tau_y \colon  \fX_0 \to X^y_\infty$ be 
  the corresponding homeomorphisms defined in   Lemma~\ref{lem-fibremodel}, each of  which is 
 equivariant with respect to the action \eqref{eqn-globalholonomy} of $G_0$.

  By the assumption that ${\rm Germ}(\Phi_0 , x)$ is locally trivial,   there exists an open set $U_x \subset V_0$ with $x \in U_x$ such that for  every $g \in K(\cG^x)$ the restriction $\Phi_0(g)|U_x$ is the trivial map. As the image $\tau_x(U_x) \subset X^x_\infty$ is open,  and contains  $(e G_i^x) = \tau_x(x)$,  there exists an index $\ell_x > 0$ such that 
    $U((e G_i^x),\ell_x)  \subset \tau_x(U_x)$, where  $U((e G_i^x),\ell_x)$  is defined in \eqref{eq-cylinderset}. 
       Note that $G_{\ell_x}^x$ is the stabilizer of $U((e G_i^x) ,\ell_x)$ for the action of $G_0$.    Then   $K(\cG^x)$  acts trivially on  $U((e G_i^x) ,\ell_x)$, so  $K(\cG^x)$   is a normal subgroup of  $G_{\ell_x}^x$ by    Lemma~\ref{lem-triv-holonomy}.

Set $V_1 = \tau_x^{-1}(U((e G_{\ell}^x),\ell_x)) \subset U_x$ and let  $z \in V_1$ with $z \ne x$. 
Then the image  $\tau_x(z) = (h_i G_{i}^x)$ where  $h_i \in G_{\ell_x}^x$ for $i \geq \ell_x$ and   $h_i =e$ for $i \leq \ell_x$.
As usual, the sequence $(h_i)$ also satisfies the compatibility condition $h_i G^x_i = h_j G^x_i$ for all $i \geq 0$ and $j > i$.
By Remark~\ref{rmk-basepoints},  we have that $\cG^z =  \{h_i G_i^x h_i^{-1}\}_{i \geq 0}$.

Note that $h_i \ K(\cG^x) \ h_i^{-1} = K(\cG^x)$ for $i \geq 0$, since $K(\cG^x)$ is normal in $G_{\ell_x}^x$, so we have   
 \begin{equation}\label{kernelinclusion}
K(\cG^x) = \bigcap_{i \geq 0} ~ G_{i}^x  = \bigcap_{i \geq 0} ~ h_i \  K(\cG^x) \ h_i^{-1}  ~ \subseteq  ~  \bigcap_{i \geq 0} ~ h_i \ G_i^x  \ h_i^{-1}  =   K(\cG^z) \ .
\end{equation}
In general, this inclusion may be proper, as illustrated in Example~\ref{ex-klein}.

Now assume that ${\rm Germ}(\Phi_0 , z)$ is locally trivial. We show that $K(\cG^z) \subseteq K(\cG^x)$.
First, note that
   there exists an open set $U_z \subset V_0$ with $z \in U_z$ such that for  every $g \in K(\cG^z)$ the restriction $\Phi_0(g)|U_z$ is the trivial map. 
Recall that  $\tau_x(z) = (h_i G_i^x) \in U((e G_{\ell}^x),\ell_x)$. Then there exists  
  $\ell_z \geq \ell_x$   such that   
 \begin{equation}\label{eq-ballatz}
 U((h_i G_i^x),\ell_z)  = \{(g_i G_{i}^x) \in X_\infty^x \mid g_i = h_i \ {\rm for}\  i \leq \ell  \} \subset \tau_x(U_z) \ ,
 \end{equation}
That is, $g \in K(\cG^z)$ acts trivially on the cylinder set $U((h_iG_i^x),\ell_z)$ in $X_\infty^x$. Let $h = h_{\ell_z} \in  G_{\ell_x}^x$, so we obtain an element $(h   G_i^x) \in X^x_\infty$. By choice of $h$ and \eqref{eq-ballatz} we have $(h   G_i^x) \in U((h_i G_i^x),\ell_z)$.  Now let  $g \in K(\cG^z)$, then  the restricted map $\Phi_0(g)|U_z$ is the identity, so 
we have  $g \cdot (h   G_i^x) =  (h   G_i^x)$. But this means that  $ h^{-1} g h G_i^x  =  G_i^x$  for all $i \geq 0$, and thus  $ h^{-1} g h   \in K(\cG^x)$, or $g \in h K(\cG^x) h^{-1}$. Since $h \in G_{\ell_z}^x \subset G_{\ell_x}^x$, and $K(\cG^x)$ is a normal subgroup of $G_{\ell_x}^x$, this implies that $K(\cG^z) \subseteq K(\cG^x)$. 

 Now suppose that $y \in V_0$ is such that ${\rm Germ}(\Phi_0 , y)$ is locally trivial. The action of $G_0$ on $V_0$ is assumed to be minimal, so there exists $g \in G_0$ such that $z = \Phi_0(g)(y) \in V_1$. Then the holonomy at $z$ is also locally trivial, so $K(\cG^z) = K(\cG^x)$ by the argument above. On the other hand, we have $K(\cG^y) = g^{-1} K(\cG^z) g$ as $K(\cG^y)$ is the isotropy subgroup of $y$. The claim of the proposition then follows.
\endproof

  \subsection{Kernels and discriminants}

We give two results concerning  the relation between the kernel  of a  group chain  and its discriminant.   

\begin{prop}\label{prop-kernel-discriminant}
Let  $(V_0,G_0,\Phi_0)$ be an equicontinuous minimal Cantor system, $x \in V_0$ be a choice of a basepoint, and $\cG^x =   \{G_{\ell}^x\}_{\ell \geq 0}$ be a group chain associated to $(V_0,G_0,\Phi_0)$ at $x$.
  Let  $\cL_0 = ker(\Phi_0)$ denote the kernel of   $\Phi_0 \colon  G_0 \to {\it Homeo}(V_0)$. 
 Then $K(\cG^x) \subset \cL_0$   if and only if the intersection $\Phi_0(G_0) \cap \overline{\Phi_0(G_0)}_x$ is the trivial group.
\end{prop}

\proof  

By Theorem~\ref{thm-quotientspace}, we can identify $\overline{\Phi_0(G_0)} \cong C_{\infty}$ and $\overline{\Phi_0(G_0)}_x \cong \cD_x$, where the image $\Phi_0(G_0)$ is identified with the elements $(g_{\ell} C_{\ell})   \in C_{\infty}$ such that  $g_{\ell} C_{\ell} = g C_{\ell}$ for all $\ell \geq 0$, for some $g \in G_0$. 

First, suppose that $g \in G_0$ satisfies $\Phi(g) \in \overline{\Phi_0(G_0)}_x$ and $\Phi(g)$ is not the trivial element. 
Then  $\whg = (g C_{\ell}) \in \cD_x$ and $(g C_{\ell}) \ne (eC_{\ell})$, so that   there exists $\ell_0 > 0$ such that $g \not\in C_{\ell_0}$. 
By the definition of $\cD_x$ in \eqref{eq-discriminantdef}, we have that $\widehat{g} $ is in the image of the map
$\ds \delta^{\ell+1}_{\ell} \colon   D_{\ell+1}^x \to D_{\ell}^x$ for all $\ell > 0$ where $D_{\ell}^x = G_{\ell}^x/C_{\ell}$. 
This implies that $ g C_{\ell}   \subset   G_{\ell}^x$, and hence $g \in G_{\ell}^x$ for all $\ell \geq 0$, and so $g \in K(\cG^x)$. We claim that $\Phi_0(g)$ is not the trivial action, so that 
$g \not\in \cL_0$. It is given that $g \not\in C_{\ell_0}$, hence $g C_{\ell_0} \ne C_{\ell_0}$.
 Then  for all $\ell \geq \ell_0$, we have $g  C_{\ell} \ne C_{\ell}$,  so $g \cdot  (eC_{\ell}) \ne (eC_{\ell})$, which implies that $g \not\in \cL_0$.  
It follows that $K(\cG^x)  \not\subset \cL_0$ as was to be shown. 
 
 Conversely, let $g \in K(\cG^x)$ and suppose that $g \not\in \cL_0$.  First note that  $g \in G_{\ell}^x$ for all $\ell \geq 0$, and so we have   $\whg = (g C_{\ell}) \in \cD_x$. 
 The assumption that $g \not\in \cL_0$ implies there exists  some $(h_{\ell} C_{\ell}) \in C_{\infty}$ such that $g \cdot (h_{\ell} C_{\ell})  \ne (h_{\ell} C_{\ell})$. 
 Thus,    there exists $\ell_0 > 0$ such that for   all $\ell \geq \ell_0$ we have $g h_{\ell} C_{\ell} \ne h_{\ell} C_{\ell}$, which implies that 
 $h_{\ell}^{-1} g h_{\ell} \not\in C_{\ell}$ and so $g \not\in C_{\ell}$ as $C_{\ell}$ is a normal subgroup of $G_0$. Thus, $(e C_{\ell}) \ne (g C_{\ell})$ for all $\ell \geq  \ell_0$, and so $(g C_{\ell}) \in \cD_x$ is non-trivial. That is, $\Phi(g) \in \overline{\Phi_0(G_0)}_x$ is a non-trivial element, as was to be shown.
\endproof
 
Compare the following   application of  Proposition~\ref{prop-kernel-discriminant}   with the conclusions of Theorem~\ref{thm-sections}.

\begin{prop}\label{prop-holodisc}
  Let $\fM$ be an equicontinuous matchbox manifold, 
     let  $V_0$ be a transverse section in $\fM$ as given in Proposition~\ref{prop-AFpres}, and let $x \in V_0$ be a choice of a basepoint.  
     Let $V_{\ell}$ be defined as in Proposition~\ref{prop-AFpres}, so that $x \in V_{\ell}$ for all $\ell \geq 0$.  
 Let  $G_0$ be the restricted holonomy group acting on $V_0$, and let $G_{\ell}^x \subset G_0$ be the stabilizer group of the set $V_{\ell}$. 
  Let $\cG^x = \{G_{\ell}^x\}_{\ell \geq 0}$ be the associated group chain in $G_0 = \pi_1(M_0, x_0)$, 
  let  $\cG^x_n = \{G_{\ell}^x\}_{\ell \geq n}$ be the  associated truncated  group chain. Assume that the leaf $L_x$ containing $x$ has non-trivial germinal holonomy, then the discriminant $\cD_n^x$ for the chain  $\cG^x_n$ is non-trivial, for all $n \geq 0$.
 \end{prop}
 \proof
 Let $n \geq 0$, and let $\cL_n \subset G_n^x$ be the kernel of the restricted action $\Phi_n \colon G_n^x \to Homeo(V_n)$.
 Then by Lemma~\ref{lem-triv-holonomy}, the kernel $K(\cG_n^x) \subset G_n^x$ is not a normal subgroup, so $\cL_n \subset K(\cG_n^x)$ is a proper inclusion. 
 Then by Proposition~\ref{prop-kernel-discriminant}, the discriminant group $\cG^x_n$ is non-trivial also.  
  \endproof
  
This yields the proof of Theorem~\ref{thm-nontrivkernel} of the Introduction, which we restate now.
 
\begin{thm}\label{cor-holodisc}
  Let $\fM$ be an equicontinuous matchbox manifold. If there exists a leaf with non-trivial holonomy for $\FfM$, then for any choice of transversal $V_0 \subset \fM$, the resulting Molino sequence \eqref{eq-molinoseqM2}  has non-trivial fiber $\cD$.
\end{thm}

  The converse to Theorem~\ref{cor-holodisc} is not true. Fokkink and Oversteegen \cite[Theorem 35]{FO2002} constructed an example of a solenoid with simply connected leaves which is non-homogeneous. Since the leaves are simply connected, they have trivial holonomy.   In Section~\ref{sec-construction} we construct further examples of actions with non-trivial Molino fibre and simply connected leaves.

\section{Strongly quasi-analytic actions}\label{sec-analytic}
 In this section, we study the condition of \emph{strong quasi-analyticity}, abbreviated as the \emph{SQA}  condition, for equicontinuous matchbox manifolds, as defined in Definition~\ref{strongly-qa} below. We  identify classes of matchbox manifolds for which this condition holds, and also give examples for which it does not.   
 The generalization of Molino theory   in \cite{ALM2016} applies  to equicontinuous foliated spaces  such that the closure of their holonomy $\psg$s satisfies the \emph{SQA} condition. Thus,   it is important to characterize the weak solenoids with this property.

\subsection{The strong quasi-analyticity condition}\label{subsec-defSQA}
  The precise notion of the \emph{SQA}   condition has evolved in the literature,   motivated by the search for a condition equivalent to quasi-analyticity condition for the $\psg$s of  smooth foliations as introduced by Haefliger \cite{Haefliger1985}. 
  {\'A}lvarez L\'opez and Candel introduced the notion of a \emph{quasi-effective} $\psg$ in the work \cite{ALC2009}  as part of their study of equicontinuous foliated spaces.  
    This terminology was replaced by the notion of a  \emph{strongly quasi-analytic} $\psg$ in the work \cite{ALM2016} by {\'A}lvarez L\'opez and Moreira Galicia.

   \begin{defn}\cite{Haefliger1985}\label{quasi-analytic}
 A $\psg$ $\cG^*$  acting on a locally compact locally connected space $\fT$ is \emph{quasi-analytic}, if for every $h \in \cG^*$ the following holds: let $U \subset {\rm Dom}(h) \subset \fT$ be an open set, and suppose $x \in \fT$ is in the closure of $U$. Suppose the restriction $h|U$ is the identity map. Then there is an open neighborhood $V$ of $x$, such that the restriction $h|V$ is the identity map.
 \end{defn}
 
 Definition~\ref{quasi-analytic} describes the properties of $\psg$s, which were discussed in Remark~\ref{rmk-pseudostar}, where the action of an element is \emph{locally determined}; that is, if $h$ is the identity on an open set, then it is the identity on a larger set. For  the case where the   space $\fT$ is   not locally connected, {\'A}lvarez L\'opez and Candel \cite{ALC2009} introduced the following modification of this notion.

 \begin{defn} \label{strongly-qa}
 A $\psg$ $\cG^*$ acting on a locally compact space $\fT$ is \emph{strongly quasi-analytic}, or \emph{SQA}, if for every $h \in \cG^*$ the following holds: let $U \subset {\rm Dom}(h)$ be a non-empty open set, and suppose the restriction $h|U$ is the identity map. Then $h$ is the identity map on its domain ${\rm Dom}(h)$.
 A matchbox manifold $\fM$ satisfies the \emph{SQA} condition if there exists a traversal $V_0 \subset \fM$ such that the induced $\psg$ $\cGF^*$ on $V_0$ satisfies the \emph{SQA} condition.
 \end{defn}

 Definition~\ref{strongly-qa} says that the action of an equicontinuous strongly quasi-analytic $\psg$ $\cG$ is locally determined. That is, if $h$ is the identity on a non-empty open subset of its domain, then it is the identity on ${\rm Dom}(h)$.  In the case where the transversal $\fT$ is locally compact and locally connected, this condition is  equivalent to quasi-analyticity by \cite[Lemma 9.8]{ALC2009}. However,   when $\fT$ is totally disconnected, the \emph{SQA}  condition becomes a statement about the algebraic properties of the group chain associated to the action, as we next discuss.

 Recall from   Proposition~\ref{prop-AFpres}(1) that if $\fM$ is an equicontinuous matchbox manifold, then we can assume that the $\psg$ action on the transversal is given by an equicontinuous minimal Cantor   action $(V_0,G_0,\Phi_0)$. Thus, for each $h \in G_0$ we have ${\rm Dom}(h) = V_0$. Moreover, the assumption that the restriction $h|U$ is the identity in the statement of Definition~\ref{strongly-qa},  means that the \emph{SQA}  condition needs only be checked for   $h \in G_0$ such that there exists $x \in V_0$ for which $\Phi_0(h)(x) = x$, that is, those elements whose action fixes at least a point.
 
Recall from  Section~\ref{subsec-ellis} that  the closure $\overline{\Phi_0(G_0)} \subset Homeo(V_0)$ in the uniform topology of the image $\Phi_0(G_0) \subset Homeo(V_0)$ is called the Ellis group of the Cantor system $(V_0, G_0, \Phi_0)$, which yields a Cantor system $(V_0, \overline{\Phi_0(G_0)}, \whPhi_0)$, where $\whPhi_0 \colon \overline{\Phi_0(G_0)} \to Homeo(V_0)$.   Given $x \in V_0$ then
$\overline{\Phi_0(G_0)}_x \subset \overline{\Phi_0(G_0)}$ denotes the isotropy subgroup at $x$ for the action, and then the \emph{SQA}  condition must be checked for all elements of  $\overline{\Phi_0(G_0)}_x$. We denote by $\whPhi_0(G_0) = \{\Phi_0(g) \mid g \in G_0\}$ which is a dense subgroup of $\overline{\Phi_0(G_0)}$.
The following result follows from the definitions.

\begin{lemma}\label{lem-closurepseudogroup-sqa}
If $(V_0,\overline{\Phi_0(G_0)},\whPhi_0)$ satisfies the \emph{SQA}  condition, then $(V_0, G_0, \Phi_0)$   also satisfies the \emph{SQA}  condition. Conversely, suppose that $\overline{\Phi_0(G_0)}_x \subset \whPhi_0(G_0)$, then 
$(V_0, G_0, \Phi_0)$ satisfies the \emph{SQA}  condition implies that $(V_0,\overline{\Phi_0(G_0)},\whPhi_0)$ satisfies the  \emph{SQA}  condition.
\end{lemma}
 
\proof 
Let $g \in G_0$ and set $\whg = \Phi_0(g) \in Homeo(V_0)$. Then $\whg \in \overline{\Phi_0(G_0)}$, so that if $(V_0,\overline{\Phi_0(G_0)},\whPhi_0)$ satisfies the \emph{SQA}  condition then so must the action of $\whg$. 
Conversely,   suppose $(V_0,G_0,\Phi_0)$ satisfies the \emph{SQA}  condition.  As noted above, the \emph{SQA}  property need only be checked for $h \in \overline{\Phi_0(G_0)}_x$. By assumption, such $h \in \Phi_0(G_0)$ and so satisfies the \emph{SQA}  condition.
 \endproof

 Note that the assumption that $\overline{\Phi_0(G_0)}_x \subset \whPhi_0(G_0)$ implies that the compact set $\overline{\Phi_0(G_0)}_x$ is contained in a countable set, hence it must be finite. Thus, by Theorem~\ref{thm-quotientspace}, this implies that the discriminant group $\cD_x$ of the action is finite. The converse need not be true, that is, if the discriminant $ \overline{\Phi_0(G_0)}_x $ is finite, then it may be possible to choose a point $y \in V_0$, such that $\overline{\Phi_0(G_0)}_y $ has trivial intersection with $\whPhi_0(G_0)$, for instance, this is the case for Example~\ref{ex-klein}. Examples in Section~\ref{sec-construction} show that it is  possible to construct   actions $(V_0,G_0,\Phi_0)$ such that  $\overline{\Phi_0(G_0)}_x $ has trivial intersection with $\whPhi_0(G_0)$ for any choice of $x \in V_0$.
 
We next consider the \emph{SQA}  property for an equicontinuous minimal Cantor system $(V_0 , G_0 , \Phi_0)$ and its associated Ellis system $(V_0,\overline{\Phi_0(G_0)},\whPhi_0)$. This condition for the system $(V_0 , G_0 , \Phi_0)$ can be formulated  in terms of the group chain model developed in  Sections~\ref{sec-inverselimitmodel} and~\ref{subsec-ellischains}, in which case Lemma~\ref{lem-triv-holonomy} and Proposition~\ref{prop-kernel-inclusion} imply that the condition is a statement about the holonomy action of the kernel $K(\cG^x)$ of the chain $\cG^x$ for each $x \in V_0$. Examples~\ref{ex-klein} and \ref{ex-schori} below  and the discussion in Section~\ref{sec-construction} illustrate the possibilities.

The \emph{SQA}  property for  the system   $(V_0,\overline{\Phi_0(G_0)}, \whPhi_0)$ can be much more subtle to check, as now it is a condition on the action of the isotropy group $\overline{\Phi_0(G_0)}_x \cong \cD_x$ which depends on the algebraic properties of the closed subgroup $\cD_x \subset C_{\infty}$. Note that in this case, for any $x, y \in V_0$ the isotropy groups $\cD_x$ and $\cD_y$ are conjugate in $C_{\infty}$, so it suffices to consider the condition for a fixed choice of basepoint $x \in V_0$.

 \subsection{Sufficient conditions for the \emph{SQA} property}
 
We next indicate a few classes of solenoids which satisfy the quasi-analyticity condition.
 
\begin{lemma}\label{lemma-homogeneous-qa}
If a matchbox manifold $\fM$ is homogeneous, then there exists a section $V_0$ with associated presentation $\cP$, such that the actions $(V_0,G_0,\Phi_0)$ and $(V_0,\overline{\Phi_0(G_0)}, \whPhi_0)$ are \emph{SQA} . 
\end{lemma}

\proof By Corollary~\ref{cor-homogeneoussolenoid} one can assume that $V_0$ and $\cP$ are chosen so that the associated group chain $\{G_{\ell}^x\}_{\ell \geq 0}$ consists of normal subgroups. Then    $K(\cG^x)$ is a normal subgroup of $G_0$, so by   Lemma~\ref{lem-triv-holonomy}, each $g \in K(\cG^x)$ defines a trivial holonomy action on $V_0$. Hence    the action of $G_0$ on $V_0$ is \emph{SQA} . 

Since $\{G_{\ell}^x\}_{\ell \geq 0}$  is a chain of normal subgroups, then by Proposition~\ref{prop-regularaction} the isotropy group $\overline{\Phi_0(G_0)}_x$ is trivial, and so the condition $\overline{\Phi_0(G_0)}_x \subset \whPhi_0(G_0)$ is trivially satisfied. Then by Lemma~\ref{lem-closurepseudogroup-sqa} the action $(V_0,\overline{\Phi_0(G_0)}, \whPhi_0)$ is \emph{SQA} .
\endproof

Note that the holonomy pseudogroups associated to homogeneous solenoids, as in  Lemma~\ref{lemma-homogeneous-qa}, satisfy a stronger condition than \emph{SQA}. 
Recall from \cite[Definition 2.22]{ALM2016}, that the action of $G_0$ on $V_0$ is \emph{strongly locally free} if for all $h \in G_0$, if $h(x) = x$, then $h(y) = y$ for all $y \in V_0$. If $\fM$ is homogeneous, then the action on a local section $V_0$, as given by Lemma~\ref{lemma-homogeneous-qa}, is strongly locally free. The actions in Lemma~\ref{lemma-homogeneous-qa} are the actions in \cite[Example 2.35]{ALM2016}.
 
 The following result gives a  class of equicontinuous matchbox manifolds which satisfy the \emph{SQA} condition. This 
 is Theorem~\ref{thm-sqaactions} of the Introduction.

\begin{thm}\label{thm-sqa-finitediscr}
Let $\fM$ be an equicontinuous matchbox manifold of finite $\pi_1$-type. Then there exists a section $V_0$ with a presentation $\cP$, such that the action $(V_0,G_0,\Phi_0)$ is \emph{SQA} . In addition, if $V_0$ can be chosen so that the discriminant group $\cD_x = \overline{\Phi_0(G_0)}_x$ is finite, then there exists an $n \geq 0$ such that the restricted action $(V_n,G_n^x,\Phi_n)$ and the action $(V_n,\overline{\Phi_n(G_n^x)},\widehat{\Phi}_n)$ are both \emph{SQA}.
\end{thm}

\proof 

 Let   $V_0$ be a transverse section in $\fM$ as given in Proposition~\ref{prop-AFpres}, and let $x \in V_0$ be a choice of a basepoint.  
By Theorem~\ref{thm-emt} we can assume that $x$ is chosen so that $L_x$ is a leaf without holonomy. As the leaves of $\FfM$ are assumed to have finite $\pi_1$-type, 
by Lemma~\ref{lem-triv-holonomy}  and Proposition~\ref{prop-kernel-inclusion}, and   restricting to a smaller section is necessary, we can assume that $V_0$ and $\{\cG_{\ell}^x\}_{\ell \geq 0}$ are chosen so that $K(\cG^x)$ is a normal subgroup of $G_0$. Then by Proposition~\ref{prop-equalkernels},  $K(\cG^x) \subseteq K(\cG^y)$ for all $y \in V_0$, and, if ${\rm Germ}(y,\Phi_0)$ is trivial, then $K(\cG^x) = K(\cG^y)$.

Since the $G_0$-orbit of $x$ is dense in $V_0$, any $g \in G_0$ which is the identity on a non-empty open set in $V_0$ must be contained in $K(\cG^x)$, and so it is the identity on $V_0$. Thus,    $(V_0,G_0,\Phi_0)$ is \emph{SQA}.

Now let $C_\infty$ be the Ellis group, associated to $(V_0,G_0,\Phi_0)$, and suppose the discriminant group $\cD_x \cong  \overline{\Phi_0(G_0)}_x $ is finite. Suppose there exists a non-trivial element $\whg \in  \overline{\Phi_0(G_0)}_x$ which fixes an open subset $U$ of $V_0$ around $x$. 

    Let $V_{\ell}$ be defined as in Proposition~\ref{prop-AFpres}, so that $x \in V_{\ell}$ for all $\ell \geq 0$.  
Choose an index $n \geq 0$ large enough so that $V_n \subset U$. Let $y \in V_n$, then $\whg = (g_iC_i) \in  \overline{\Phi_0(G_0)}_y$, and it follows that
  $$\whg \in \bigcap_{y \in V_n}  \overline{\Phi_0(G_0)}_y,$$
that is, the intersection $\bigcap_{y \in V_n} \cD_y$ is non-trivial.  

Consider the truncated chain $\{G_\ell^x\}_{\ell \geq n}$, and the corresponding action $(V_n,G_n^x,\Phi_n)$. Recall from Section~\ref{subsec-stability} that  $E_\ell^n = {\rm core}_{G_n^x} G_\ell^x$ is a maximal normal subgroup of $G_\ell^x$ in $G_n^x$, and there is an inclusion
   \begin{align}\label{eq-cosetincl}C_\ell \subset E_\ell^n \subset G_\ell^x,\end{align}
 where $C_\ell$ is the maximal normal subgroup of $G_\ell^x$ in $G_0$.  The Ellis group $E^n_\infty$ of the restricted action $(V_n,G_n^x,\Phi_n)$ is defined by \eqref{Dinfty-definen} as the inverse limit of coset spaces $G_n^x/E_\ell^n$. The inclusions \eqref{eq-cosetincl} yield a commutative diagram
   \begin{align}\label{diag-triangle}  \xymatrix{ G_n^x/C_\ell \ar[rr]^{\phi_{n,\ell}} \ar[rd] & & G_n^x/E^n_\ell \ar[ld] \\ & G_n^x/G_\ell^x& }\end{align}
 which is equivariant with respect to the natural action of $G_n^x$ on its coset spaces. Taking the inverse limits, we obtain the commutative diagram
  \begin{align}\label{diag-triangle-limit}  \xymatrix{ C_\infty^n \ar[rr]^{\phi_{n,\infty}} \ar[rd] & & E^n_\infty \ar[ld] \\ & G_\infty^n \cong V_n& }\end{align}
 where $C^n_\infty$ is the profinite subgroup of $C_\infty$, defined by \eqref{Ginfty-definen}, which is again equivariant with respect to the action of $G_n^x$ on the inverse limits, and $\phi_{n,\infty}$ is a surjective group homomorphism.
 
Let $\whg_n = \phi_{n,\infty}(\whg)$. We will show that $\whg_n$ acts trivially on $V_n$. Indeed, let $\whg = (g_\ell C_\ell)$, where $g_\ell \in G_n^x$. Then $\whg_n = (g_\ell E^n_\ell)$ for $\ell \geq n$. Since $C_\ell$ and $E_\ell^n$ are normal subgroups of $G_n^x$, the actions of $g_\ell C_\ell$ and $g_\ell E_\ell^n$ on $G_n^x/G_\ell$ are well-defined, for example, for any $h \in G_n^x$ we have
  $$g_\ell C_\ell hG_\ell^x = g_\ell hC_\ell h^{-1} hG_\ell^x = gh G_\ell^x,$$
and similarly for $g_\ell E_\ell^n$.   Since   Diagram~\ref{diag-triangle} is a commutative diagram of equivariant maps, we obtain that 
  $$g_\ell C_\ell h G_\ell^x = h G_\ell^x \ \implies \   g_\ell E^n_\ell h G_\ell^x = h G_\ell^x,$$
and it follows that if $\whg$ acts trivially on $y = (h_iG_\ell^x) \in V_n$, then $\whg_n$ acts trivially on $y$ as well.  

Then by an argument similar to the one at the beginning of this proof, we obtain that $\whg_n \in \bigcap_{y \in V_n} \cD_y^n$, where $\cD_y^n$ is the discriminant group of the truncated action $(V_0,G_n^x,\Phi_n)$ at $y \in V_n$. We note that $\bigcap_{y \in V_n} \cD_y^n$ is the maximal normal subgroup of $\cD_x^n$, and so by Proposition~\ref{prop-discrisnotnormal}  it must be trivial. Therefore, $\whg_n = \phi_{n,\infty}(\whg)$ is the identity in $E^n_\infty$.

We note that the restricted group action $(V_n,G_n^x,\Phi_n)$ is \emph{SQA}  since $(V_0,G_0,\Phi_0)$ is \emph{SQA} . By restricting to a smaller section and applying the above argument a finite number of times we may assume that no element of the discriminant group $\cD_x^n$ fixes an open subset of $V_n$. It follows that the action $(V_n,\overline{\Phi_n(G_n^x)},\widehat{\Phi}_n)$ of the closure is \emph{SQA} .
\endproof

 \subsection{\emph{SQA} counter-examples}
 
 We   give two classes of   examples to illustrate the above results.

\begin{ex}\label{ex-klein}
{\rm

We first give an example of a group action, corresponding to the holonomy of a solenoid with leaves of finite $\pi_1$-type, that is not strongly locally free.

Let $K$ be the Klein bottle, with the fundamental group $G_0 = \langle a,b \mid bab^{-1}= a^{-1}\rangle$, and let ${\ds K_\infty = \lim_{\longleftarrow}\{p \colon  K \to K\}}$ be the inverse limit space, as described in Example~\ref{ex-RT}. The solenoid $K_\infty$ contains one non-orientable leaf with one end, and every other leaf is an open two-ended cylinder. Thus,  each leaf is homotopic to a circle, and thus has  finite $\pi_1$-type.

The group chain $\cG^x$, associated to the choice of a basepoint as in Example~\ref{ex-RT} consists of subgroups $G_{\ell}^x = \langle a^{2^\ell},b \rangle$, $\ell \geq 0$, and $K(\cG^x) = \langle b\rangle$. This leaf has non-trivial holonomy, with ${\rm Germ}(x,\Phi) \cong \mZ_2$. Fokkink and Oversteegen \cite{FO2002} computed, that the kernel of a group chain based at any point which is not in the orbit of $x$, is $K(\cG^y) = \langle b^2 \rangle$, which is easily seen to be a normal subgroup of $G_0$. Thus for the chosen section $V_0$, for every point $y$ with trivial ${\rm Germ}(y,\Phi)$ the kernel $K(\cG^y)$ is a normal subgroup of $G_{\ell}^y$, $\ell \geq 0$, and the section satisfies Proposition~\ref{prop-equalkernels}. So the action $(V_0,G_0,\Phi)$ satisfies the \emph{SQA}  condition.

This action is not strongly locally free. Indeed, the action of the element $b$ fixes $x$, but it does not fix any $y$ with trivial ${\rm Germ}(y,\Phi)$. The non-trivial element in $\overline{\Phi(G_0)}_x$ acts non-trivially on any open subset of $V_0$, and so the action $(V_0,\overline{\Phi(G_0)}_x,\whPhi)$ satisfies the \emph{SQA}  condition.

}
\end{ex}

\begin{ex}\label{ex-schori}
{\rm
 We next give an example of a solenoid, for which the action of the holonomy group on the fibre is not \emph{SQA}  for any choice of a transverse section $V_0$.
 This example is the Schori solenoid \cite{Schori1966}. We now recall its construction, as described in \cite{CFL2014}.
 
 Let $X_0$ be a genus $2$ surface. Recall that a $1$-handle is a $2$-torus without an open disc, and note that the genus two surface $X_0$ can be seen as the union of two $1$-handles $H_0$ and $F_0$ intersecting along the boundaries of the open discs taken out. Let $x_0$ be a point in the intersection of the handles. Recall that the fundamental group of the genus $2$ surface can be presented as
   $$\pi_1(X_0,x_0) = \langle a,b,c,d \mid  aca^{-1}c^{-1}bdb^{-1}d^{-1}=1\rangle ,$$
 where $a$ and $b$ are longitudinal loops in $X_0$.  
 
Cut the handle $H_0$ (resp. $F_0$) along a closed curve $C_0$ (resp. $D_0$), as shown in Figure~\ref{fig:M0-pic}, a). Pull the cut handles apart to obtain the surface with boundary $\overline{X}_0$, see Figure~\ref{fig:M0-pic}, b). Take three copies of $\overline{X}_0$, denoted by $\overline{X}_0^1,\overline{X}_0^2,\overline{X}_0^3$, and identify their boundaries as shown in Figure~\ref{fig:M0-pic}, c). The resulting surface $X_1$, see Figure~\ref{fig:M0-pic}, d), has genus $4$, and there is an obvious $3$-to-$1$ covering map $f^1_0 \colon X_1 \to X_0$. Let $x_1$ be the preimage of $x_0$ in the second copy of the handle. We note that the covering $f^1_0$ is not regular; that is, the image $(f^1_0)_*\pi_1(X_1,x_1)$ of the fundamental group of $X_1$ is not a normal subgroup of $\pi_1(X_0,x_0)$. Geometrically, we can see that $f^1_0$ is irregular as follows: take a longitudinal loop $\gamma$ in $X_0$, which represents an equivalence class of loops in $\pi_1(X_0,x_0)$. The fibre of $f^1_0$ consists of three points, and we see from Figure~\ref{fig:M0-pic}, c), that depending on the initial point of the lift, $\gamma$ may lift to a loop or to a non-closed curve \cite{Schori1966}.

\begin{figure}[!htbp]
\centering
\includegraphics [width=130mm] {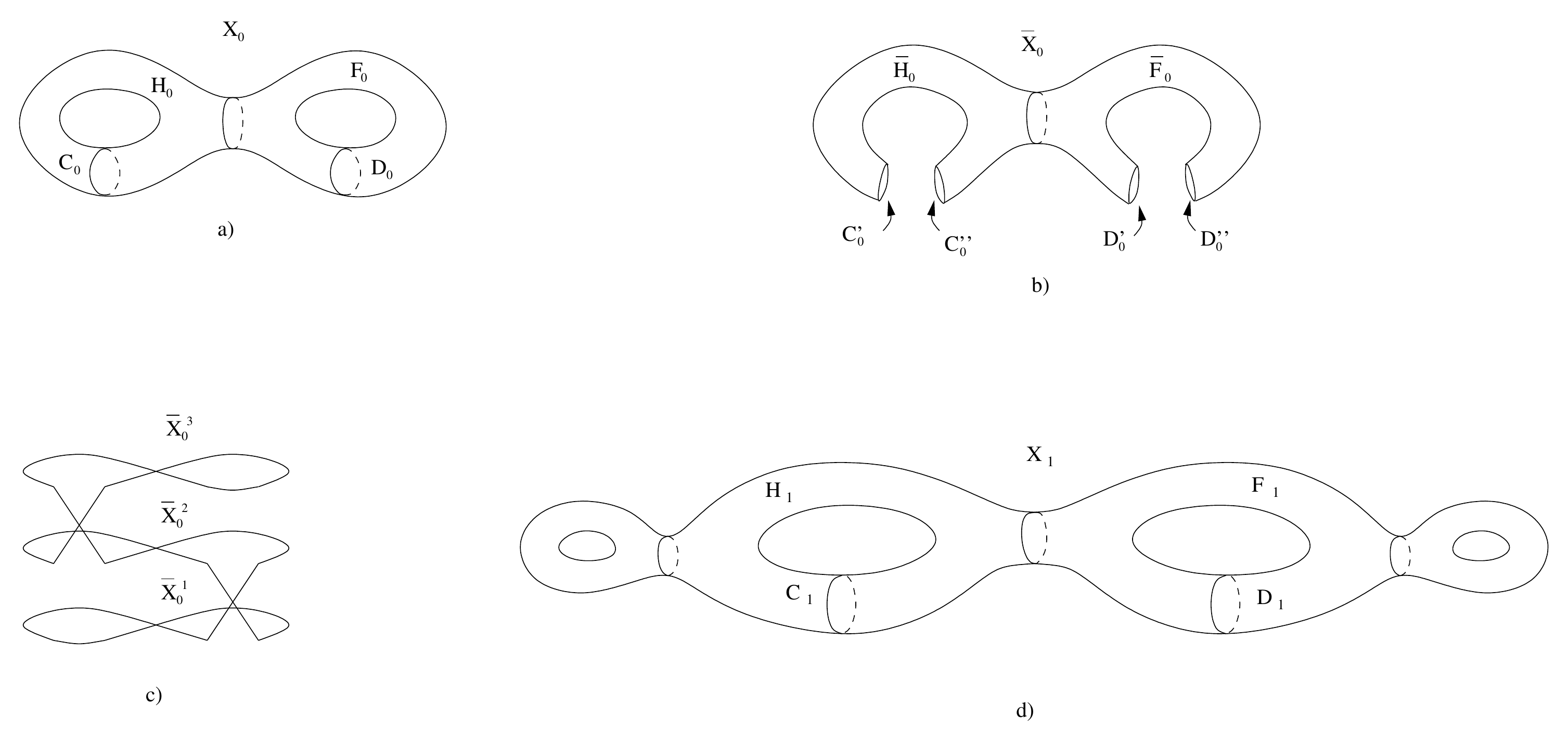}
\caption{Construction of the Schori example: a) Choice of the handles $H_0$ and $F_0$ and closed curves $C_0$ and $D_0$ in $X_0$, b) The cut surface $\overline{X}_0$, c) Identifications between $\overline{X}^\ell_0$, $\ell=1,2,3$. Each $\overline{X}^i_0$ is represented by a cut copy of figure $8$, and identifications are depicted with straight lines, d) The surface $X_1$ and the choice of the handles $H_1$ and $F_1$ and closed curves $C_1$ and $D_1$.}
 \label{fig:M0-pic}
\end{figure}

Then proceed inductively to obtain a collection of $3$-to-$1$ coverings $f^{\ell+1}_{\ell} \colon X_{\ell+1} \to X_{\ell}$. That is, we can see $X_{\ell}$ as the union of two handles $H_{\ell}$ and $F_{\ell}$, intersecting along their boundaries, see Figure~\ref{fig:M0-pic}, d), for $\ell=1$.  We cut the handle $H_{\ell}$ (resp. $F_{\ell}$) along a closed curve $C_{\ell}$ (resp. $D_{\ell}$), pull the handles apart to obtain the surface with boundary $\overline{X}_{\ell}$, take three copies of $\overline{X}_{\ell}$, denoted by $\overline{X}_{\ell}^1,\overline{X}_{\ell}^2,\overline{X}_{\ell}^3$, and identify their boundaries in a way similar to Figure~\ref{fig:M0-pic}, c). The resulting surface $X_{\ell+1}$ is a $3$-to-$1$ non-regular cover of $X_{\ell}$. This defines a presentation $\cP = \{f^{\ell+1}_{\ell} \colon X_{\ell+1} \to X_{\ell}, \ell \geq 0\}$ of the Schori solenoid $\cS_\cP$. Let $\fX_0$ be the fibre of $\cS_\cP$ at $x_0$.

For each $\ell \geq 0$, we choose $x_{\ell+1}$ to be a preimage of $x_{\ell}$ under the covering map $f^{\ell+1}_{\ell}$ in the second copy of $X_{\ell}$. Denote by $\cG^x = \{G_{\ell}^x\}_{\ell \geq 0}$ the corresponding group chain, and recall that there is a conjugacy $\phi \colon \fX_0 \to {\ds X_\infty^x = \lim_{\longleftarrow}\{G_0/G^x_{\ell+1} \to G_0/G^x_{\ell}, \ell \geq 0\}}$. As before, we denote by $U(x,\ell)$ the cylinder set in $X^x_\infty$ containing $(eG_{\ell}^x)$. If $y \in \fX_0$ is a point with $\phi(y) = (g_{\ell} G_{\ell}^x)$, then $g_{\ell} \cdot U(x,\ell) = \Phi(g_{\ell})(U(x,\ell))$ is a cylinder set containing $\phi(y)$. Denote $U_{\ell}^y = \phi^{-1}(g_{\ell} \cdot U(x,\ell))$. The group chain $\cG^y = \{G_{\ell}^y = g_{\ell} G_{\ell}^x g_{\ell}^{-1}\}_{\ell \geq 0}$ corresponds to a presentation $\cP'$ of the Schori solenoid with basepoint $y$.
 
The following theorem is Theorem~\ref{thm-schorinotsqa} of the Introduction.
 \begin{thm}\label{Schori-notqa}
In the Schori solenoid, for any choice of a basepoint $y \in \fX_0$, and any choice of a section $U^y_n$, $n \geq 0$, the holonomy action $(U^y_n,G_n^y,\Phi_n)$ is not \emph{SQA} . 
 \end{thm}
 
 \proof Let $y \in \fX_0$, and let $(U^y_n,G_n^y,\Phi_n)$ be the holonomy action. At the end of Section~\ref{sec-foliatedspaces} we described the procedure of restricting to a smaller section, which gives us a presentation $\cP_n' = \{f^{\ell+1}_{\ell} \colon  X_{\ell+1} \to X_{\ell}, \ \ell \geq n\}$. By a slight abuse of notation, we now denote $G_n^y = \pi_1(X_n, y_n)$, and $G_\ell^y = (f^\ell_n)_*\pi_1(X_{\ell}, y_{\ell})$ (these groups are isomorphic to the groups $(f^\ell_0)_*\pi_1(X_{\ell}, y_{\ell})$, which we denoted by $G_{\ell}^y$ earlier). Thus we have a homeomorphism
   $$\phi_n'  \colon  U^y_n \to X_{\infty,n}^y = \lim_{\longleftarrow}\{G^y_n/ G^y_{\ell+1} \to G^y_n/G^y_{\ell}\} \ ,$$
 which commutes with the action of $G_n^y$ on $U_y^n$ and $X_{\infty,n}^y$. Denote by $U(y,\ell)$ the cylinder neighborhoods of $(eG_{\ell}^y)$ in $X_{\infty,n}^y$. In particular, $U(y,n) = X_{\infty,n}$.
 
The surface $X_{\ell}$ in the presentation $\cP'$ has genus $m_{\ell} = 3^\ell +1$ (see \cite{CFL2014}), so $G_{\ell}^y $ has $m_{\ell}$ generators, represented by longitudinal loops. In particular, there are  loops $\gamma_n$ and $\delta_n$, which wind around the handles $H_n$ and $F_n$ in $X_n$ respectively. Denote by $g_\gamma$ and $g_\delta$ the elements represented by $\gamma_n$ and $\delta_n$ in $G_n^y$ respectively. 

Now consider the construction of the surface $X_{n+1}$. It is obtained by the identification of three copies $\overline{X}_n^{1,2,3}$ of $X_n$ similar to the identification in Figure~\ref{fig:M0-pic}, c).  There is a point $y_{n+1}$ in one of the copies, which satisfies $f^{n+1}_n(y_{n+1})= y_{n}$, and which corresponds to our choice of the basepoint $y$. Denote by $z_{n+1}$ and $v_{n+1}$ the other two points such that  $f^{n+1}_n(z_{n+1}) = f^{n+1}_n(v_{n+1})= y_n$.

Denote by $\gamma_{y_{n+1}}$, $\gamma_{z_{n+1}}$, and $\gamma_{v_{n+1}}$ the copies of $\gamma_n$ in  $\overline{X}_n^{1,2,3}$ with respective basepoints $y_{n+1}$, $z_{n+1}$ and $v_{n+1}$. Note that these loops are cut when constructing $\overline{X}_n^{1,2,3}$. We now proceed to identify the boundaries of $\overline{X}_n^{1,2,3}$ according to the construction, which would close one of the loops back, and would intertwine the boundaries of the other two loops, so as to create a single loop of twice the length of $\gamma_n$. 

We have the following alternatives: first, suppose $\gamma_{z_{n+1}}$ is identified into a loop, and $\gamma_{y_{n+1}}$ and $\gamma_{v_{n+1}}$ are identified to make a single loop of twice the length. Then the lift of $\gamma_{n}$ with the starting point $y_{n+1}$ is the curve $\gamma_{y_{n+1}}$ which is not closed and has $v_{n+1}$ as its ending point. This means that the action of $g_{\gamma}$ on the coset space $G_n^y/G_{n+1}^y$ maps $eG_{n+1}^y$ onto $g_\gamma G_n^y$, and so maps the cylinder neighborhood $U(y,n+1)$ onto the clopen set $g_\gamma (U(y,n+1))$. At the same time, the lift of
 $\gamma_n$ with the starting point $z_{n+1}$ is a closed loop. So the action of $g_\gamma$ fixes the coset $\gamma_\delta G_{n+1}^y$, and the clopen set $g_\delta (U(y,n+1))$. We note that on the subsequent steps of the construction, when creating $X_{n+i}$, the lifts of the loop $\gamma_{z_{n+1}}$ are never cut and identified, which means that the action of $g_\gamma$ is the identity on $g_\delta (U(y,n+1))$.
 
 Another alternative is that $\gamma_{y_{n+1}}$ is identified into a loop, and $\gamma_{z_{n+1}}$ and $\gamma_{v_{n+1}}$ are identified to make a single loop. Arguing similarly, in this case we obtain that the action of $g_{\gamma}$ is the identity on $U(y,n+1)$, and it permutes the sets $g_\delta(U(y,n+1))$ and $g_\gamma \circ g_\delta(U(y,n+1))$. Thus in both cases we obtain an element which is the identity on a clopen subset of the section $U^y_n$, which permutes two other subsets of $U_y^n$, which means that $(U^y_n,G_n^y,\Phi)$ is not \emph{SQA} . Since the choice of $y$ and $n$ was arbitrary, we conclude that the holonomy pseudogroup for the Schori solenoid is not \emph{SQA} . \endproof

}
\end{ex} 

From the proof of Theorem~\ref{Schori-notqa} we obtain the following corollary, which shows that the hypotheses of  Proposition~\ref{prop-kernel-inclusion} are necessary.
 
  \begin{cor}\label{cor-nouniformsection}
  In the Schori solenoid, for any choice of a transverse section $V_0$, and any choice of a point $x$, $Germ(\Phi_0,x)$ is not locally trivial.
  \end{cor}

\proof From the proof of Theorem~\ref{Schori-notqa} we conclude that, for any choice of a basepoint $y \in \fX_0$, and any choice of a group chain $\cG^y_n=\{G_n^y\}_{i \geq 0}$, the kernel $ K(G^y_n)$  is not a normal subgroup of $G_n^x$. It follows that, even if $Germ(\Phi_0,x)$ is trivial, it is not locally trivial.
\endproof

\section{A universal construction}\label{sec-construction}

In this section, we give a general method to construct  examples of group chains with prescribed discriminant groups. 
This construction is inspired by  the proof of Lemma~37 in Section~8 of  Fokkink and Oversteegen \cite{FO2002}, which they attribute to Hendrik Lenstra. The construction of Lenstra is given in Section~\ref{subsec-lenstra}, and Section~\ref{subsec-lenstraremks} discusses some properties of this construction. Then in Section~\ref{subsec-stable} we give   criteria for when the resulting group chains are stable.

Section~\ref{subsec-lubotzky} recalls two basic results   of Lubotzky from the paper \cite{Lubotzky1993}. The first, given here as Theorem~\ref{thm-lubotzky1},    realizes    any given  finite group $F$  embedded into the profinite completion  of a finitely-generated, torsion-free group $G$.  
A second result of Lubotzky, given here as Theorem~\ref{thm-lubotzky2},  embeds the infinite product $\bH$ of a collection of finite groups  as a subgroup of the profinite completion of a finitely-generated, torsion-free group $G$.    Then in Section~\ref{subsec-stableexamples}, these constructions of Lubotzky are used to construct the examples used in the proofs of Theorems~\ref{thm-stablefinitefibre} and \ref{thm-stableCantorfibre}   of the introduction.

There is an extensive literature on embedding   groups into the profinite completion of a given torsion-free, finitely-generated group (see \cite{RZ2010} for a discussion of this topic, and further references).  The methods of this section apply in this generality to yield an enormous range of equicontinuous minimal Cantor actions with infinite, hence Cantor  discriminant groups.  

\subsection{A profinite construction}\label{subsec-lenstra}
We first give a reformulation of the constructions in Sections~\ref{subsec-ellischains} and \ref{subsec-disc}, in analogy with  the construction of Lenstra in \cite{FO2002}.  This alternate formulation is of strong interest in itself, as it gives a deeper understanding  of the Molino spaces introduced in this work.

Let $G_0$ be a finitely-generated group, $\cG = \{G_{\ell}\}_{\ell \geq 0}$ a group chain in $G_0$, and let $\cC = \{C_{\ell}\}_{\ell \geq 0}$ be the core group chain associated to $\cG$, with $C_{\infty}$ the core group associated with $\cC$. Also assume that the kernel $K(\cG)$ is trivial, so the map $\whPhi \colon    G_0 \to C_{\infty}$ is an injective homomorphism with dense image  $\whG_0 = \whPhi(G_0) \subset C_{\infty}$.  Then the discriminant   of  $\cG$ is a compact subgroup $\cD \subset C_{\infty}$, whose rational core as defined in \eqref{eq-rationalcore} is trivial by Proposition~\ref{prop-discrisnotnormal}.

Let $C_{n,\infty} \subset C_{\infty}$ be the clopen normal subgroup neighborhood of the identity $\{e\}$ defined in \eqref{Cinfty-definen}. As  the intersection ${\ds \bigcap_{n \geq 1} \ C_{n,\infty}  = \{e\}}$,  the collection $\{C_{n,\infty} \mid n \geq 1\}$ is a clopen neighborhood system about the identity in $C_{\infty}$.
Observe that from the definition \eqref{cinfty-coords}, we have that 
$\ds C_{\infty}/C_{n,\infty} \cong G_0/C_n$ and   $\whG_0 \cap C_{n,\infty}  \cong G_n$.
As each subgroup $C_{n,\infty} $ is normal and $\cD$ is compact, the product $V_n = \cD \cdot C_{n,\infty}  \subset C_{\infty}$ is a clopen subgroup of $C_{\infty}$ containing $\cD$, and we have 
${\ds \cD = \bigcap_{n \geq 1} \ V_n}$. Thus, $\cD$ is realized as the countable intersection of clopen subgroups of $C_{\infty}$. It is an exercise to show that this formulation of $\cD$ agrees the definition of $\cD$ as an inverse limit in \eqref{eq-discriminantdef}.

We now turn the order of the above remarks around, to obtain a construction of a group chain with prescribed discriminant group.

\begin{prop}\label{prop-lenstra}
Let $C_{\infty}$ be a profinite group, and let $G \subset C_{\infty}$ be a finitely-generated dense subgroup.
Let $\cD \subset C_{\infty}$ be a compact subgroup of infinite index which has trivial rational core,  
\begin{equation}\label{eq-rationalcorepro}
{\rm core}_{G} \  \cD =    \bigcap_{k \in G} k   \cD   k^{-1} = \{e\} \ .
\end{equation}
Then there exists a group chain $\cG = \{G_{\ell}\}_{\ell \geq 0}$ with $G_0 = G$, with discriminant group     $\cD$.
\end{prop}

\proof
By the assumption that $C_{\infty}$ is a profinite group, there exists a group chain  $\{U_{\ell}  \mid \ell \geq 1\}$ which is a clopen neighborhood system about the identity in $C_{\infty}$, so that: 
\begin{enumerate}
\item each $U_{\ell}$ is normal in $C_{\infty}$ ;
\item for each $\ell \geq 0$ there is a proper inclusion $U_{\ell+1}  \subset U_{\ell}$ ;
\item the  intersection ${\ds \bigcap_{\ell \geq 1} \ U_{\ell}   = \{e\}}$ .
\end{enumerate}
 In particular, each quotient $H^{\ell} \equiv C_{\infty}/U_{\ell}$ is a finite group. Let $\iota^{\ell+1}_{\ell} \colon H^{\ell+1} \to H^{\ell}$ be the map induced by inclusion of cosets.
 Then   there is a natural identification 
 \begin{equation}
C_{\infty} \cong \lim_{\longleftarrow}\, \left \{\iota^{\ell+1}_{\ell} \colon  H^{\ell +1}  \to H^{\ell}  \right\} \ .
\end{equation}
Next, for each $\ell \geq 1$, set $W_{\ell} = \cD \cdot U_{\ell}$ which is a   subgroup of $C_{\infty}$, as $U_{\ell}$ is normal. Moreover, the assumption that $\cD$ is compact implies that each $W_{\ell}$ is a clopen subset of $C_{\infty}$. Then set  $G_{\ell} = G \cap W_{\ell}$ which is a subgroup of finite index in $G$, and so $\cG = \{G_{\ell}\}_{\ell  \geq 0}$ is a subgroup chain in $G$. Note that 
\begin{equation}\label{eq-kernelpro}
K(\cG) = \bigcap_{\ell \geq 0} \  G_{\ell} =  \bigcap_{\ell \geq 0} \  G \cap W_{\ell} = G \cap  \bigcap_{\ell \geq 0} \   W_{\ell} = G \cap \cD \ . 
\end{equation}

We next calculate the discriminant   of the chain $\cG$. Let $\pi^{\ell} \colon C_{\infty} \to H^{\ell}$ be the quotient map. 
As  each $H^{\ell}$ is finite, the image $\cD^{\ell} \equiv \pi^{\ell}(\cD)$ is a finite set. 
The group $G$ is dense in $C_{\infty}$ so has non-trivial intersection with  each clopen set $g  U_{\ell}$. Thus,  
\begin{equation}\label{eq-discell}
\cD^{\ell} = \pi^{\ell}(\cD) = \pi^{\ell}(W_{\ell}) = \pi^{\ell}(G \cap W_{\ell}) = \pi^{\ell}(G_{\ell}) \subset H^{\ell} \ .
\end{equation}
The core of the group $G_{\ell}$ is the   group
$\ds  C_{\ell}   \equiv {\rm core}_{G} \  G_{\ell}   =    \bigcap_{g \in G} \ g G_{\ell} g^{-1}$. Then   we have
\begin{equation}\label{eq-coreell}
\pi^{\ell}(C_{\ell}) = \pi^{\ell}(\bigcap_{g \in G} \ g \ G_{\ell} \ g^{-1}) = \bigcap_{g \in G} \ \pi^{\ell}(g) \ \pi^{\ell}(G_{\ell}) \ \pi^{\ell}(g)^{-1} = \bigcap_{g \in G} \ \pi^{\ell}(g) \ \pi^{\ell}(\cD) \ \pi^{\ell}(g)^{-1}  = \{ e^{\ell} \} \ ,
\end{equation}
where $e^{\ell} \in H^{\ell}$ is the identity, and the last equality follows since $G$ is dense in $C_{\infty}$ and the core of $\cD$ is trivial.
It follows that $C_\ell = G \cap U_\ell$, and thus we obtain induced maps on the quotients, $\ovpi^{\ell} \colon G/C_{\ell} \to H_{\ell}$.
Then note that  $\pi^{\ell}(G_{\ell}/C_{\ell}) = \pi^{\ell}(\cD) = \cD^{\ell}$ for all $\ell \geq 0$.

The map $\iota^{\ell+1}_{\ell} \colon H^{\ell+1} \to H^{\ell}$  induces a map (denoted   the same),  
$\iota^{\ell+1}_{\ell} \colon \cD^{\ell+1} \to \cD^{\ell}$.
Then for the inverse limits we have
\begin{equation}\label{eq-discpro}
\lim_{\longleftarrow}\, \left \{\delta^{\ell+1}_{\ell} \colon  G_{\ell+1}/C_{\ell+1} \to G_{\ell}/C_{\ell}  \right\} = \lim_{\longleftarrow}\, \left \{\iota^{\ell+1}_{\ell} \colon  \cD^{\ell+1} \to \cD^{\ell}   \right\}   \   .
\end{equation}
The term on the left-hand-side of \eqref{eq-discpro} is by definition the discriminant of the chain $\cG$, while the term on the right hand side of \eqref{eq-discpro} is homeomorphic  to the subgroup $\cD$, as $\{U_{\ell} \mid \ell \geq 1\}$   is a clopen neighborhood system about the identity in $C_{\infty}$. 
\endproof

\subsection{Properties of the Lenstra construction} \label{subsec-lenstraremks}
 We make some remarks about the   construction in Proposition~\ref{prop-lenstra}. First,  note that  
 the    proof of   \cite[Lemma~37]{FO2002}       defined the chain $\cG_n$ using a collection of clopen neighborhoods of $e \in C_{\infty}$. However,   the proof   in \cite{FO2002}  that the chain $\cG_n$ is not weakly regular used Proposition~\ref{prop-conjclasses}, that is, the fact that if the number of conjugacy classes of the kernel $K(\cG_n)$ is infinite, then $\cG_n$ cannot be weakly regular. Our approach is to calculate the discriminant group   for the chain directly.

 Assume there is given a profinite group $C_{\infty}$, a compact subgroup $\cD \subset C_{\infty}$, and a dense subgroup $G \subset C_{\infty}$ satisfying the hypotheses of Proposition~\ref{prop-lenstra}. Set $X = C_{\infty}/\cD$ which is a Cantor space. The left action of $G$ on $X$ defines a map $\Phi \colon G \to Homeo(X)$, which is a minimal action as $G$ is dense in $C_{\infty}$. Thus, the construction yields   an equicontinuous minimal Cantor system $(X, G , \Phi)$.
 
 Next, given a   clopen neighborhood system $\{U_{\ell} \mid \ell \geq 1\}$   about the identity in $C_{\infty}$, which satisfies the conditions in the proof of Proposition~\ref{prop-lenstra},   let $\cG \equiv \{G_{\ell}\}_{\ell \geq 0}$ be the group chain in $G$ constructed with respect to this clopen neighborhood system.  Then it is an exercise, using the techniques of the proof of Proposition~\ref{prop-lenstra}, to show that there is a $G$-equivariant homeomorphism of spaces  
 $$\tau \colon X  \cong \lim_{\longleftarrow}\, \left \{\iota_{\ell +1} \colon  G/C_{\ell +1} \to G/G_{\ell}  \right\} \equiv X_{\infty} \ .$$
Now suppose that  $\{V_{\ell} \mid \ell \geq 1\}$  is another  clopen neighborhood system about the identity in $C_{\infty}$, which also satisfies the conditions in the proof of Proposition~\ref{prop-lenstra}, and let $\cH \equiv \{H_{\ell}\}_{\ell \geq 0}$ be the group chain in $G$ constructed with respect to this second clopen neighborhood system. 
A basic property of neighborhood systems is that given any $\ell \geq 0$ there exists $\ell' \geq 0$ such that $V_{\ell'} \subset U_{\ell}$, and $\ell'' \geq 0$ such that $U_{\ell''} \subset V_{\ell}$.  It follows from their definitions that the group chains $\cG$ and $\cH$ are equivalent in the sense of Definition~\ref{defn-greq}.

 Suppose that $G \cap \cD = \{e\}$, then the calculation  \eqref{eq-kernelpro} shows that the kernel   $K(\cG) = \{e\}$ is trivial. 
 Moreover, suppose the choice of $\cD$ is made so that 
 $G \cap \whg \cD \whg^{-1} = \{e\}$ for all $\whg \in C_{\infty}$. Given $y \in X$ let    $\tau(y) = (g_{\ell} G_{\ell}) \in X_{\infty}$, and let  $\cG^y = \{g_{\ell} G_{\ell} g_{\ell}^{-1}\}_{\ell \geq 0}$ be the conjugate group chain. Choose $\whg \in C_{\infty}$ such that $\tau(\whg \cD) = (g_{\ell} G_{\ell})$. Then 
  \begin{equation}
K(\cG^y) =  G  \ \cap \ (\whg \ \cD \ \whg^{-1}) =\{e\}
\end{equation}
 so that $\cG^y$ also has trivial kernel. Thus, if we chose $\cD$ so that  $G \cap \whg \cD \whg^{-1} = \{e\}$ for all $\whg \in C_{\infty}$ is satisfied, then 
 the Cantor system $(X, G , \Phi)$ has trivial kernel for the group chain $\cG^y$ at $y$, for all points $y \in X$.
  For example, suppose that $G$ is a torsion-free group, and $\cD$ is a torsion group. Then  the   condition  $G \cap \whg \cD \whg^{-1} = \{e\}$ for all $\whg \in C_{\infty}$ is automatically satisfied, as each non-trivial element of $\cD$, and hence $\whg \cD \whg^{-1}$, has finite order. We use this observation in Theorem~\ref{thm-finiterealization}  below.
 
 On the other hand, given $G \subset C_{\infty}$ as in Proposition~\ref{prop-lenstra}, and suppose that the compact subgroup $\cD \subset C_{\infty}$ is chosen so that  $G \cap \whg \cD \whg^{-1} \ne \{e\}$ for some $\whg \in C_{\infty}$, then by  Proposition~\ref{prop-kernel-discriminant} there exists $y \in X$ such that the Cantor system $(X, G , \Phi)$ has non-trivial kernel $K(\cG^y)$ for the group chain $\cG^y$ about $y$. It then follows that the germinal holonomy group ${\rm Germ}(\Phi , y)$ is non-trivial, so this method can also be used to construct examples with non-trivial germinal holonomy groups.

\subsection{Stable actions}\label{subsec-stable}

Recall from Definition~\ref{def-stableGC} that a  group chain $\cG  = \{G_{\ell}\}_{\ell \geq 0}$ is said to be \emph{stable} if there exists $n_0 \geq 0$ such that the maps $\ds \psi_{n,m} \colon \cD^n  \to    \cD^m$ defined in \eqref{eq-discmapsnm2} are isomorphisms for all $m \geq n \geq n_0$.  We consider the problem of when a group chain $\cG =  \{G_{\ell}\}_{\ell \geq 0}$ constructed using the method of  proof of Proposition~\ref{prop-lenstra} is stable. 

 We assume the hypotheses of  Proposition~\ref{prop-lenstra}, and the   constructions of its proof. Fix $n > 0$, and consider the truncated group chain $\cG_n = \{G_{\ell}\}_{\ell \geq n}$.  Then the calculation of the kernel $K(\cG_n) = G \cap \cD$ is   the same as   \eqref{eq-kernelpro}. 
Also, note that $\cD \subset W_{\ell}$ for all $\ell \geq 0$,   so the calculations in  \eqref{eq-discell} also proceed analogously. However, the last equality in \eqref{eq-coreell} 
requires   the additional assumption
\begin{equation}\label{eq-rationalcorepro2}
{\rm core}_{U} \  \cD =    \bigcap_{k \in U} k   \cD   k^{-1} = \{e\}  
\end{equation}
for the clopen neighborhoods $U = U_{\ell}$ of the identity, in order to conclude that $\cD$ is the discriminant group for $\cG_n$.
In other words, we require that the subgroup $\cD$ is ``totally not-normal'' for every neighborhood of the identity in $\whG$.
 The above remarks   yield: 
\begin{prop}\label{prop-lenstran}
Let $C_{\infty}$ be a profinite group,   let $G \subset C_{\infty}$ be a finitely-generated dense subgroup, and 
let $\cD \subset C_{\infty}$ be a compact subgroup of infinite index, such that \eqref{eq-rationalcorepro2} holds for every clopen neighborhood $\{e\} \in U \subset C_{\infty}$.
Choose a  group chain  $\{U_{\ell} \mid \ell \geq 1\}$ which is a clopen neighborhood system about the identity in $C_{\infty}$, then the associated group chain $\cG = \{G_{\ell}\}_{\ell \geq 0}$ with $G_0 = G$ has discriminant group  $\cD$ and is stable.
\end{prop}

Finally, in the case where $\cD \subset C_{\infty}$ is a compact subgroup of infinite index, but need not satisfy the condition that its   core is trivial, then noting that the core is a normal subgroup, we can modify the construction above as follows to obtain a minimal Cantor action.

\begin{cor}\label{cor-anydisc}
Let $C_{\infty}'$ be a profinite group,     $G' \subset C_{\infty}'$ be a finitely-generated dense subgroup, and 
  $\cD' \subset C_{\infty}'$ be a non-trivial compact subgroup of infinite index , and let ${\rm core}_{G'}    \cD' $ denote the rational core of $\cD'$ as in \eqref{eq-rationalcorepro2},  which is a normal subgroup of $C_{\infty}'$ as $G'$ is dense. Set:
  $$C_{\infty} = C_{\infty}'/({\rm core}_{G'}    \cD') ~ ; ~ G = G'/(G' \cap {\rm core}_{G'}    \cD') ~ ; ~ \cD = \cD'/{\rm core}_{G'}    \cD' \ . $$
Then there exists a group chain $\cG = \{G_{\ell}\}_{\ell \geq 0}$ with $G_0 = G$, and discriminant group  $\cD$.
\end{cor}

 \subsection{Constructing embedded groups}\label{subsec-lubotzky}
We next recall the remarkable constructions  of  Lubotzky, which when combined with the techniques of Proposition~\ref{prop-lenstra}, makes possible the construction of a wide class of equicontinuous minimal Cantor actions by a    finitely-generated, torsion-free, residually finite group $G$, with prescribed discriminant group $\cD$. 
There are two cases of the construction. 
  
  \begin{thm}[Theorem~2(b), Lubotzky \cite{Lubotzky1993}] \label{thm-lubotzky1}
  Let $F$ be a non-trivial finite group, and set $\bF_i = F$ for all integers $i \geq 1$. Let $\ds \bF = \prod \ \bF_i$ denote the infinite cartesian product of $F$.
  Then there exists a finitely-generated, residually-finite, torsion-free group $G \subset {\bf SL}_n(\mZ)$ for $n \geq 3$ sufficiently large,  whose profinite completion $\whG$ contains  $\bF$.
  \end{thm}
   \proof 
  We give just an outline   of the construction used in the proof of Theorem~2(b) in  \cite{Lubotzky1993}, with   details as required for the constructions of our examples. 
  First recall some basic facts.   For $n \geq 3$, let $\G_n = {\bf SL}_n(\mZ)$ denote the $n \times n$ integer matrices.
       The group $\G_n$ is finitely-generated and residually finite, and hence so are all finite index subgroups of $\G_n$.  
Let $\G_n(m)$ denote the congruence subgroup 
    $$\G_n(m) \equiv  {\rm Ker} \left\{\varphi_m \colon  {\bf SL}_n(\mZ) \to {\bf SL}_n(\mZ/m \mZ) \right\} \ .$$
For $m \geq 3$, $\G_n(m)$ is torsion-free.  
 Moreover,  by the congruence subgroup property, every finite index subgroup of $\G_n$ contains $\G_n(m)$ for some non-zero $m$.    
 Then this implies 
 \begin{equation}\label{eq-prodprimes}
\widehat{{\bf SL}_n(\mZ)} \cong \lim_{\longleftarrow} ~ {\bf SL}_n(\mZ/m \mZ) \cong {\bf SL}_n(\widehat{\mZ}) \cong  \prod_p \ {\bf SL}_n(\mZ_p)
\end{equation}
where $\widehat{\mZ}$ is the profinite completion of $\mZ$, and we use that $\ds \widehat{\mZ} = \prod_p \ \mZ_p$ where $\mZ_p$ is the ring of $p$-adic integers, and the product is taken over all primes.  Note that the factors in the cartesian product on the right hand side of  \eqref{eq-prodprimes} commute with each other.
   
  Let $G \subset \G_n$ be a   finite index, torsion-free subgroup, which is then finitely-generated, 
  and its profinite completion $\whG$ is an open subgroup of $\widehat{{\bf SL}_n(\mZ)}$. Then there exists a cofinite subgroup $\cP(G)$ of the primes such that 
 \begin{equation}\label{eq-embedp}
 \prod_{p \in \cP(G)} \ {\bf SL}_n(\mZ_{p}) \subset \whG \ .
\end{equation}

 Let $d_F= |F|$ denote the cardinality of $F$, and let $n \geq |F|+2$. Then   $F$ embeds in the alternating group ${\bf Alt}(n)$ on $n$ symbols.
  Then $F$ is non-trivial implies that $n \geq 4 > 3$.
 For each $p \in \cP(F)$, the     group ${\bf Alt}(n)$ embeds into ${\bf SL}_{n}(\mZ_p)$, and thus  we obtain an embedding 
    \begin{equation}\label{eq-embedpH}
\iota_{\infty} \colon  \bF \cong  \prod_{p \in \cP(G)} \ F_p  \subset  \prod_{p \in \cP(G)} \ {\bf Alt}(n)    \subset  \prod_{p \in \cP(G)} \ {\bf SL}_{n}(\mZ_{p}) \subset \whG \ ,
\end{equation}
where $F_p = F$ for each $p \in \cP(G)$.  This completes the construction. 
 \endproof

In  \cite[Theorem~1]{Lubotzky1993} Lubotzky extended the above construction, to obtain an embedding  for a   group $\cD$ which is an infinite product of possibly distinct finite groups $\{\bH_i \mid i = 1,2, \ldots\}$. The extension is highly non-trivial, as if all of the groups $\bH_i$ are distinct, then  the degrees $|\bH_i|$ must tend to infinity, and so the above straightforward strategy for embedding no longer works.

    \begin{thm}[Theorem~1, Lubotzky \cite{Lubotzky1993}] \label{thm-lubotzky2}
  Let $\{\bH_i \mid i =1,2, \ldots\}$ be an infinite     collection of non-trivial finite groups, and let $\ds \bH = \prod \ \bH_i$ denote their cartesian product.
  Then there exists a finitely-generated, residually-finite, torsion-free group $G$ whose profinite completion $\whG$ contains  $\bH$.
  \end{thm}
   \proof
 Again, we only sketch some key aspects of the proof from \cite{Lubotzky1993}.
  Let $G \subset \G_n$ be the   finitely-generated, torsion-free, residually-finite   group constructed  on  \cite[page 330]{Lubotzky1993}, and $\whG$ its profinite completion. 
   Lubotzky constructs by induction an increasing  sequence of primes $\{p_n \mid n \geq 3\}$ such that 
   \begin{equation}\label{eq-embedpp}
 \prod_{n=3}^{\infty} \ {\bf SL}_n(\mZ_{p_n}) \subset \whG \ .
\end{equation}
   For $i \geq 1$, let $d_i = |\bH_i|$ denote the cardinality of $\bH_i$. Then each $\bH_i$ embeds in the alternating group ${\bf Alt}(d_i+2)$ on $d_i+2$ symbols.
   Now choose an increasing sequence of integers $\{n_i \mid i \geq 1\}$ such that $n_i \geq d_i +2$. Then for each $i \geq 1$, the group ${\bf Alt}(d_i+2)$ embeds into the alternating group ${\bf Alt}(n_i)$  by taking only the   permutations on the first $(d_i+2)$ symbols. 
 For each $i \geq 1$ the     group ${\bf Alt}(n_i)$ embeds into ${\bf SL}_{n_i}(\mZ_{p_{n_i}})$. Thus,  we have  embeddings $\bH_i \subset {\bf Alt}(n_i) \subset {\bf SL}_{n_i}(\mZ_{p_{n_i}})$. 
 
The product in \eqref{eq-embedpp} is over all $n \geq 3$, while the group $\bH_n = \bH_{n_i}$ if $n = n_i$ for some $n_i$ as chosen above. For $n \ne n_i$ for some $i$, let $\bH_n$ be the trivial group.  Set  $\bA_n = {\bf Alt}(n_i)$ if $n = n_i$ for some $n_i$ and let $\bA_n$ be the trivial group otherwise.
Then we obtain an embedding of the infinite product $\bH$, 
   \begin{equation}\label{eq-embedallH}
\iota_{\infty} \colon  \cD \cong  \prod_{n \geq 3} \ \bH_n  \subset  \prod_{n \geq 3} \ \bA_n  \subset   \prod_{i \geq 1} \ {\bf SL}_{n_i}(\mZ_{p_{n_i}})\subset \whG \ .
\end{equation}
This completes the construction.
\endproof

 \subsection{Constructing   stable actions}\label{subsec-stableexamples}
 
We next use  Theorems~\ref{thm-lubotzky1} and \ref{thm-lubotzky2}, and observations from their proofs in  \cite{Lubotzky1993},  to construct 
 examples of stable equicontinuous minimal Cantor group actions.

We first require a simple observation. For $n \geq 2$, the alternating group ${\bf Alt}(n)$ on $n$ symbols embeds into the the alternating group ${\bf Alt}(4n)$ on $4n$ symbols, by consider ${\bf Alt}(n)$ as acting on the first $n$ symbols, and fixing the remaining $3n$ symbols.  We thus consider $ {\bf Alt}(n) \subset {\bf Alt}(4n)$ as a subgroup.
\begin{lemma}\label{lem-altcore}
The core of ${\bf Alt}(n)$ in ${\bf Alt}(4n)$ is the trivial group.
\end{lemma}
\proof
There exists an element    $\sigma \in {\bf Alt}(4n)$ which swaps the first $2n$ symbols for the last $2n$ symbols. Then $\sigma^{-1} \   {\bf Alt}(n) \ \sigma$ is contained in the alternating group which permutes the last $2n$ symbols, and hence is disjoint from the subgroup  ${\bf Alt}(n)$.
\endproof

Lemma~\ref{lem-altcore} is used   to ensure that the chains constructed below satisfy the conditions of Section~\ref{subsec-stable}.

\begin{thm} \label{thm-finiterealization} 
Let $F$ be  a finite group. Then 
there exists   a  finite index,  torsion-free group $G \subset {\bf SL}_n(\mZ)$ and an embedding  of $F$ into the profinite completion $\whG$,  
so that the resulting group chain $\cG_K = \{G_{\ell}\}_{\ell \geq 0}$ constructed as in Section~\ref{subsec-lenstra}      yields    
  an equicontinuous minimal Cantor system $(X_{\infty}, G , \Phi)$   whose     
     discriminant group for the     truncated group chain $\{G_{\ell}\}_{\ell \geq k}$ is   isomorphic to $F$ for all $k \geq 0$. Hence the action is stable and irregular. 
Moreover,   the germinal holonomy group for each $x \in X$ is trivial. 
\end{thm}

\proof
As noted in the proof of    Theorem~\ref{thm-lubotzky1}, if $F$ is a non-trivial  finite group of order $d_F = |F|$, then $F$ embeds in the alternating group ${\bf Alt}(d_F +2)$. 
We then    embed ${\bf Alt}(d_F +2)$ in the alternating group ${\bf Alt}(n)$ for $n \geq 4(d_F +2)$, by considering ${\bf Alt}(n)$ as acting on the first $n$ symbols,   as in the proof of Lemma~\ref{lem-altcore}. We identify $F$ with its image, and then note that the core of $F$ in ${\bf Alt}(n)$ is the trivial group. Note that $d_F \geq 2$, so we have that $n \geq 16$. Also note that if $F'$ is any other finite group of order at most $d_F$, then it also embeds into ${\bf Alt}(d_F +2)$, and hence the following construction is universal for all   finite groups $F'$ with $|F'|\leq |F|$.

For $n \geq 4(d_F +2)$, let $G \subset \G_n = {\bf SL}_n(\mZ)$ be the 
  finite index, torsion-free subgroup  constructed in the proof of Theorem~\ref{thm-lubotzky1}. 
 Set $\bH_{\ell} = {\bf Alt}(n)$ for all integers $\ell \geq 1$, and let $\ds \bH = \prod \ \bH_{\ell}$ denote their cartesian product. Then the embedding \eqref{eq-embedpH} becomes
  \begin{equation} \label{eq-embedH}
\iota_{\infty} \colon  \bH \cong     \prod_{p \in \cP(G)} \ {\bf Alt}_p(n)    \subset  \prod_{p \in \cP(G)} \ {\bf SL}_{n}(\mZ_{p}) \subset \whG \ ,
\end{equation}
where ${\bf Alt}_p(n) =  {\bf Alt}(n)$ for each prime $p$. 

For  each $i \geq 1$, we have the embedding $F \subset {\bf Alt}(d_F +2) \subset {\bf Alt}(n) = \bH_{\ell}$. Let $F \to \bH$ be the diagonal embedding into the infinite product, which then yields an embedding $\iota_F \colon F \to \whG$ into the profinite completion of $G$, with image denoted by $\cD = \iota_F(F)$.

Next, use the method of Section~\ref{subsec-lenstra} to construct a group chain in $G$. The group $G$ is residually finite, so there exists a clopen neighborhood system 
 $\{U_{\ell}  \mid \ell \geq 1\}$   about the identity in $\whG$, where each $U_{\ell}$ is normal in $\whG$.  Note that $G$ is dense in $\whG$ and each $U_{\ell}$ is closed, so the closure of $G \cap U_{\ell}$ in $\whG$ is equal to $U_{\ell}$. 
Set $W_{\ell} = \cD \cdot U_{\ell}$ for $\ell \geq 1$, and $G_{\ell} = G \cap W_{\ell}$. Let $\cG_{F} = \{G_{\ell}\}_{\ell \geq 0}$ denote the resulting group chain.

Let $\{e\} \in U \subset \whG$ be a normal clopen neighborhood of the identity, so that $\whG/U$ is a finite group with cardinality $|\whG/U|$.
We claim    $\ds {\rm core}_{U} \  \cD = \{e\}$. 
The normal subgroup $U$ has finite index, hence as argued in the proof of   \cite[Theorem~2]{Lubotzky1993}, there exists  a cofinite subset of primes $\cP(G,U) \subset \cP(G)$ of the list in the product in \eqref{eq-embedH} such that 
$$ \prod_{p \in \cP(G,U)} \ {\bf Alt}_p(n)  \subset  \prod_{p \in \cP(G,U)} \ {\bf SL}_{n}(\mZ_{p})  \subset  U \subset    \whG \ . $$
For $p \in \cP(G,U)$, note that for the diagonal embedding of $F$ into $\bH$, the projection to each factor of $\bH$ is an isomorphism. For the image of $F$ in the $p-th$ factor, 
we have 
$$F \subset {\bf Alt}(d_F +2) \subset {\bf Alt}(n) = {\bf Alt}_p(n) \subset {\bf SL}_{n}(\mZ_{p}) \ .$$
The image group has trivial core by Lemma~\ref{lem-altcore}.   The projection of $\cD$ to $F \subset {\bf Alt}_p(n)$ is an isomorphism, so this implies that   $\cD$ has trivial core in $U$ as well.
Then by Proposition~\ref{prop-lenstran}, for all $k \geq 0$, the discriminant group for the truncated group chain  $\{G_{\ell}\}_{\ell \geq k}$    is   isomorphic to $F$. In particular, $\cG_F$ is a stable group chain.

Next, observe that $\cD$   compact implies that the closure $\overline{G_{\ell}}$ of $G_{\ell}$ in $\whG$ equals $W_{\ell}$, and $\ds \cD = \cap \ \overline{G_{\ell}}$. 
For the kernel of $\cG_F$ as defined in Section~\ref{subsec-kernel}, we then have
\begin{equation}
K(\cG_F) = \bigcap_{\ell \geq 0} \ G_{\ell} \subset \bigcap_{\ell \geq 0} \ \overline{G_{\ell}} = \cD \ .
\end{equation}
The group $\cD$ is finite, hence every element of $\cD$ has finite order, while $K(\cG_F) \subset G$ which is torsion-free. Thus, 
$K(\cG_F)  \subset \cD \cap G = \{e\}$, hence $K(\cG_F)$ is the trivial group. Moreover,  for each $\whg \in \whG$ let $\cG_F^{\whg} = \{\whg \  G_{\ell} \ \whg^{-1}\}_{\ell \geq 0}$ denote the conjugate group chain. Then by the same reasoning, we also have  $K(\cG_F^{\whg}) = \{e\}$, as $\whg^{-1} \ \cD \ \whg \subset \whG$ is again a finite subgroup, hence has trivial intersection with $G$.

 Let $(X, G , \Phi)$ be the  equicontinuous minimal Cantor system with $X = \whG/\cD$ with the associated group chain $\cG_F$, as discussed in Section~\ref{subsec-lenstraremks}. 
 The discriminant group of $\cG_F$ is $\cD$, and  each non-trivial element $h \in \cD$ is torsion, hence its image in $\whG$ is torsion, and thus any conjugate of it is not contained in the torsion-free subgroup $G$. Thus, for each $y \in X$, the action $\Phi$ has trivial germinal holonomy at $y$.
  
The  discriminant group of the truncated chain $\{G_{\ell}\}_{\ell \geq k}$ is   isomorphic to $F$ for all $k \geq 0$. Thus, $\cG_F$ cannot be a weakly regular group chain. This establishes all of the claims of Theorem~\ref{thm-finiterealization}.
 \endproof
 
  Note that the action $(X, G , \Phi)$ satisfies the SQA condition by default, as all germinal holonomy groups are trivial. The action of $\whG$ on $X = \whG/\cD$ satisfies the SQA condition by Theorem~\ref{thm-sqa-finitediscr}.
 Corollary~\ref{cor-stablefinitefibre} now follows by using the construction in Section~\ref{subsec-suspensions} to obtain a matchbox manifold with section $V_0 \cong X$ and induced holonomy action $(X, G , \Phi)$. 
 
 We remark that it is tempting to use the fact that $G \subset   {\bf SL}_n(\mZ) \subset {\bf SL}_n(\mR)$ is a torsion-free subgroup, and then use the quotient space $M_0 = {\bf SL}_n(\mR)/G$ as the base of a presentation for a weak solenoid $\cS_{\cP}$. However, this quotient space is not compact, and we do not have a ``theory of weak solenoids'' over non-compact manifolds.

 We next use Theorem~\ref{thm-lubotzky2} to construct two types of embeddings of Cantor groups into profinite groups. 
 Theorem~\ref{thm-stableCantor} embeds a profinite group such that the resulting action is stable. 
 Theorem~\ref{thm-nvr} embeds a Cantor group such that the resulting action is not virtually regular.  
 
  \begin{thm}\label{thm-stableCantor}
Let  $K$ be a  separable profinite  group. 
There exists a finitely-generated, residually-finite, torsion-free group $G$, and an  embedding  of $K$ into its profinite completion $\whG$, 
so that the resulting group chain $\cG_K = \{G_{\ell}\}_{\ell \geq 0}$ constructed as in Section~\ref{subsec-lenstra}  yields    
  an equicontinuous minimal Cantor system $(X, G , \Phi)$ whose     
     discriminant group for the     truncated group chain $\{G_{\ell}\}_{\ell \geq k}$ is   isomorphic to $K$ for all $k \geq 0$. Hence the action is stable and irregular. 
   \end{thm}
\proof

Let $G \subset \G_n$ be the   finitely-generated, torsion-free, residually-finite   group in the proof of Theorem~\ref{thm-lubotzky2}, and as constructed  on  \cite[page 330]{Lubotzky1993}, and let $\whG$ be its profinite completion.

The assumption that $K$ is a separable profinite group implies that $K$ is isomorphic to an inverse system of finite groups
\begin{equation}\label{eq-Kinvlim}
K \cong  \lim_{\longleftarrow} \, \left\{\phi^{\ell+1}_{\ell} \colon   \bK_{\ell+1}  \to \bK_{\ell}  \mid \ell \geq 0 \right\} \subset \bK \cong \prod \ \bK_{\ell} \ ,
\end{equation}
where each $\bK_{\ell}$ is a finite group, and the bonding maps $\phi^{\ell+1}_{\ell}$ are  epimorphisms for all $\ell \geq 0$, but not isomorphisms. 
Thus, their cardinalities $\{|\bK_{\ell}| \mid \ell \geq 0\}$ form an increasing sequence of integers. Note that we have isomorphisms for all $k > 0$, induced by the shift map $\sigma_i$ on indices,
\begin{equation}
\sigma_i \colon K \cong  \lim_{\longleftarrow} \, \left\{\phi^{\ell+1}_{\ell} \colon   \bK_{\ell+1}  \to \bK_{\ell}  \mid \ell \geq k \right\} \ .
\end{equation}

For each $\ell \geq 0$, set $d_{\ell} = 4(|\bK_{\ell}| +2)$.  Then as in the construction in Theorem~\ref{thm-finiterealization}, there is an embedding of $\bK_{\ell}$ into the alternating group,  $\ds \bK_{\ell} \subset {\bf Alt}(|\bK_{\ell}| +2) \subset {\bf Alt}(d_{\ell})$. 
 Choose an increasing sequence of integers $\{n_{\ell} \mid \ell \geq 1\}$ so that $n_{\ell} \geq d_{\ell}$ for all $\ell \geq 1$.

Then as in the proof of Theorem~\ref{thm-lubotzky2}, we set  $\bH_n = {\bf Alt}(d_{\ell})$ if $n = n_{\ell}$ for some $n_{\ell}$ as chosen above. If $n \ne n_{\ell}$ for all $\ell$, let $\bH_n$ be the trivial group.  Set  $\bA_n = {\bf Alt}(n_\ell)$ if $n = n_{\ell}$ for some $n_{\ell}$,  and let $\bA_n$ be the trivial group otherwise.
Then we obtain an embedding of the infinite product, 
   \begin{equation}\label{eq-embedallHalt}
   \bH \equiv   \prod_{n \geq 3} \ \bH_n  \subset  \bA \equiv \prod_{n \geq 3} \ \bA_n  \subset   \prod_{\ell \geq 1} \ {\bf SL}_{n_{\ell}}(\mZ_{p_{n_{\ell}}})\subset \whG \ .
\end{equation}
Now observe that the inverse limit presentation in \eqref{eq-Kinvlim}, along with the above embedding \eqref{eq-embedallHalt}, gives an embedding 
\begin{equation}
\Delta_K \colon  K   \subset \prod \ \bK_{\ell} \subset  \prod_{n \geq 3} \ \bH_n   \subset  \prod_{n \geq 3} \ \bA_n  \subset \whG \ .
\end{equation}

Set $\cD = \Delta_K(K) \subset \whG$.
 Then as in the proof of Theorem~\ref{thm-finiterealization},  use the method of Section~\ref{subsec-lenstra} to construct a group chain in $G$. 
 The group $G$ is residually finite, so there exists a clopen neighborhood system 
 $\{U_{\ell}  \mid \ell \geq 1\}$   about the identity in $\whG$, where each $U_{\ell}$ is normal in $\whG$.   Set $W_{\ell} = \cD \cdot U_{\ell}$ for $\ell \geq 1$, and $G_{\ell} = G \cap W_{\ell}$. Let $\cG_{K} = \{G_{\ell}\}_{\ell \geq 0}$ denote the resulting group chain.

Let $\{e\} \in U \subset \whG$ be a normal clopen neighborhood of the identity, so that $\whG/U$ is a finite group with cardinality $|\whG/U|$.
We claim    $\ds {\rm core}_{U} \  \cD = \{e\}$.    
Note that for  $m \geq 5$, the alternating group ${\bf Alt}(m)$ is simple, and its cardinality $|{\bf Alt}(m)| = m!/2$ tends to infinity as $m$ increases. As the sequence $\{n_{\ell}\}$ is increasing, for some $\ell_{0} > 0$, then $\ell \geq \ell_{0}$ and $\bA_m  = {\bf Alt}(n_{\ell})$ non-trivial implies    that $\bH_m$ has order 
$|\bH_m| = (n_{\ell})!/2 > |\whG/U|$. Thus,   the projection $\bA_m \subset \whG   \to \whG/U$ cannot be an injection, and as $\bA_m$ is a simple group, it must be contained in the kernel, so  $\bA_m \subset U$.
Let $\pi_m \colon \bA \to \bA_m$ be the projection onto the $m-th$ factor. We have that $\cD \subset \bH \subset \bA$ and let $\cD_m \subset \bA_m$ denote its image. By the choice of $m$, and that  $n_{\ell} \geq d_{\ell} = 4(|\bK_{\ell}| +2)$,   Lemma~\ref{lem-altcore} implies the subgroup $\cD_m$ has trivial core in $\bA_m$.  It follows that $\cD$ has trivial core in $U$.

 Then by Proposition~\ref{prop-lenstran}, for all $k \geq 0$, the discriminant group for the truncated group chain  $\{G_{\ell}\}_{\ell \geq k}$    is   isomorphic to $K$. In particular, $\cG_K$ is a stable group chain and is not weakly normal.

The rest of the proof proceeds as for that of Theorem~\ref{thm-finiterealization}. 
\endproof

  Note that in the above proof, we cannot assert that all leaves of the suspended foliation $\FfM$ have trivial holonomy, as examples show that    some  
   conjugate of $\cD$ in $\whG$ may intersect  $G$ non-trivially. 
     
Our final example, which is again based on the application of Theorem~\ref{thm-lubotzky2}, answers a question posed in the work \cite{DHL2016a}. In that work, the notion of a virtually regular action $(X, G , \Phi)$ with group chain $\cG = \{G_{\ell}\}_{\ell \geq 0}$ was introduced:
\begin{defn}\cite[Definition~1.12]{DHL2016a}\label{def-vitreg}
A group chain $\cG = \{G_{\ell}\}_{\ell \geq 0}$ is said to be  \emph{virtually regular} if there exists a normal subgroup $G_0' \subset G_0$ of finite index such that the restricted chain $\cG' = \{G_{\ell}'\}_{\ell \geq 0}$, where $\ds G_{\ell}' = G_{\ell} \cap G_0'$, is weakly normal in $G_0'$.
\end{defn}

There is an alternate definition of this concept, which is shown in  \cite{DHL2016a} to be equivalent: a matchbox manifold $\fM$ is \emph{virtually regular}, if there exists a homogeneous matchbox manifold $\fM'$ and a finite-to-one normal covering map $h \colon \fM' \to \fM$. Thus, the notion of virtually regular is a natural property of a matchbox manifold $\fM$, and can be checked by considering a group chain model for the holonomy action of the foliation $\FfM$.

  The following example is the first known to the authors which is not virtually regular, and gives a natural paradigm for the  construction of   group chains which are not virtually regular.

 \begin{thm}\label{thm-nvr}
  There exists a finitely-generated, residually-finite, torsion-free group $G$ with profinite completion $\whG$, such that for any 
 infinite     collection     $\{F_{\ell} \mid \ell =1,2, \ldots\}$ of  non-trivial finite simple groups, 
   their cartesian product    $\ds \bF = \prod \ F_{\ell}$  can be embedded   into  $\whG$, 
so that the resulting group chain $\cG_{\bF} = \{G_{\ell}\}_{\ell \geq 0}$ constructed as in Section~\ref{subsec-lenstra}  yields    
  an equicontinuous minimal Cantor system $(X, G , \Phi)$ whose     
     discriminant group for the       group chain $\cG_{\bF}$ is   isomorphic to $\bF$. Moreover, $\cG_{\bF}$ is not virtually regular. 
   \end{thm}
\proof
The proof follows the same approach as that used in the proof of Theorem~\ref{thm-stableCantor}.

Let $G \subset \G_n$ be the   finitely-generated, torsion-free, residually-finite   group as used in the proof of Theorem~\ref{thm-lubotzky2}, and as constructed  on  \cite[page 330]{Lubotzky1993}. Let $\whG$ denote its profinite completion.

For each $\ell \geq 0$, set $d_{\ell} = 4(|F_{\ell}| +2)$.  Then   there is an embedding of $F_{\ell}$ into the alternating groups,  $\ds F_{\ell} \subset {\bf Alt}(|F_{\ell}| +2) \subset {\bf Alt}(d_{\ell})$ as in the proof of Theorem~\ref{thm-stableCantor}. 
 Choose an increasing sequence of integers $\{n_{\ell} \mid \ell \geq 1\}$ so that $n_{\ell} \geq d_{\ell}$ for all $\ell \geq 1$. Let ${\bf Alt}(d_{\ell}) \subset {\bf Alt}(n_{\ell})$ be the embedding as the permutations on the first $d_{\ell}$ symbols.  Then we obtain an embedding $\iota_{\bF} \colon \bF \to \whG$, of the infinite product $\bF$ into $\whG$, given by the composition
   \begin{equation}\label{eq-embedallHalt2}
 \iota_{\bF} \colon   \bF \cong  \prod_{\ell \geq 1} \ F_{\ell}  \subset  \prod_{\ell \geq 1} \ {\bf Alt}(d_{\ell})  \subset  \prod_{\ell \geq 1} \ {\bf Alt}(n_{\ell})  \subset   \prod_{\ell \geq 1} \ {\bf SL}_{n_{\ell}}(\mZ_{p_{n_{\ell}}})\subset \whG \ .
\end{equation}

Set $\cD = \iota_{\bF}(\bF) \subset \whG$.
Use the method of Section~\ref{subsec-lenstra} to construct a group chain in $G$. 
 The group $G$ is residually finite, so there exists a clopen neighborhood system 
 $\{U_{\ell}  \mid \ell \geq 1\}$   about the identity in $\whG$, where each $U_{\ell}$ is normal in $\whG$.   Set $W_{\ell} = \cD \cdot U_{\ell}$ for $\ell \geq 1$, and $G_{\ell} = G \cap W_{\ell}$. Let $\cG_{\bF} = \{G_{\ell}\}_{\ell \geq 0}$ denote the resulting group chain.

Let $U \subset \whG$ be a normal clopen neighborhood of the identity. For example, given a normal subgroup $G' \subset G$ with finite index, we can take $U$ to be the profinite completion of $G'$ in $\whG$. Let $\cG_{\bF}^U = \{G_{\ell}'\}_{\ell \geq 0}$ be the group chain defined by $G_{\ell}' = G_{\ell} \cap U$ for $\ell \geq 0$.
 Then $\cD \cap U = \cap (U_{\ell} \cap U)$. 
 
 We next show that the normal core,   $core_U \cD \subset \cD$,  of $\cD \cap U$ in $U$ is a finite subgroup, and then apply Corollary~\ref{cor-anydisc} to conclude that the discriminant of the action defined by the group chain $\cG_{\bF}^U$ is a non-trivial Cantor group.  
 The following argument is similar to that used in the proof of Theorem~\ref{thm-stableCantor}, and uses that   the alternating group ${\bf Alt}(m)$ is simple for $m \geq 5$, and has order $|{\bf Alt}(m)| = m!/2$.
 Let $d_U = |\whG/U|$ be the order of the finite group. 
 
   Choose $\ell_U \geq 1$   such that $n_{\ell_U} \geq 5$ and $|{\bf Alt}(n_{\ell_U})| = (n_{\ell_U})!/2 > d_U$.
   
   The for all $\ell \geq \ell_U$,  the factor ${\bf Alt}(n_{\ell})$ in the product in \eqref{eq-embedallHalt2}  is contained in the kernel of  the projection $\whG   \to \whG/U$,   and thus, $F_{\ell} \subset {\bf Alt}(n_{\ell})   \subset U$.  Consequently, we have that 
\begin{equation}
  \cD_{\ell_U} \equiv    \prod_{\ell \geq \ell_U} \ F_{\ell}  \subset \bA_{\ell_U} \equiv  \prod_{\ell \geq \ell_U}   \ {\bf Alt}(n_{\ell})  \subset  \cD \cap  U \ .
\end{equation}
In particular, this shows that $\cD \cap  U$ contains a non-trivial Cantor group. Moreover,  by applying  Lemma~\ref{lem-altcore} to each factor of  the product in $\cD_{\ell_U}$, we see that $\cD_{\ell_U}$ has   trivial core in $U_{\ell}$ as well. Thus, we have
\begin{equation}
core_U \cD  \subset  \prod_{1 \leq \ell < \ell_U} \ F_{\ell} \ ,
\end{equation}
and so is a finite normal subgroup of $U$. 

By Corollary~\ref{cor-anydisc}, the quotient group chain $\{(G_{\ell} \cap U)/(core_U \cD)\}_{\ell \geq 0}$ has non-trivial discriminant group 
$\ds \cD/(core_U \cD)$ which contains a subgroup isomorphic to the non-trivial Cantor group   $\cD_{\ell_U}$. 
For $\ell > 0$, apply this to the case $U = U_{\ell}$ to obtain that the quotient chain
$\{(G_{\ell} \cap U_{\ell})/(core_{U_{\ell}} \cD)\}_{\ell \geq 0}$   is not equivalent to a normal chain.
Now suppose that the restricted group chain $\cG_{\bF}^{U_{\ell}}$ is  equivalent to a    normal chain, then as $core_{U_{\ell}} \cD$ is a normal subgroup of $U_{\ell}$, this implies that the 
quotient group chain $\{(G_{\ell} \cap U_{\ell})/(core_{U_{\ell}} \cD)\}_{\ell \geq 0}$ is equivalent to a normal chain, hence has trivial discriminant, which is a contradiction. Thus, the group chain   $\cG_{\bF}$ is not virtually regular. 
\endproof

   \subsection{Open problems}\label{subsec-open}
  
     There are many variations of the above method that can be considered, and open questions about the resulting minimal Cantor actions.
  First, it is interesting to understand the answer to the following.
   
   \begin{prob}\label{prob-lubotzky}
   Given    a separable profinite group $\whH$ and an embedding into a profinite group $\whG$ with trivial rational core, constructed using the methods of \cite{Lubotzky1993}, give criteria for when  the resulting equicontinuous minimal Cantor system $(X, G , \Phi)$ is weakly normal, and whether the action is stable or wild. Furthermore, when do the resulting actions satisfy the \emph{SQA} condition of Section~\ref{subsec-defSQA}?
   \end{prob}  
    There is also an extensive literature for the construction of   embeddings of groups $H$ into the profinite completions of torsion-free, finitely generated nilpotent and solvable groups. For example,  the work  \cite{CBKL1988} shows that if $G$ is a finitely-generated, torsion-free nilpotent group, then the profinite completion $\whG$ is torsion free, so if $D \subset \whG$ is a closed subgroup, then it must be a Cantor group. 
    
    On the other hand, the works \cite{Evans1990,KrophollerWilson1993} shows that there exists a countable, torsion-free, residually finite, metabelian group $G$ such that its profinite completion contains a non-trivial torsion subgroup. The work \cite{Quick2001} studies the profinite topology of nilpotent groups of class two and finitely generated centre-by-metabelian groups, and uses this to construct embeddings of finite groups into the profinite completions of these classes of groups. However, the embedding obtained in \cite{Quick2001} is contained in the center of $G$, so does not satisfy the trivial core condition. We conclude with an open question, suggested by the examples and results of the works \cite{Dyer2015,DHL2016a,DHL2016b,FO2002,RT1971b,Schori1966}.
  \begin{prob}\label{prob-amenable}
Determine which groups $H$ can be embedded as a closed subgroup of $\whG$ with trivial rational core, where $G$ is a finitely-generated, torsion-free amenable group.   
\end{prob}


\end{document}